\documentclass[a4paper,11pt]{article}

\usepackage{ucs}
\usepackage[utf8]{inputenc}

\usepackage{graphicx}
\usepackage{amsfonts}
\usepackage{dsfont}
\usepackage{amssymb}
\usepackage{amsmath}
\usepackage{amsthm}
\usepackage{enumerate}
\usepackage{stmaryrd}
\usepackage{fullpage}
\usepackage{ifthen}
\usepackage{subfigure}
\usepackage{epic}
\usepackage{authblk}
\usepackage{textcomp}
\usepackage{mathrsfs}
\usepackage{bm}
\usepackage[small]{caption}
\usepackage{tikz,tkz-tab}
\usetikzlibrary{matrix}
\usepackage{pgfplots}
\pgfplotsset{compat=newest} 
\pgfplotsset{plot coordinates/math parser=false}

\usepackage[hypertexnames=false,colorlinks=true,linkcolor=blue,citecolor=blue]{hyperref}
\usepackage[numbers,comma,square,sort&compress]{natbib}
\usepackage[letterpaper,text={6.5in,9in},centering]{geometry}

\usepackage{color}
\usepackage{titlesec}

\graphicspath{{eps/}{pdf/}{png/}}
\newcommand{\bqq}{\begin{equation}}
\newcommand{\eqq}{\end{equation}}
\newcommand{\bqs}{\begin{equation*}}
\newcommand{\eqs}{\end{equation*}}

\newcommand{\R}{\mathbb{R}} 
\newcommand{\N}{\mathbb{N}}

\newcommand{\bu}{\mathbf{u}}

\newcommand{\bI}{\mathbf{I}}

\newcommand{\E}{\mathcal{E}}
\newcommand{\V}{\mathcal{V}}

\newcommand{\cG}{\mathcal{G}}

\newcommand{\md}{\mathrm{d}}

\newtheorem{lem}{Lemma}[section]
\newtheorem{thm}{Theorem}

\newtheorem{rmk}[lem]{Remark}

\newtheorem{hyp}[lem]{Hypothesis}

\newenvironment{Proof}[1][.]%
 {\begin{trivlist}\item[]\textbf{Proof#1 }}%
 {\hspace*{\fill}$\rule{0.3\baselineskip}{0.35\baselineskip}$\end{trivlist}}

  {\begin{trivlist}\item[]{\bf Hypothesis #1 }\em}{\end{trivlist}}

\numberwithin{equation}{section}

\title{Dynamics of epidemic spreading on connected graphs}

\author[1]{Christophe Besse \& Gr\'egory Faye}
\affil[1]{\small CNRS, UMR 5219, Institut de Math\'ematiques de Toulouse, 31062 Toulouse Cedex, France}

\begin{document}
\maketitle

\begin{abstract}
We propose a new model that describes the dynamics of epidemic spreading on connected graphs. Our model consists in a PDE-ODE system where at each vertex of the graph we have a standard SIR model and connexions between vertices are given by heat equations on the edges supplemented with Robin like boundary conditions at the vertices modeling exchanges between incident edges and the associated vertex. We describe the main properties of the system, and also derive the final total population of infected individuals. We present a semi-implicit in time numerical scheme based on finite differences in space which preserves the main properties of the continuous model such as  the uniqueness and positivity of solutions and the conservation of the total population. We also illustrate our results with a selection of numerical simulations for a selection of connected graphs.
\end{abstract}

\bigskip


\bigskip

\section{Introduction}

Classical SIR compartment models are cornerstone models of epidemiology which
allow one to study the evolution of an infected population at a given spatial
scale (e.g. whole countries, regions, counties or cities). Such models date back
to the pioneer work of Kermack and McKendrick \cite{KMK27} and describe the
evolution of susceptible (S) and infected (I) populations which eventually
become removed (R) via systems of ordinary differential equations which typically take the form
\bqq
\left\{
\begin{split}
S'(t)&=-\tau S(t)I(t),\\
I'(t)&=\tau S(t)I(t)-\eta I(t),\\
R'(t)&=\eta I(t),
\end{split}
\right.
\label{SIR}
\eqq
where $\tau>0$ is a contact rate between susceptible and infected populations, and $1/\eta>0$ is the average infectious period; see \cite{Heth00} for a review on SIR models. These models have been used in the past to reproduce data of epidemic outbreaks such as the bubonic plague \cite{KMK27}, malaria \cite{MSS11}, SARS influenza \cite{CDC03,MSB16} and most recently COVID-19 \cite{NJE20,MW18,LMSW20}; see also \cite{MW18} for other applications.

\begin{figure}[t!]
  \centering
  \includegraphics[width=.4\textwidth]{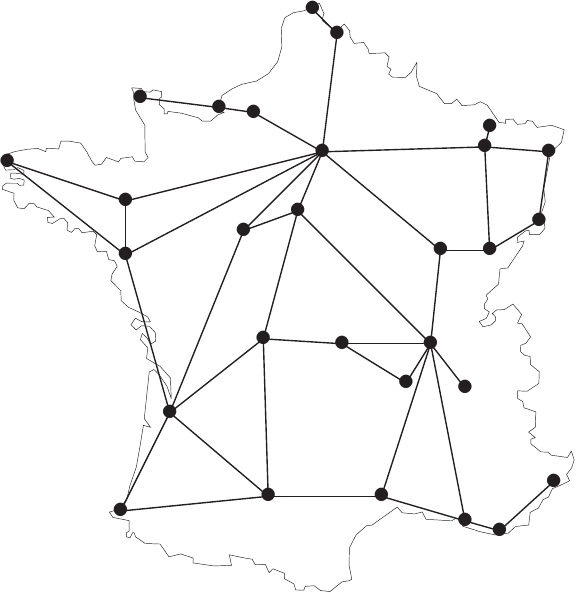}
  
  \caption{Map of France with an illustration of a connected graph connecting major cities.}
  \label{fig:map}
\end{figure}

In classical SIR models such as \eqref{SIR}, interactions among the infected population are oversimplified, and when taken into account they typically involve transfer matrices of populations of infected between various uniform patches \cite{VanW02,MSB16,MSB18}. Our interest lies in the understanding of the intricate interplay between spatial effects and heterogeneous interactions among infected populations. Schematically, we propose a model composed of cities linked by a given transportation network (roads, railroads or rivers), see Figure~\ref{fig:map} for an illustration in the case of France. It will turn out that the appropriate theoretical framework will be graph theory where each vertices of the graph will be thought as the cities and the edges the lines of transportation. In a first approximation, we will assume that infected populations are only subject to spatial diffusion along the lines, as it is traditionally assumed in classical spatial SIR models \cite{A77,D78,Reluga,BRR20}. As a consequence, in our model, the dynamics of the epidemic only takes place in the cities. Interactions are then modeled by flux exchanges between cities and lines where we assume that some fraction of infected individuals can either leave a city to be on a line, or leave a line and stop in a city, or pass from one line to another through a city. The typical question that we address here can easily be stated as follows. Given a connected graph of cities linked by roads and an initial configuration of infected individuals, how does the epidemic spread into the network and what is the eventual final configuration of the infected population? Our aim here is to gain insight into this spreading aspect at the fundamental mathematical level of a SIR type model that incorporates the possibility of infected individuals to travel along a specific given transportation network.

Our framework is at the crossroad of models that take into account lines of transportation such as recent reaction-diffusion models that study propagation of epidemics along lines with fast diffusion \cite{BRR20} and models that incorporate networks with more sophisticated interactions dynamics \cite{BDL08,SB19,BB20,BG18}. On a formal level, our proposed model can be thought of as being a one-dimensional version of the planar reaction-diffusion system of \cite{BRR20} with a line of fast diffusion in the case of one city and one line of transportation. Actually, the graph structure of the transportation network provides a natural embedding into a planar spatial domain. From a mathematical point of view, our model shares also some similarities with the PDE-ODE model of \cite{DIWB20} which studies the spread of airborne diseases where the movement of pathogens in the air is assumed to follow a linear diffusion.

\section{Model formulation and main results}

Throughout, we denote by $\cG=(\V,\E)$ a compact metric graph, \textit{i.e.} a collection of vertices $\V$ and edges $\E$ and further assume that $\cG$ is finite and connected. Each edge $e\in\E$ is identified with a segment $\Omega_e=[0,\ell_e]$ with $\ell_e\in(0,\infty)$, where $\ell_e$ is the finite length of the edge. A real valued function $u:\cG\longrightarrow\R$ is a collection of one dimensional maps defined for each edge $e\in\E$:
\bqs
u_e:\Omega_e\longrightarrow \R.
\eqs
For future references, we define $BC(\cG,\R)$ the space of bounded continuous functions on $\cG$ as
\bqs
BC(\cG,\R):=\underset{e\in\E}{\bigoplus}~BC(\Omega_e,\R),
\eqs
and similarly $BC^k(\cG,\R)$ with $k\geq1$. We define the $L^\infty$ norm on
$\mathcal{G}$ for $u \in BC(\cG,\R)$ as
\[
  \|u\|_\infty:=\max_{e\in \mathcal{E}} \sup_{x\in \Omega_e} |u_e(x)|.
\]

\subsection{A SIR model on compact connected graph}

Given a graph $\cG$, we let $X_v(t):=(S_v(t),I_v(t),R_v(t))\in\R^3$, for each $v\in\V$, where $S_v(t)$ represents the population of susceptible individuals, $I_v(t)$ the population of infected individuals and $R_v(t)$ the population of susceptible individuals at vertex $v\in\V$ and time $t>0$. We assume that $X_v$ evolves according to a SIR model of the form
\bqq
\left\{
\begin{split}
S'_v(t)&=-\tau_v S_v(t)I_v(t),\\
I_v'(t)&=\tau_v S_v(t)I_v(t)-\eta_v I_v(t)+\sum_{e\sim v} \alpha^v_e u_e(t,v)-\overline{\lambda}_vI_v(t),\\
R_v'(t)&=\eta_v I_v(t),
\end{split}
\right.
\label{SIRgraph}
\eqq
where $\tau_v,\eta_v>0$ are the intrinsic parameters of the epidemic which may depend on the vertex $v$. The contribution $-\overline{\lambda}_vI_v(t)$ in the right-hand side of the equation for the infected population traduces the fact that infected individuals can leave the vertex $v$ to incident edges whereas $\sum_{e\sim v} \alpha^v_e u_e(t,v)$ reflects the contribution of incoming infected individuals from incident edges. Here, $e\sim v$ denotes the edges incident to the vertex $v$ and
\bqs
\overline{\lambda}_v:=\sum_{e\sim v}\lambda_{e}^v, % \in (0,1),
\eqs
such that $\lambda_{e}^vI_v(t)$ infected individuals leave vertex $v$ to edge
$e$. We have assumed that only the infected population is subject to movement,
and we think of $S_v$ being an ambiant population whose movement does not affect
its distribution. We recover the standard SIR model \eqref{SIR} by considering
the trivial graph $\mathcal{G}=(\{v\},\emptyset)$. Throughout the manuscript, we will assume the following standing hypothesis on the coefficients $\alpha_e^v$ and $\lambda_e^v$ in \eqref{SIRgraph}.

\begin{hyp}\label{hypPos} For each $(e,v)\in \E\times\V$ we assume that
\bqs
\alpha_e^v \in(0,1) \quad \text{ and }\quad  \lambda_e^v\in(0,1),
\eqs
together with
\bqs
\sum_{e\sim v}\lambda_{e}^v \in (0,1) \quad \text{ and } \quad \sum_{e\sim v}\alpha_{e}^v \in (0,1).
\eqs
\end{hyp}

Next, for each $e\in\E$, we let $d_e>0$ and we assume that $u_e$ evolves according to
\bqq
\partial_t u_e(t,x) = d_e \partial_x^2 u_e(t,x),\quad t>0, \quad x\in\overset{\circ}{\Omega}_e.
\label{heatgraph}
\eqq
Assuming that infected individuals have local diffusion along the edges of the graph is a first approximation, and this can be viewed as a limiting Brownian movement of individuals. We shall come back to this modeling hypothesis later in the manuscript, but possible extensions could be to incorporate nonlocal diffusion or transport terms.

It now remains to model the exchanges of infected individuals at the vertices. Fo each $v\in\V$, we associate an integer $\delta_v \geq 1$ which we refer to as its degree (\textit{i.e.} number of edges incident to the vertex $v$). We define $\bu_v(t)\in\R^{\delta_v}$ as the column vector function
\bqs
\bu_v(t):=(u_e(t,v))_{e\sim v},
\eqs
where we recall that $e\sim v$ denotes the edges incident to the vertex $v$, and thus $u_e(t,v)$ is the corresponding limit value of $u_e$ at $x=v$. Define also $\partial_n \bu_v(t)\in\R^{\delta_v}$ as the column vector function
\bqs
\partial_n \bu_v(t):=(\partial_n u_e(t,v))_{e\sim v},
\eqs
where $\partial_n u_e(t,v)$ is the outwardly normal derivative of $u_e$ at the vertex $v$. Our boundary conditions at the vertex $v$ that link \eqref{SIRgraph} and \eqref{heatgraph} are described by
\bqq
D_v \partial_n \bu_v(t)+ K_v \bu_v(t) = \Lambda_v  \bI_v(t),
\label{BCgraph}
\eqq
where $D_v \in \mathscr{M}_{\delta_v}(\R)$ is the diagonal matrix $D_v=\mathrm{diag}[(d_e)_{e\sim v}]$ and $K_v \in \mathscr{M}_{\delta_v}(\R)$ whose structure will be specified below. Formally, \eqref{BCgraph} traduces the  balance of fluxes of infected individuals at the vertex $v$, and we will demonstrate this heuristic rigorously by showing in the forthcoming Subsection~\ref{subsec_cons} the conservation of total population.

\subsection{Assumptions on the connectivity matrices $K_v$}

\begin{figure}[ht!]
  \centering
  \begin{tikzpicture}[scale=2]
    \node at (0,0) {$\bullet$};
    \node[above right] at (0,0) {$v$};
    \node[right] at (-1,2) {$e_1$};
    \node[right] at (-1,-2) {$e_2$};
    \node[above right] at (2,0) {$e_3$};
    \draw (0,0)--(2,0);
    \draw (0,0)--(-1,2);
    \draw (0,0)--(-1,-2);
    \draw[->,>=latex] (1.3,0.3) to[bend right] (-0.3,1.3);
    \draw[->,>=latex] (-0.4,1.5) to[bend left] (1.5,0.3);
    \node at (0.6,0.8) {$\nu_{e_3,e_1}^v$};
    \node at (1,1.3) {$\nu_{e_1,e_3}^v$};

    \draw[->,>=latex] (-1  ,1.2) to[bend right] (-1,-1.2);
    \draw[->,>=latex] (-0.85,-1.05) to[bend left] (-0.85,1.05);
    \node at (-0.85,0) {$\nu_{e_2,e_1}^v$};
    \node at (-1.6,0) {$\nu_{e_1,e_2}^v$};

    \draw[->,>=latex] (1.3,-0.3) to[bend left] (-0.3,-1.3);
    \draw[->,>=latex] (-0.4,-1.5) to[bend right] (1.5,-0.3);
    \node at (0.4,-0.9) {$\nu_{e_3,e_2}^v$};
    \node at (1,-1.3) {$\nu_{e_2,e_3}^v$};

    \draw[->,>=latex] (0.2,0.1)--(1.,0.1);
    \node[above] at (0.6,0.1) {$\lambda_{e_3}^v$};
    \draw[->,>=latex] (1.,-0.1)--(0.2,-0.1);
    \node[below] at (0.6,-0.1) {$\alpha_{e_3}^v$};

    \draw[->,>=latex] (0.,0.2)--(-0.3,0.8);
    \node[right] at (-0.15,0.6) {$\lambda_{e_1}^v$};
    \draw[<-,>=latex] (-0.15,0.1)--(-0.45,0.7);
    \node[left] at (-0.3,0.4) {$\alpha_{e_1}^v$};

    \draw[<-,>=latex] (0.,-0.2)--(-0.3,-0.8);
    \node[right] at (-0.15,-0.6) {$\alpha_{e_2}^v$};
    \draw[->,>=latex] (-0.15,-0.1)--(-0.45,-0.7);
    \node[left] at (-0.3,-0.4) {$\lambda_{e_2}^v$};
    
  \end{tikzpicture}  
  \caption{Schematic illustration of the exchanges at a given vertex $v$ with $\delta_v=3$.}
  \label{fig:exchange}
\end{figure}
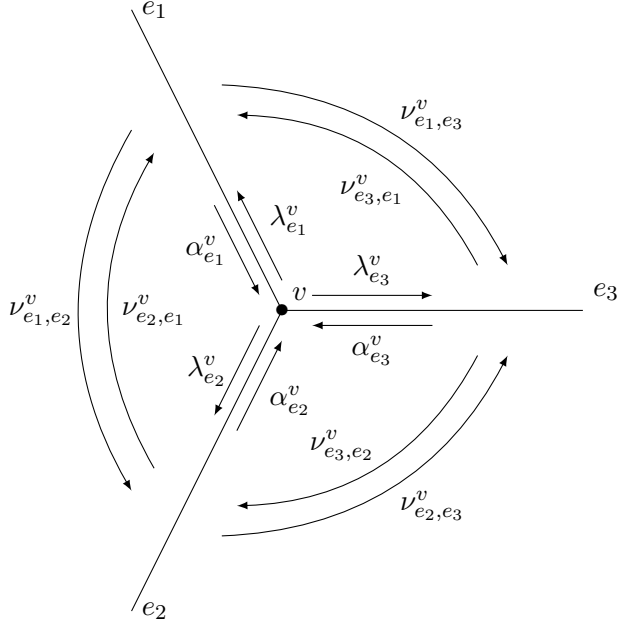

We now precise the form of the matrix $K_v$ entering in the boundary condition \eqref{BCgraph}. Essentially, $K_v$ gathers two contributions. One contribution comes from the exchanges between infected individuals at the vertex with the incoming infected individuals for the incident edges. The second contribution traduces exchanges between edges. Indeed we allow infected individuals to pass from one edge to another one. More precisely, we have that $K_v$ splits into two parts
\bqs
K_v :=A_v + N_v,
\eqs
where the matrix $A_v \in \mathscr{M}_{\delta_v}(\R)$ is the diagonal matrix $A_v=\mathrm{diag}(\alpha_{e}^v)_{e\sim v}$ while the matrix $N_v \in \mathscr{M}_{\delta_v}(\R)$ is such that the sum of each column is zero. More precisely, if we label by $e\sim v=(e_1,\cdots,e_{\delta_v})$ the edges incident to the vertex $v$, we have that for all $i=1,\cdots,\delta_v$
\bqs
(N_v)_{i,i} = \sum_{j \neq i} \nu^v_{e_i,e_j}
\eqs
and for $i\neq j=1,\cdots,\delta_v$
\bqs
(N_v)_{i,j}=-\nu^v_{e_j,e_i}.
\eqs
In the case $\delta_v=3$, we get
\bqs
N_v = \left(
\begin{matrix}
\nu^v_{e_1,e_2}+\nu^v_{e_1,e_3} & -\nu^v_{e_2,e_1}  & -\nu^v_{e_3,e_1}\\ 
-\nu^v_{e_1,e_2} & \nu^v_{e_2,e_1}+\nu^v_{e_2,e_3} & -\nu^v_{e_3,e_2}\\
-\nu^v_{e_1,e_3} & -\nu^v_{e_2,e_3} & \nu^v_{e_3,e_1}+\nu^v_{e_3,e_2}
\end{matrix}
\right),
\eqs
see Figure \ref{fig:exchange} for an illustration in that case.

Furthermore, for the diagonal term we will use the shorthand notation
\bqs
(N_v)_{e,e} = \sum_{e'\neq e} \nu^v_{e,e'}.
\eqs
The fact that $N_v \in \mathscr{M}_{\delta_v}(\R)$ is such that the sum of each column is zero precisely traduces the fact that there is the conservation of infected individuals through exchanges between incident edges at each vertex. And, we remark that it implies that the matrix $K_v$ has a strict column diagonal dominance in the sense that for each $i=1,\cdots,\delta_v$
\bqs
\sum_{j=1}^{\delta_v} (K_v)_{e_i,e_j}=\alpha_{e_i}^v>0,
\eqs
because of this specific structure of $N_v$. From now on we also assume that $K_v$ has a diagonal dominance with respect to its lines. This property will be crucial later on in the proof of existence of solutions. As a consequence, we impose the following running assumptions on the matrices $K_v$.

\begin{hyp}\label{hypDD}
For each $v\in\V$ and $(e,e')\in\E\times\E$, we assume that
\bqs
\nu_{e,e'}^v\in[0,1).
\eqs
Furthermore, we impose that for all $v\in\V$
\bqs
\sum_{e'\neq e} \nu^v_{e',e} < \alpha_e^v+ \sum_{e'\neq e} \nu^v_{e,e'},
\eqs
together with
\bqs
(K_v)_{e,e} := \alpha_e^v+ \sum_{e'\neq e} \nu^v_{e,e'} \in(0,1),
\eqs
for each $e\sim v$.
\end{hyp}

\begin{rmk}
If the exchanges between the edges are symmetric, that is for each $v\in\V$ the matrices $N_v$ are symmetric, that is
\bqs
 \nu^v_{e,e'}= \nu^v_{e',e}, \quad \forall (e,e')\in\E\times\E,
\eqs
then Hypothesis~\ref{hypDD} is automatically satisfied.
\end{rmk}

\subsection{Initial configuration}

Finally, we complement our coupled PDE-ODE \eqref{SIRgraph}-\eqref{heatgraph}-\eqref{BCgraph} with some initial conditions. We assume that at $t=0$, we have
\bqs
u(t=0,\cdot) = u^0 \in BC(\cG,\R),
\eqs 
such that for $e\in\E$,
\bqs
u^0_e(x)\geq 0, \quad x\in \Omega_e.
\eqs
On the other hand, for the ODE system \eqref{SIRgraph}, we suppose that
\bqs
(S_v(t=0),I_v(t=0),R_v(t=0))=(S_v^0,I_v^0,0) \geq \mathbf{0}, \quad \forall v\in\V.
\eqs
We further assume that \eqref{BCgraph} is satisfied at $t=0$
\bqs
D_v \partial_n \bu_v^0+ K_v \bu_v^0 = \Lambda_v  \bI_v^0,
\eqs
with obvious notations $\bu_v^0:=(u_e^0(v))_{e\sim v}$ and $\bI_v^0:=( I_v^0,\cdots,I_v^0)^{\mathbf{t}}$. Last, we impose that the initial total population of infected individuals is strictly positive, 
\bqs
\sum_{v\in\V} I_v^0>0,
\eqs
and that susceptible individuals are initially present at each vertex of the graph
\bqs
S_v^0>0, \quad \forall v\in\V.
\eqs
This in turn implies that the total population is initially 
\bqs
M^0:=\sum_{e\in\E} \int_{\Omega_e} u_e^0(x)\md x+\sum_{v\in\V}\left(S_v^0+I_v^0\right)>0.
\eqs

\subsection{Conservation of total population}\label{subsec_cons}

Assuming that there is a solution to $(u,(X_v)_{v\in\V})$ to \eqref{SIRgraph}-\eqref{heatgraph}-\eqref{BCgraph}, we have that the total mass of the system $M(t)$ defined as
\bqs
M(t):=\sum_{e\in\E} \int_{\Omega_e} u_e(t,x)\md x+\sum_{v\in\V}\left(S_v(t)+I_v(t)+R_v(t)\right)
\eqs
is  a conserved quantity and thus independent of $t$.

To see that, we first remark that
\bqs
\left(S'_v(t)+I'_v(t)+R'_v(t)\right) =  \sum_{e\sim v} \alpha^v_e u_e(t,v)-  \overline{\lambda}_vI_v(t)= \langle A_v  \bu_v(t) ,\mathbf{1}_{\delta_v}\rangle -\langle \Lambda_v  \bI_v(t) ,\mathbf{1}_{\delta_v}\rangle
\eqs
with 
\bqs\mathbf{1}_{\delta_v}:=(1,\cdots,1)^{\mathbf{t}}\in\R^{\delta_v},
\eqs
and $\langle \cdot,\cdot \rangle$ is the standard Euclidean inner product on $\mathbb{R}^{\delta_v}$.
On the other hand let us define
\bqs
m(t):= \sum_{e\in\E} \int_{\Omega_e} u_e(t,x)\md x,
\eqs
and assume that $u$ is a classical solution of $\eqref{heatgraph}$, which we will prove in the next section, and compute
\begin{align*}
m'(t)&=\sum_{e\in\E} \int_{\Omega_e} \partial_t u_e(t,x)\md x=\sum_{e\in\E} d_e \left[ \partial_x u_e(t,x) \right]_{\partial \Omega_e} \\
&=\sum_{v\in\V} \langle D_v \partial_n \bu_v(t),\mathbf{1}_{\delta_v} \rangle \\
&=\sum_{v\in\V} \langle  \Lambda_v  \bI_v(t) - K_v \bu_v(t),\mathbf{1}_{\delta_v} \rangle \\
&=\sum_{v\in\V} \langle  \Lambda_v  \bI_v(t) - A_v \bu_v(t),\mathbf{1}_{\delta_v} \rangle - \underbrace{\sum_{v\in\V} \langle N_v \bu_v(t),\mathbf{1}_{\delta_v} \rangle}_{=0}\\
&=-\sum_{v\in\V}\left(S'_v(t)+I'_v(t)+R'_v(t)\right).
\end{align*}
The fact that
\bqs
\sum_{v\in\V} \langle N_v \bu_v(t),\mathbf{1}_{\delta_v} \rangle=0
\eqs
is a direct consequence on the specific structure of each matrix $N_v$ and the
fact that the sum of each column is zero. We therefore conclude that $M'(t)=0$ and
\[
  \sum_{e\in\E} \int_{\Omega_e} u_e(t,x)\md x+\sum_{v\in\V}\left(S_v(t)+I_v(t)+R_v(t)\right)=M^0, \quad \forall t\geq 0.
\]

\paragraph{Biological interpretation.} Our model is thus consistent with the conservation of the total population as it is traditionally the case for SIR model in the case of zero natality/mortality rate. The exchanges between the vertices and the edges exactly compensate each other as is natural.

\subsection{Main results and outline of the paper}

We now present our main results regarding our model \eqref{SIRgraph}-\eqref{heatgraph}-\eqref{BCgraph}. At this stage of the presentation, we remain formal and refer to the following sections for precise statements and assumptions.

\paragraph{Main result 1: Existence and uniqueness of classical solutions.} We prove in Theorem~\ref{thmexun} below that for each well prepared initial condition our model \eqref{SIRgraph}-\eqref{heatgraph}-\eqref{BCgraph} admits a unique positive classical solution which is global in time. We remark that the system \eqref{SIRgraph}-\eqref{heatgraph}-\eqref{BCgraph} is not standard as it couples a system of PDEs to ODEs at each vertices through inhomogeneous Robin boundary conditions. As a consequence, the existence and uniqueness of classical solutions has to be proved. This analysis is conducted in Section~\ref{secExUn}.

\paragraph{Main result 2: Long time behavior of the solutions.} We fully characterize the long time behavior of the unique solution of our model. More precisely, we prove that the final total population of infected individuals at each vertex, denoted by $\mathcal{I}_v^\infty$, is a well defined quantity: $0<\mathcal{I}_v^\infty<\infty$ for $v\in\V$ and $(\mathcal{I}_v^\infty)_{v\in\V}$ are solutions of a system of $c_\V+1$ implicit equations, where $c_\V$ stands for the cardinal of $\V$, which belong to the parametrized submanifold
\bqs
\sum_{v\in\V}\left(S_v^0e^{-\tau_v \mathcal{I}_v^\infty}+\eta_v \mathcal{I}_v^\infty \right)=M^0,
\eqs
where $M^0$ is the initial total mass. We refer to Theorem~\ref{thmFTI} for a precise statement. We also present further qualitative results on the final total configuration $(\mathcal{I}_v^\infty)_{v\in\V}$ in the fully symmetric case where we obtain closed form formula (see Lemma~\ref{lem_sym}) and in the case of two vertices where we manage to obtain sharp bounds on the final total populations of infected individuals (see Lemma~\ref{lem_2vert}). In each case, we manage to relate these quantities to standard basic and effective reproductive number for classical SIR model. The aforementioned results are proved in Section~\ref{secLong}.

\paragraph{Main result 3: A mass preserving semi-implicit numerical scheme.} We propose and analyze a semi-implicit in time numerical scheme based on finite differences in space which has the property to preserve a {\it discrete} total mass associated to the discretization. We prove that if the time discretization constant is smaller than a universal constant depending only on the parameters of the system (and not on the space discretization constant) and if $N_v$ is symmetric for each $v\in\V$, then our mass preserving semi-implicit numerical scheme is well-posed and preserves the positivity of the solutions. We refer to Section~\ref{secNumSc} for a presentation of the numerical scheme and Theorem~\ref{thmNumSc} for a precise statement of our main result.

\paragraph{Main result 4: Numerical results for various types of graphs.} We illustrate our theoretical findings with selection of numerical simulations for various types of graphs in Section~\ref{secSim}. We respectively study the case of 2 vertices and 1 edge, 3 vertices and 3 edges (closed graph), 4 vertices and 3 edges (star-shape graph) and $N+1$ vertices and $N$ edges with $N$ being arbitrarily large (lattice graph). Most notably, in the last case, we show the propagation of the epidemics across the vertices of the graph in the form of a traveling wave.

\section{The Cauchy problem: existence and uniqueness of classical solutions}\label{secExUn}

This section is devoted to the proof of the following main theorem which guarantees that our model is well-posed.

\begin{thm}\label{thmexun}
For each $(S_v^0,I_v^0)\geq \mathbf{0}$ with $S_v^0>0$, $\sum_{v\in\V} I_v^0>0$ and $u^0\in BC(\cG,\R^+)$ that satisfy the boundary condition \eqref{BCgraph}, there exists a unique positive global solution $(S_v,I_v,R_v)\in \mathscr{C}^1(\R^+,\R^+\times\R^+\times\R^+)$ and $u \in \mathscr{C}^{1,2}(\R^+_*\times \cG,\R^+)$.
\end{thm}

The proof of Theorem~\ref{thmexun} is divided into two parts. We first prove the existence of positive global classical solutions and then show that such constructed solutions are unique. We look for solutions that satisfy \eqref{SIRgraph}-\eqref{heatgraph}-\eqref{BCgraph} in the classical sense, and we always assume that $(S_v^0,I_v^0)\geq \mathbf{0}$ with $S_v^0>0$, $\sum_{v\in\V} I_v^0>0$ and $u^0\in BC(\cG,\R^+)$, that is for all $e\in\E$, $u_e^0\geq 0$ is bounded continuous on $\Omega_e$. We further assume that the initial conditions satisfy the boundary condition \eqref{BCgraph}. We remark that the system \eqref{SIRgraph}-\eqref{heatgraph}-\eqref{BCgraph} is not standard as it couples a system PDEs to ODEs at each vertices through inhomogeneous Robin boundary conditions. As a consequence, the well-posedness of the Cauchy problem has to be proved.

\begin{rmk}
Our existence and uniqueness result extends trivially in the case that parameters $\tau_v>0$, $\alpha_e^v\in(0,1)$, $\lambda_e^v\in(0,1)$ and $\nu_{e,e'}^v\in[0,1)$ are continuous functions of time satisfying $\tau_v(t)>0$, $\alpha_e^v(t)\in(0,1)$, $\lambda_e^v(t)\in(0,1)$ and $\nu_{e,e'}^v(t)\in[0,1)$ together with Hypotheses~\ref{hypPos}-\ref{hypDD} verified at all times $t>0$.
\end{rmk}

\subsection{Existence} In this section, we construct a classical solution to \eqref{SIRgraph}-\eqref{heatgraph}-\eqref{BCgraph} through a limiting argument. We will obtain a solution $(u,(X_v)_{v\in\V})$ has the limit of a subsequence of solution $((u^n,(X_v^n)_{v\in\V}))_{n\geq0}$ of the following problems
\bqq
\left\{
\begin{split}
\frac{\md S^n_v(t)}{\md t}&=-\tau_v S_v^n(t)I_v^n(t),\\
\frac{\md I^n_v(t)}{\md t}&=\tau_v S_v^n(t)I_v^n(t)-(\eta_v+\overline{\lambda}_v) I_v^n(t)+\sum_{e\sim v} \alpha^v_e u_e^{n-1}(t,v),\\
\frac{\md R^n_v(t)}{\md t}&=\eta_v I_v^n(t),
\end{split}
\right. \quad t>0, \quad \forall v\in\V,
\label{SIRgraph_n}
\eqq
with
\bqq
D_v \partial_n \bu_v^n(t)+ K_v \bu_v^n(t) = \Lambda_v  \bI_v^n(t), \quad t>0,  \quad \forall v\in\V,
\label{BCgraph_n}
\eqq
and
\bqq
\partial_t u_e^n(t,x) = d_e \partial_x^2 u_e^n(t,x), \quad t>0, \quad x\in\overset{\circ}{\Omega}_e, \quad \forall e\in\E.
\label{heatgraph_n}
\eqq
starting from $u^0\in BC(\cG,\R^+)$ and $(X^0_v)_{v\in\V}$. Note that \eqref{SIRgraph_n}-\eqref{BCgraph_n}\eqref{heatgraph_n} is supplemented by the same initial condition $(u^0,(X_v^0)_{v\in\V})$ at each step. We proceed along three main steps.

\paragraph{Step \#1: solvability of \eqref{SIRgraph_n}-\eqref{BCgraph_n}-\eqref{heatgraph_n}.} We first show that \eqref{SIRgraph_n}-\eqref{BCgraph_n}-\eqref{heatgraph_n} admits a unique solution. It can be done by induction. Assume that at step $n-1$, we have constructed a solution $(u^{n-1},(X_v^{n-1})_{v\in\V})$ such that for each $t\mapsto u_e^{n-1}(t,v)$ is continuous, then we get the existence of a unique solution of \eqref{SIRgraph_n} which is $\mathscr{C}^1$ in time. Next we solve the system of PDEs \eqref{heatgraph_n}-\eqref{BCgraph_n} whose coupling comes from the boundary conditions and owing that now the right-hand side of \eqref{BCgraph_n} can be seen as given inhomogeneous term of class $\mathscr{C}^1$ in time. As both $D_v$ and $K_v$ are invertible matrices, we get the existence of a classical solution $u^n \in \mathscr{C}^{1,2}$ which then ensures that $t\mapsto u_e^{n}(t,v)$ is continuous.

\paragraph{Step \#2: a priori estimates.} Let $0<T < 1$ be fixed. We first show by a recursive argument that $0< S_{v}^n$, $0\leq I_{v}^n$, $0\leq R_{v}^n$ for each $v\in\V$  and $0\leq u^n_e$ for each $e\in\E$. It is trivial at $n=0$. Let assume that is it true at $n-1$. We start from \eqref{SIRgraph_n} and a direct integration gives
\begin{align*}
S_v^n(t)&=S_v^0 e^{-\tau_v \int_0^t I_v^n(s)\md s}>0,\\
I_v^n(t)&=I_v^0e^{-(\eta_v+\overline{\lambda}_v)t+\int_0^t S_v^n(s)\md s}+\sum_{e\sim v} \alpha^v_e \int_0^t e^{-(\eta_v+\overline{\lambda}_v)(t-s)+\int_s^t S_v^n(\tau)\md\tau} u_e^{n-1}(s,v)\md s \geq 0,\\
R_v^n(t)&=\eta_v \int_0^t I_v^n(s)\md s \geq 0.
\end{align*}
Now owing that $0\leq I_{v}^n$ for each $v\in\V$, the maximum principle implies that $u_e^n\geq 0$ for each $e\in\E$. Assume by contradiction that $e_*\in\E$ is the component which reaches a negative minimum, namely $u_{e_*}^n(t_*,x_*)=-\delta<0$ with $u_{e_*}^n(t,x)>-\delta$ for $t<t_*$ and $x\in \Omega_{e_*}$ and for each $e\neq e_*$ we have $u_e^n(t,x)>-\delta$ for $t\leq t_*$ and $x\in \Omega_{e}$. We know that $x_* \in \partial \Omega_{e_*}$ and let denote $v_*=x_*\in\V$ the vertex where this occurs. The Hopf lemma implies that $\partial_n u_{e_*}^n(t_*,v_*)<0$. Inspecting the boundary condition \eqref{BCgraph_n} at $v_*$, we obtain that
\bqs
d_{e_*} \partial_n u_{e_*}^n(t_*,v_*)+\alpha_{e_*}^{v_*}u_{e_*}^n(t_*,v_*)+\left(\sum_{e\sim e_*} \nu^{v_*}_{e_*,e} \right)u_{e_*}^n(t_*,v_*)- \sum_{e\sim e_*}\nu^{v_*}_{e,e_*} u_{e}^n(t_*,v_*)=\lambda_{e_*}^{v_*}I_{v_*}^n(t),
\eqs
which writes
\bqs
0> d_{e_*} \partial_n u_{e_*}^n(t_*,v_*)+\delta\left( \sum_{e\sim e_*}\nu^{v_*}_{e,e_*}-\alpha_{e_*}^{v_*}-\sum_{e\sim e_*} \nu^{v_*}_{e_*,e}\right) - \sum_{e\sim e_*}\nu^{v_*}_{e,e_*} (\delta+ u_{e}^n(t_*,v_*)) = \lambda_{e_*}^{v_*}I_{v_*}^n(t) \geq0,
\eqs
and leads to a contradiction. Here we have used the fact that 
\bqs
\sum_{e\sim e_*}\nu^{v_*}_{e,e_*}\leq \alpha_{e_*}^{v_*}+\sum_{e\sim e_*} \nu^{v_*}_{e_*,e},
\eqs 
from Hypothesis~\ref{hypDD} on the matrices $(K_v)_{v\in\V}$. 

Next, from the positivity of solutions, we obtain some uniform $L^\infty$ bounds. More precisely, we claim that there exists a constant $K>0$ depending only on $(T,I_v^0,S_v^0,\|u^0\|_\infty)$ such that
\bqs
0 \leq S_v^n(t), I_v^n(t), R_v^n(t) \leq K, \text{ and } \left| \frac{\md S_v^n(t)}{\md t}\right|, \left| \frac{\md I_v^n(t)}{\md t}\right|, \left| \frac{\md R_v^n(t)}{\md t}\right|\leq K \quad 0<t\leq T,
\eqs
and
\bqs
0\leq  u^n_e (t,x) \leq K , \quad 0<t\leq T, \quad x\in \Omega_e.
\eqs
First, using \eqref{SIRgraph_n} we obtain that
\bqs
\frac{\md S^n_v(t)}{\md t}+\frac{\md I^n_v(t)}{\md t}+\frac{\md R^n_v(t)}{\md t}=-\overline{\lambda}_v I_v^n(t)+\sum_{e\sim v} \alpha^v_e u_e^{n-1}(t,v),
\eqs
which gives that
\bqs
0\leq S_v^n(t)+I_v^n(t)+R_v^n(t) \leq S_v^0+I_v^0+ T \| u^{n-1}\|_\infty, \quad 0<t\leq T,
\eqs
together with
\bqs
\left| \frac{\md S_v^n(t)}{\md t}\right|, \left| \frac{\md I_v^n(t)}{\md t}\right|, \left| \frac{\md R_v^n(t)}{\md t}\right|\leq C\| u^{n-1}\|_\infty (1+ T +T^2\| u^{n-1}\|_\infty) \quad 0<t\leq T,
\eqs
which in turn implies that
\bqs
\| u^{n}\|_\infty \leq \tilde{C} T \| u^{n-1}\|_\infty (1+ T +T^2\| u^{n-1}\|_\infty),
\eqs
for $C,\tilde{C}>0$ only depend on the initial condition $(u^0,(X_v^0)_{v\in\V})$ and the parameters of the system. We now claim that by induction, we have for all $0<t\leq T$,
\begin{align*}
0\leq S_v^n(t)+I_v^n(t)+R_v^n(t) & \leq \hat{C} \sum_{p=0}^{a_n}T^p, \\
\left| \frac{\md S_v^n(t)}{\md t}\right|, \left| \frac{\md I_v^n(t)}{\md t}\right|, \left| \frac{\md R_v^n(t)}{\md t}\right| & \leq \hat{C} \sum_{p=0}^{2 a_n}T^p,\\
\| u^{n}\|_\infty \leq \hat{C} \sum_{p=0}^{2 a_n}T^{p+1},
\end{align*}
with $a_n=2+2a_{n-1}$ for $n\geq2$ and $a_1=1$ for some $\hat{C}>0$ only depending on $(u^0,(X_v^0)_{v\in\V})$. As $0<T<1$, we get that
\bqs
0\leq S_v^n(t)+I_v^n(t)+R_v^n(t) \leq \hat{C}_T , \quad 0<t\leq T,
\eqs
together with
\bqs
\| u^{n}\|_\infty, \left| \frac{\md S_v^n(t)}{\md t}\right|, \left| \frac{\md I_v^n(t)}{\md t}\right|, \left| \frac{\md R_v^n(t)}{\md t}\right| \leq  \hat{C}_T,
\eqs
for some constant $\hat{C}_T>0$ depending on $(T,u^0,(X_v^0)_{v\in\V})$.

\paragraph{Step \#3: existence of a solution.} Parabolic Schauder estimates give that the time derivative and the space derivatives up to order 2 of $u^n$ are uniformly H\"older continuous in compact sets. As a consequence, we can use the Arzela-Ascoli theorem to show that $(u^n,(X_v^n)_{v\in\V})$ converges (up to sequences) toward $(u,(X_v)_{v\in\V})$ in $ \mathcal{C}^{1,2}_{loc}((0,T)\times\cG)\times\mathcal{C}^1_{loc}((0,T))\times\mathcal{C}^1_{loc}((0,T))\times\mathcal{C}^1_{loc}((0,T))$. Passing to the limit $n\rightarrow +\infty$ in \eqref{SIRgraph_n}-\eqref{BCgraph_n}-\eqref{heatgraph_n} we get that $(u,(X_v)_{v\in\V})$   satisfies \eqref{SIRgraph}-\eqref{heatgraph} subject to boundary conditions \eqref{BCgraph}. 

As a by product of the proof we get that for the just constructed solution $(u,(X_v)_{v\in\V})$ we have the uniform bounds:
\bqs
0<S_v(t) \leq S_v^0, \quad 0<t\leq T, \quad v\in\V,
\eqs
and
\bqs
0\leq I_v(t),R_v(t) \quad 0<t\leq T,
\eqs
together with
\bqs
0\leq  u_e (t,x) , \quad 0<t\leq T, \quad x\in \Omega_e, \quad e\in\E.
\eqs
The fact that $I_v(t)\geq0$ implies thanks to the strong maximum principle that in fact 
\bqs
0<  u_e (t,x) , \quad 0<t\leq T, \quad x\in \Omega_e, \quad e\in\E,
\eqs
which in turn gives that $I_v(t)>0$ for each $v\in\V$ since
\bqs
I_v^n(t)=I_v^0e^{-(\eta_v+\overline{\lambda}_v)t+\int_0^t S_v^n(s)\md s}+\sum_{e\sim v} \alpha^v_e \int_0^t e^{-(\eta_v+\overline{\lambda}_v)(t-s)+\int_s^t S_v^n(\tau)\md\tau} u_e^{n-1}(s,v)\md s > 0.
\eqs
Finally, we use the conservation of mass which tells us that
\bqs
\sum_{e\in\E} \int_{\Omega_e} u_e(t,x)\md x+\sum_{v\in\V}\left(S_v(t)+I_v(t)+R_v(t)\right)=M^0>0, \quad 0<t\leq T,
\eqs
such that both $I_v(t)$ and $R_v(t)$ are uniformly bounded in time, together with their derivatives. This also implies that there exists a constant $M>0$, depending only $(u^0,(X_v^0)_{v\in\V})$
such that
\bqs
0<  u (t,x)\leq M, \quad 0<t\leq M, \quad x\in \Omega_e.
\eqs
Using again parabolic regularity, we obtain the solution $(u,(X_v)_{v\in\V})$ is global in time and satisfies \eqref{SIRgraph}-\eqref{heatgraph}-\eqref{BCgraph} in the classical sense.

\subsection{Uniqueness}

Let assume that $(u,(X_v)_{v\in\V})$ and $(\tilde{u},(\tilde{X}_v)_{v\in\V})$ are two classical solutions  to \eqref{heatgraph}-\eqref{BCgraph}-\eqref{SIRgraph} starting from the same initial datum $(u^0,(X_v^0)_{v\in\V})$. We denote $(\hat{u},(\hat{X}_v)_{v\in\V})$ where for each $e\in\E$
\bqs
\hat{u}_e=u_e-\tilde{u}_e,
\eqs
and each $v\in\V$
\bqs
\hat{X}_v=(\hat{S}_v,\hat{I}_v,\hat{R}_v)=(S_v-\tilde{S}_v,I_v-\tilde{I}_v
,R_v-\tilde{R}_v).
\eqs
By linearity, we get that for $e\in\E$
\bqs
\partial_t \hat{u}_e = d_e \partial_x^2 \hat{u}_e,\quad t>0, \quad x\in\overset{\circ}{\Omega}_e,
\eqs
together with
\bqs
D_v \partial_n \hat{\bu}_v(t)+ K_v \hat{\bu}_v(t) = \Lambda_v  \hat{\bI}_v(t), \quad t>0, \quad v\in\V.
\eqs
On the other, one computes that $ \hat{X}_v$ satisfies for each $v\in\V$ 
\bqs
\left\{
\begin{split}
\hat{S}'_v(t)&=-\tau_v \left(S_v(t)\hat{I}_v(t)+\hat{S}_v(t)\tilde{I}_v(t)\right),\\
\hat{I}_v'(t)&=\tau_v \left(S_v(t)\hat{I}_v(t)+\hat{S}_v(t)\tilde{I}_v(t)\right)-\eta_v \hat{I}_v(t)+\sum_{e\sim v} \alpha^v_e \hat{u}_e(t,v)-\overline{\lambda}_v\hat{I}_v(t),\\
\hat{R}_v'(t)&=\eta_v \hat{I}_v(t),
\end{split}
\right.
\eqs
We define the energy 
\bqs
\mathscr{E}(t):=\frac{1}{2}\sum_{e\in\E} \int_{\Omega_e} \left(\hat{u}_e(t,x)\right)^2\md x+\frac{1}{2}\sum_{v\in\V}\left(\hat{S}_v(t)^2+\hat{I}_v(t)^2+\hat{R}_v(t)^2\right),
\eqs
and note that $\mathscr{E}(0)=0$ by definition. Next, differentiating $\mathscr{E}(t)$, we obtain
\begin{align*}
\mathscr{E}'(t)&=\sum_{e\in\E} \int_{\Omega_e} \hat{u}_e(t,x)\partial_t \hat{u}_e(t,x) \md x+\sum_{v\in\V}\left(\hat{S}_v(t)\hat{S}'_v(t)+\hat{I}_v(t)\hat{I}'_v(t)+\hat{R}_v(t)\hat{R}'_v(t)\right)\\
&:=\mathscr{E}_u(t)+\mathscr{E}_X(t).
\end{align*}
On the one hand, we have
\begin{align*}
\mathscr{E}_u(t)&=\sum_{e\in\E} \int_{\Omega_e} \hat{u}_e(t,x)\partial_t \hat{u}_e(t,x) \md x=\sum_{e\in\E}d_e \int_{\Omega_e} \hat{u}_e(t,x)\partial_x^2 \hat{u}_e(t,x) \md x\\
&=-\sum_{e\in\E}d_e \int_{\Omega_e} \left(\partial_x\hat{u}_e(t,x)\right)^2 \md x+\sum_{e\in\E}d_e \left[ \hat{u}_e(t,x)\partial_t \hat{u}_e(t,x)\right]_{\partial \Omega_e}\\
&\leq \sum_{v\in\V} \langle D_v \partial_n \hat{\bu}_v(t),\hat{\bu}_v(t) \rangle \\
&= \sum_{v\in\V} \langle  \Lambda_v  \hat{\bI}_v(t) - K_v \hat{\bu}_v(t),\hat{\bu}_v(t) \rangle \\
&= \sum_{v\in\V} \langle  \Lambda_v  \hat{\bI}_v(t) - A_v \hat{\bu}_v(t),\hat{\bu}_v(t) \rangle -  \sum_{v\in\V} \langle  N_v \hat{\bu}_v(t),\hat{\bu}_v(t) \rangle \\
&\leq  \sum_{v\in\V} \langle  \Lambda_v  \hat{\bI}_v(t) - A_v \hat{\bu}_v(t),\hat{\bu}_v(t) \rangle,
\end{align*}
as $N_v$ is symmetric positive. On the other hand, we compute
\begin{align*}
\mathscr{E}_u(t)&=\sum_{v\in\V}\left(\hat{S}_v(t)\hat{S}'_v(t)+\hat{I}_v(t)\hat{I}'_v(t)+\hat{R}_v(t)\hat{R}'_v(t)\right)\\
&= \sum_{v\in\V}\tau_v\left(S_v(t)\hat{I}_v(t)\hat{S}_v(t)+\hat{S}_v(t)^2\tilde{I}_v(t)+S_v(t)\hat{I}_v(t)^2+\hat{S}_v(t)\tilde{I}_v(t)\hat{I}_v(t) \right)\\
&~~~- \sum_{v\in\V}\eta_v\hat{I}_v(t)^2-\sum_{v\in\V} \overline{\lambda}_v \hat{I}_v(t)^2+\sum_{v\in\V} \hat{I}_v(t)\sum_{e\sim v} \alpha^v_e \hat{u}_e(t,v)+\sum_{v\in\V}\eta_v\hat{R}_v(t) \hat{I}_v(t)\\
&\leq C\mathscr{E}(t)+\sum_{v\in\V} \hat{I}_v(t)\sum_{e\sim v} \alpha^v_e \hat{u}_e(t,v),
\end{align*}
where $C>0$ is some large positive constant. Next, we see that
\bqs
\sum_{v\in\V} \hat{I}_v(t)\sum_{e\sim v} \alpha^v_e \hat{u}_e(t,v)=\sum_{v\in\V} \langle  \hat{\bI}_v(t) , A_v \hat{\bu}_v(t)\rangle,
\eqs
such that we obtain
\bqs
\mathscr{E}'(t)\leq C\mathscr{E}(t)+\sum_{v\in\V} \langle  \Lambda_v  \hat{\bI}_v(t) - A_v \hat{\bu}_v(t),\hat{\bu}_v(t) \rangle +\sum_{v\in\V} \langle  \hat{\bI}_v(t) , A_v \hat{\bu}_v(t)\rangle.
\eqs
Next, if we denote $\hat{\bf{w}}_v(t):=\frac{1}{2}A_v^{-1/2} \left(\Lambda_v +A_v\right) \hat{\bI}_v(t) - A_v^{1/2} \hat{\bu}_v(t)$, we compute
\begin{align*}
0\leq \sum_{v\in\V} \langle  \hat{\bf{w}}_v(t),\hat{\bf{w}}_v(t) \rangle&=\frac{1}{4}\sum_{v\in\V} \langle A_v^{-1} \left(\Lambda_v +A_v\right)^2 \hat{\bI}_v(t)  , \hat{\bI}_v(t) \rangle-\sum_{v\in\V} \langle \left(\Lambda_v +A_v\right) \hat{\bI}_v(t),\hat{\bu}_v(t)\rangle\\
&~~~+\sum_{v\in\V} \langle  A_v \hat{\bu}_v(t),\hat{\bu}_v(t) \rangle\\
&=\frac{1}{4}\sum_{v\in\V} \langle A_v^{-1} \left(\Lambda_v +A_v\right)^2 \hat{\bI}_v(t) , \hat{\bI}_v(t) \rangle -\sum_{v\in\V} \langle  \Lambda_v  \hat{\bI}_v(t) - A_v \hat{\bu}_v(t),\hat{\bu}_v(t) \rangle\\
&~~~ -\sum_{v\in\V} \langle  \hat{\bI}_v(t) , A_v \hat{\bu}_v(t)\rangle.
\end{align*}
As a consequence, we get
\begin{align*}
\mathscr{E}'(t)&\leq C\mathscr{E}(t) +\frac{1}{4}\sum_{v\in\V} \langle A_v^{-1} \left(\Lambda_v +A_v\right)^2 \hat{\bI}_v(t) , \hat{\bI}_v(t) \rangle - \sum_{v\in\V} \langle  \hat{\bf{w}}_v(t),\hat{\bf{w}}_v(t) \rangle\\
&\leq \tilde{C}\mathscr{E}(t),
\end{align*}
for some $\tilde{C}>0$ and we conclude that $\mathscr{E}(t)=0$ for all time which then implies that $\hat{u}=0$ and $\hat{X}_v=0$.

\section{Long-time behavior of the solutions}\label{secLong}

Throughout this section, we denote by $(u,(X_v)_{v\in\V})$ the unique positive bounded classical solution of the Cauchy problem \eqref{SIRgraph}-\eqref{heatgraph}-\eqref{BCgraph} as given by Theorem~\ref{thmexun} and which further satisfies the conservation of total population, namely
\bqs
\sum_{e\in\E} \int_{\Omega_e} u_e(t,x)\md x+\sum_{v\in\V}\left(S_v(t)+I_v(t)+R_v(t)\right)=M^0>0, \quad \forall t>0.
\eqs

\subsection{Final total populations: general results}

As $0<S_v(t) < S_v^0$ and $S_v(t)$ is strictly decreasing, it asymptotically converges towards a limit that we denote
\bqs
S_v^\infty:=\underset{t\rightarrow+\infty}{\lim}S_v(t), \quad v\in\V.
\eqs
Furthermore, as $R_v(t)$ is strictly increasing and uniformly bounded, 
it asymptotically converges towards a limit that is denoted 
\bqs
0<R_v^\infty:=\underset{t\rightarrow+\infty}{\lim}R_v(t)<\infty, \quad v\in\V.
\eqs
But as for each $t>0$
\bqs
R_v(t)=\eta_v \int_0^t I_v(s)\md s,
\eqs
this implies that
\bqs
\mathcal{I}_v(t):=\int_0^tI_v(s)\md s \longrightarrow \mathcal{I}_v^\infty = \frac{R_v^\infty}{\eta_v}<\infty \text{ as } t\rightarrow +\infty,
\eqs
which in turn proves that
\bqs
I_v^\infty= \underset{t\rightarrow+\infty}{\lim}I_v(t)=0.
\eqs

If one recall the notation $m(t)$ for the total population on the edges then we have
\bqs
m(t)=\sum_{e\in\E} \int_{\Omega_e} u_e(t,x)\md x,
\eqs
and it verifies
\bqs
\sum_{v\in\V} \left(S_v(t)+I_v(t)+R_v(t)\right)+m(t)= \sum_{v\in\V} \left(S_v^0+I_v^0\right)+\sum_{e\in\E} \int_{\Omega_e} u_e^0(x)\md x.
\eqs
The above computations shows that $m(t)$ has a limit as $t\longrightarrow+\infty$, that we denote $m_\infty$ and which satisfies
\bqq
\sum_{v\in\V} \left(S_v^\infty+R_v^\infty\right)+m_\infty=
\sum_{v\in\V}\left(S_v^0+I_v^0\right)+\sum_{e\in\E} \int_{\Omega_e} u_e^0(x)\md
x.
\label{Minfini}
\eqq
We shall also keep in mind that
\bqs
S_v^\infty=S_v^0 e^{-\tau_v \mathcal{I}_v^\infty}=S_v^0 e^{-\frac{\tau_v}{\eta_v} R_v^\infty}, \text{ or } R_v^\infty = -\frac{\eta_v}{\tau_v}\ln \frac{S_v^\infty}{S_v^0}, \quad v\in\V
\eqs
And so if we introduce the function $\Psi_v(x):=x-\frac{\eta_v}{\tau_v}\ln x$, then the above conservation of mass can be written as 
\bqs
\sum_{v\in\V}\Psi_v(S_v^\infty)+m_\infty=\sum_{v\in\V}\left(I_v^0+\Psi_v(S_v^0)\right)+\sum_{e\in\E} \int_{\Omega_e} u_e^0(x)\md x.
\eqs

On the other, one can compute that
\bqs
\frac{\md m(t) }{\md t}=\sum_{v\in\V}\overline{\lambda}_vI_v(t)-\sum_{v\in\V} \langle A_v  \bu_v(t) ,\mathbf{1}_{\delta_v}\rangle,
\eqs
such that
\bqs
m(t)+\sum_{v\in\V} \int_0^t \langle A_v  \bu_v(s) ,\mathbf{1}_{\delta_v}\rangle \md s = \sum_{e\in\E} \int_{\Omega_e} u_e^0(x)\md x+\sum_{v\in\V}\overline{\lambda}_v\mathcal{I}_v(t).
\eqs
Now, as $m(t)$ and each $\mathcal{I}_v(t)$ are convergent we deduce that all $\int_0^t \bu_v(s)\md s$ are also convergent so that
\bqq
m_\infty+\sum_{v\in\V} \int_0^\infty \langle A_v  \bu_v(s) ,\mathbf{1}_{\delta_v}\rangle \md s = \sum_{e\in\E} \int_{\Omega_e} u_e^0(x)\md x+\sum_{v\in\V}\overline{\lambda}_v\mathcal{I}_v^\infty,
\label{eqmass2}
\eqq
and
\bqs
\int_0^\infty\bu_v(s)\md s <\infty, \quad v\in\V,
\eqs
which proves that
\bqs
\bu_v(t)\longrightarrow 0 \text{ as } t\rightarrow+\infty, \quad v\in\V.
\eqs
And the boundary conditions imply that
\bqs
\partial_n \bu_v(t)\longrightarrow 0 \text{ as } t\rightarrow+\infty, \quad v\in\V.
\eqs
Next, we define the sequence of functions $u^n_e(t,x)=u_e(t+n,x)$ for each $e\in\E$ and $X^n_v(t)=X_v(t+n)$ for each $v\in\V$ which are uniformly bounded such that one can extract a convergent subsequence. On the one hand we have that $\underset{n\rightarrow\infty}{\lim}X^n_v(t)=X_v^\infty=(S_v^\infty,0,R_v^\infty)$ and on the other if  $u^\infty_e(t,x)=\underset{n\rightarrow\infty}{\lim}u^n_e(t,x)$ it is solution of
\bqs
\partial_t u_e^\infty(t,x) = d_e\partial_x^2 u_e^\infty(t,x),
\eqs
with the boundary conditions 
\bqs
\partial_n \bu_v^\infty(t)= \bu_v^\infty(t)=\mathbf{0}_{\delta_v}, \quad v\in\V.
\eqs
This then shows that $u_e^\infty(t,x)=0$, $t>0$ and $x\in\overset{\circ}{\Omega}_e$ for each $e\in\E$. As there is unicity of the limit, we deduce that
\bqs
\underset{t\rightarrow+\infty}{\lim}u_e(t,x)=0, \quad e\in\E.
\eqs
From which we also get that $m_\infty=0$ and that
\bqs
\sum_{v\in\V}\Psi_v(S_v^\infty)=\sum_{v\in\V}\left(I_v^0+\Psi_v(S_v^0)\right)+\sum_{e\in\E} \int_{\Omega_e} u_e^0(x)\md x.
\eqs
This implies that each $\Psi_v(S_v^\infty)$ is bounded, we get that $S_v^\infty>0$ for all $v\in\V$.

We also get from \eqref{eqmass2}, that 
\bqs
\sum_{v\in\V} \int_0^\infty \langle A_v  \bu_v(s) ,\mathbf{1}_{\delta_v}\rangle \md s = \sum_{e\in\E} \int_{\Omega_e} u_e^0(x)\md x+\sum_{v\in\V}\overline{\lambda}_v\mathcal{I}_v^\infty.
\eqs

Finally, we use the fact that
\bqs
\frac{\md I_v(t)}{\md t}+\frac{\md S_v(t)}{\md t}-\frac{\eta_v+\overline{\lambda}_v}{\tau_v} \frac{\md \ln S_v(t)}{\md t} = \sum_{e\sim v} \alpha^v_e u_e(t,v), \quad v\in\V,
\eqs
to obtain that
\bqs
 I_v(t)+S_v(t)-\frac{\eta_v+\overline{\lambda}_v}{\tau_v}  \ln S_v(t) -I_v^0-S_v^0 +\frac{\eta_v+\overline{\lambda}_v}{\tau_v}  \ln S_v^0= \sum_{e\sim v} \alpha^v_e \int_0^tu_e(s,v)\md s, \quad v\in\V.
\eqs

As a consequence, the final total populations of infected individuals at each vertices satisfy the following scalar differential equation
\bqq
\frac{\md \mathcal{I}_v(t)}{\md t}=S_v^0\left(1-e^{-\tau_v \mathcal{I}_v(t)}\right)-\eta_v \mathcal{I}_v(t)+I_v^0 +\sum_{e\sim v} \alpha^v_e \int_0^tu_e(s,v)\md s -\overline{\lambda}_v \mathcal{I}_v(t), \quad v\in\V.
\eqq

Passing to the limit as $t\rightarrow+\infty$, we get
\bqs
0= S_v^0\left(1-e^{-\tau_v \mathcal{I}_v^\infty}\right)-\eta_v \mathcal{I}_v^\infty+I_v^0 +\sum_{e\sim v} \alpha^v_e \int_0^\infty u_e(s,v)\md s -\overline{\lambda}_v \mathcal{I}_v^\infty, \quad v\in\V.
\eqs
%\bqs
%S_v^\infty-\frac{\eta+\overline{\lambda}_v}{\tau}  \ln S_v^\infty= I_v^0+S_v^0-\frac{\eta+\overline{\lambda}_v}{\tau}  \ln S_v^0+ \sum_{e\sim v} \alpha^v_e \int_0^\infty u_e(s,v)\md s, \quad v\in\V,
%\eqs
%which also writes
%\bqs
%\Psi(S_v^\infty)+\overline{\lambda}_v \mathcal{I}_v^\infty=I_v^0+\Psi(S_v^0)+ \sum_{e\sim v} \alpha^v_e \int_0^\infty u_e(s,v)\md s, \quad v \in\V.
%\eqs

To summarize, we have proved the following result.

\begin{thm}\label{thmFTI}
For each $(S_v^0,I_v^0)\geq \mathbf{0}$ with $S_v^0>0$, $\sum_{v\in\V} I_v^0>0$ and $u^0\in BC(\cG,\R^+)$ that satisfy the boundary condition \eqref{BCgraph}, the long time behavior of the unique corresponding solution $(u,(X_v)_{v\in\V})$ is given by
\bqs
\underset{t\rightarrow+\infty}{\lim}u_e(t,x)=0, \quad x\in\Omega_e, \quad e\in\E, \quad \text{ with }\quad \int_0^\infty u_e(s,v)\md s < +\infty,  \quad (v,e)\in\cG,
\eqs
and
\bqq
\underset{t\rightarrow+\infty}{\lim}~(S_v(t),I_v(t),R_v(t)) =\left(S_v^0e^{-\tau_v
    \mathcal{I}_v^\infty},0,\eta_v \mathcal{I}_v^\infty \right), \quad v\in\V,
\label{ValInfinity}
\eqq
where the final total populations of infected individuals $0<\mathcal{I}_v^\infty<\infty$ at each vertices $v\in\V$ are solutions of the system
\bqq
\left\{
\begin{split}
S_v^0e^{-\tau_v \mathcal{I}_v^\infty}+\eta_v \mathcal{I}_v^\infty&=I_v^0+S_v^0+ \sum_{e\sim v} \alpha^v_e \int_0^\infty u_e(s,v)\md s - \overline{\lambda}_v\mathcal{I}_v^\infty, \quad v \in\V,\\
\sum_{v\in\V} \int_0^\infty \langle A_v  \bu_v(s) ,\mathbf{1}_{\delta_v}\rangle \md s &= \sum_{e\in\E} \int_{\Omega_e} u_e^0(x)\md x+\sum_{v\in\V}\overline{\lambda}_v\mathcal{I}_v^\infty.
\end{split}
\right.
\label{systemSvinf}
\eqq
As a consequence, $(\mathcal{I}_v^\infty)_{v\in\V}$ belongs to the parametrized submanifold given by
\bqq
\sum_{v\in\V}\left(S_v^0e^{-\tau_v \mathcal{I}_v^\infty}+\eta_v \mathcal{I}_v^\infty \right)=M^0.
\label{manifoldIvinf}
\eqq
\end{thm}

\begin{rmk}
Equivalently, $(S_v^\infty)_{v\in\V}$ belongs to the parametrized submanifold given by
\bqq
  \sum_{v\in\V}\left(S_{v}^\infty-\frac{\eta_v}{\tau_v} \log\left({S_{v}^\infty}\right)+\frac{\eta_v}{\tau_v} \log\left({S_{v}^0}\right)\right)=M^0,
\label{manifoldSvinf}
\eqq
and $(R_v^\infty)_{v\in\V}$ belongs to the parametrized submanifold given by
\bqq
  \sum_{v\in \V} \left(S^0_v\exp(-\tau_v/\eta_v\, R_v^\infty) + R_v^\infty\right) = M^0.
\label{manifoldRvinf}
\eqq
The equations \eqref{manifoldIvinf}, \eqref{manifoldSvinf}, and \eqref{manifoldRvinf} also read
$\sum_{v\in\V}\left(S_v^\infty+R_v^\infty\right)=M^0$, which is nothing but
\eqref{Minfini} since we have proved that $m_\infty=0$.
\end{rmk}

\begin{rmk}
If we assume that $\tau=\tau_v>0$ and $\eta=\eta_v>0$ are independent of $v\in\V$ and let $\widetilde{S_{v}}={\tau}/{\eta}\, {S_{v}}$, $\widetilde{R_v^\infty} =
\exp(-\tau/\eta\, R_v^\infty)$ and
$\widetilde{\mathcal{I}_v^\infty}=\exp(-\tau \mathcal{I}_v^\infty) $. Then,
equations \eqref{manifoldIvinf}, \eqref{manifoldSvinf}, and
\eqref{manifoldRvinf} are respectively equivalent to
\[
  \prod_{v\in\V}\exp\left(\widetilde{S_{v}^\infty}\right) \frac{\widetilde{S_{v}^0}}{\widetilde{S_{v}^\infty}}=\exp\left(\frac{\tau}{\eta}M^0\right),
\]
\[
  \prod_{v\in \V}
  \frac{\exp\left(\widetilde{S_v^0}\widetilde{R_v^\infty}\right)}{\widetilde{R_v^\infty}}=\exp\left(\frac{\tau}{\eta}M^0\right),
\]
and
\[
\prod_{v\in\V}\frac{\exp\left(\widetilde{S_v^0}\widetilde{\mathcal{I}_v^\infty}\right)}{\widetilde{\mathcal{I}_v^\infty}}
=\exp\left(\frac{\tau}{\eta}M^0\right).
\]
The common right hand side features $\frac{\tau}{\eta}\, M^0$ that is nothing but the
traditional basic reproductive number $\mathscr{R}_0$.
\end{rmk}

\subsection{Final total populations of infected individuals: further properties}

The aim of this section is to present further qualitative results on the final total configuration $(\mathcal{I}_v^\infty)_{v\in\V}$ in the fully symmetric case where one can obtain closed form formula and in the case of two vertices where we manage to obtain sharp bounds on the final total populations of infected individuals. In each case, we manage to relate these quantities to standard basic and effective reproductive number for classical SIR model \cite{DHM90}.

\paragraph{Fully symmetric case.}
We  assume that the
length $\ell_e$ of every edge $e \in \mathcal{E}$ is equal to a
reference length $\ell$. For every $e\in \mathcal{E}$, the diffusion coefficient
$d_e$ is equal to $d$. We moreover suppose that for every vertex $v\in 
\mathcal{V}$, $S_v^0=S^0$, $I_v^0=I^0$ and $R_v^0=R^0$. We also assume that $\tau=\tau_v>0$ and $\eta=\eta_v>0$ are independent of $v\in\V$.  In the same spirit,
$\lambda_e^v=\lambda$ and $\alpha_e^v=\alpha$ for every $e\in \mathcal{E}$ and $v\in
\mathcal{V}$. We also assume $\nu_{e_i,e_j}^v=\nu$ for every edges incident to the
vertex $v$. Finally, the components $u_e^0$ of initial condition on each edges
$e\in \mathcal{E}$ are supposed to be even with respect to the center of the
interval $\Omega_e=[0,\ell]$. Thanks to all these assumptions,
$\mathcal{I}_v^\infty$ does not depend on the vertex $v\in \mathcal{V}$ and we
set $\mathcal{I}_v^\infty=\mathcal{I}^\infty$. Let us recall the notation $c_{\mathcal{V}}$
 for the cardinal of the set $\mathcal{V}$. The parametrized submanifold given by
\eqref{manifoldIvinf} becomes
\begin{equation*}
  S^\infty+R^\infty=S^0e^{-\tau \mathcal{I}^\infty}+\eta \mathcal{I}^\infty =\widetilde{M^0},
\end{equation*}
where $\widetilde{M^0}=M^0/c_{\mathcal{V}}$. We can transform this relation as
\begin{equation}
  \label{manifoldIvinfsym}
  S^0e^{-\tau \mathcal{I}^\infty}+\frac{\eta}{\tau} \tau\mathcal{I}^\infty -\widetilde{M^0}=0.
\end{equation}
Let $\mathscr{I}=-\tau \mathcal{I}^\infty$. We have to solve
\[
  S^0 e^{\mathscr{I}}-\frac{\eta}{\tau}\mathscr{I}-\widetilde{M^0}=0.
\]
The solutions are given in terms of Lambert W function that is the multivalued
inverse relation of the function $f(w) = we^w$ for $w\in 
\mathbb{C}$ \cite{lambert}. Let us recall how to compute the real
solutions of the equation $\alpha e^x+\beta x+\gamma=0$ for
$(\alpha,\beta,\gamma)\in \mathbb{R}^*\times \mathbb{R}^*\times \mathbb{R}$. Let
$\Delta=\alpha/\beta\, \exp(-\gamma/\beta)$ be the discriminant. If $\Delta\geq 0$ or
$\Delta=-\exp(-1)$, the solution is unique and $x=-W_0(\Delta)-\gamma/\beta$ where $W_0$
is the principal branch. If $\Delta \in (-\exp(-1),0)$, there are two solutions
$x_0=-W_0(\Delta)-\gamma/\beta$ and $x_{-1}=-W_{-1}(\Delta)-\gamma/\beta$, where
$W_{-1}$ is another branch. When $\Delta<-\exp(-1)$, there is no solution.

In our symmetric case, the discriminant writes
\[
  \Delta=-\frac{S^0\tau}{\eta}\exp\left(-\frac{\widetilde{M^0}\tau}{\eta}\right).
\]
Since $\Delta<0$, there exist solutions to \eqref{manifoldIvinfsym} if
$\Delta\geq -\exp(-1)$, which is equivalent to
\begin{equation}
  \exp\left(\frac{\widetilde{M^0}\tau}{\eta}-1\right) \geq
  \frac{S^0\tau}{\eta}.\label{discriminant}
\end{equation}
We recall that when we consider the standard SIR model (meaning in the context
of this paper that we consider an isolated vertex), we can define the effective
reproductive number $\mathscr{R}_e$ and the basic reproductive number $\mathscr{R}_0$
respectively given by
\bqq
  \mathscr{R}_e:=\frac{S^0\tau}{\eta}, \quad \text{and} \quad \mathscr{R}_0:=\frac{{M^0}\tau}{\eta},
  \label{basicRN}
\eqq
see \cite{DHM90,VanW02,Heth00} for further properties of effective and basic 
reproductive numbers. If we denote $\widetilde{\mathscr{R}_0}=\widetilde{M^0}\tau/\eta$, the equation
\eqref{discriminant} reads
\[
  \exp\left(\widetilde{\mathscr{R}_0}-1\right)\geq \mathscr{R}_e.
\]
This inequality is satisfied as long as $S^0\leq \widetilde{M^0}$, which is
always true since $M^0= \sum_{e\in \mathcal{E}}\int_{\Omega_e}u_e^0(x) \md
x+c_{\V}\left(S^0+I^0\right) \geq c_{\V} S^0$. Since $\Delta=-\mathscr{R}_e \exp(-\widetilde{\mathscr{R}_0})$, the solutions are
\[
  \mathscr{I}_{0,-1}=-W_{0,-1}\left(-\mathscr{R}_e \exp(-\widetilde{\mathscr{R}_0})\right)-\widetilde{\mathscr{R}_0},
\]
and so
\[
  \mathcal{I}^\infty_{0,-1}=W_{0,-1}\left(-\mathscr{R}_e
    \exp(-\widetilde{\mathscr{R}_0})\right)/\tau+\widetilde{\mathscr{R}_0}/\tau.
\]
Both $W_{0,-1}\left(-\mathscr{R}_e \exp(-\widetilde{\mathscr{R}_0})\right)<0$. However, we can show
that $\mathcal{I}_0^\infty>0$ and $\mathcal{I}_{-1}^\infty<0$. Thus, the only possibility is
\[
  \mathcal{I}^\infty=W_{0}\left(-\mathscr{R}_e
    \exp(-\widetilde{\mathscr{R}_0})\right)/\tau+\widetilde{\mathscr{R}_0}/\tau.
\]
We also have access to $S^\infty$ and $R^\infty$ thanks to
\eqref{ValInfinity}. Since $\exp(-W_0(x))=W_0(x)/x$, we obtain
\[
  S^\infty=-\frac{\eta}{\tau}W_{0}\left(-\mathscr{R}_e\exp(-\widetilde{\mathscr{R}_0})\right),
\]
and
\[
  R^\infty=\frac{\eta}{\tau} W_{0}\left(-\mathscr{R}_e
    \exp(-\widetilde{\mathscr{R}_0})\right)+\widetilde{M^0}.
\]

We can summarize these results in the following lemma.

\begin{figure}[!htbp]
  \centering
  \begin{tabular}{cc}
    \includegraphics[height=.29\textheight]{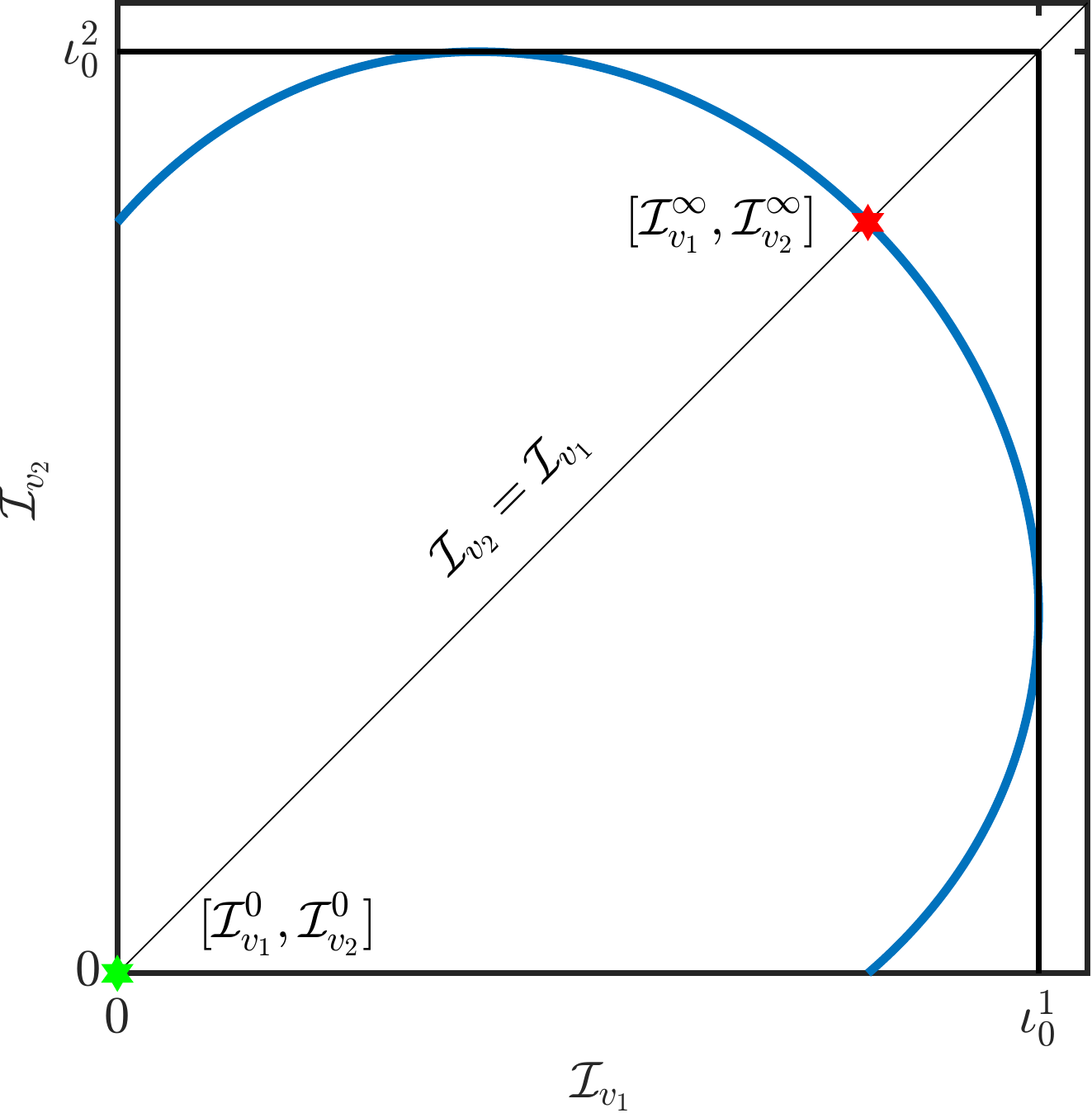} &
   \includegraphics[height=.29\textheight]{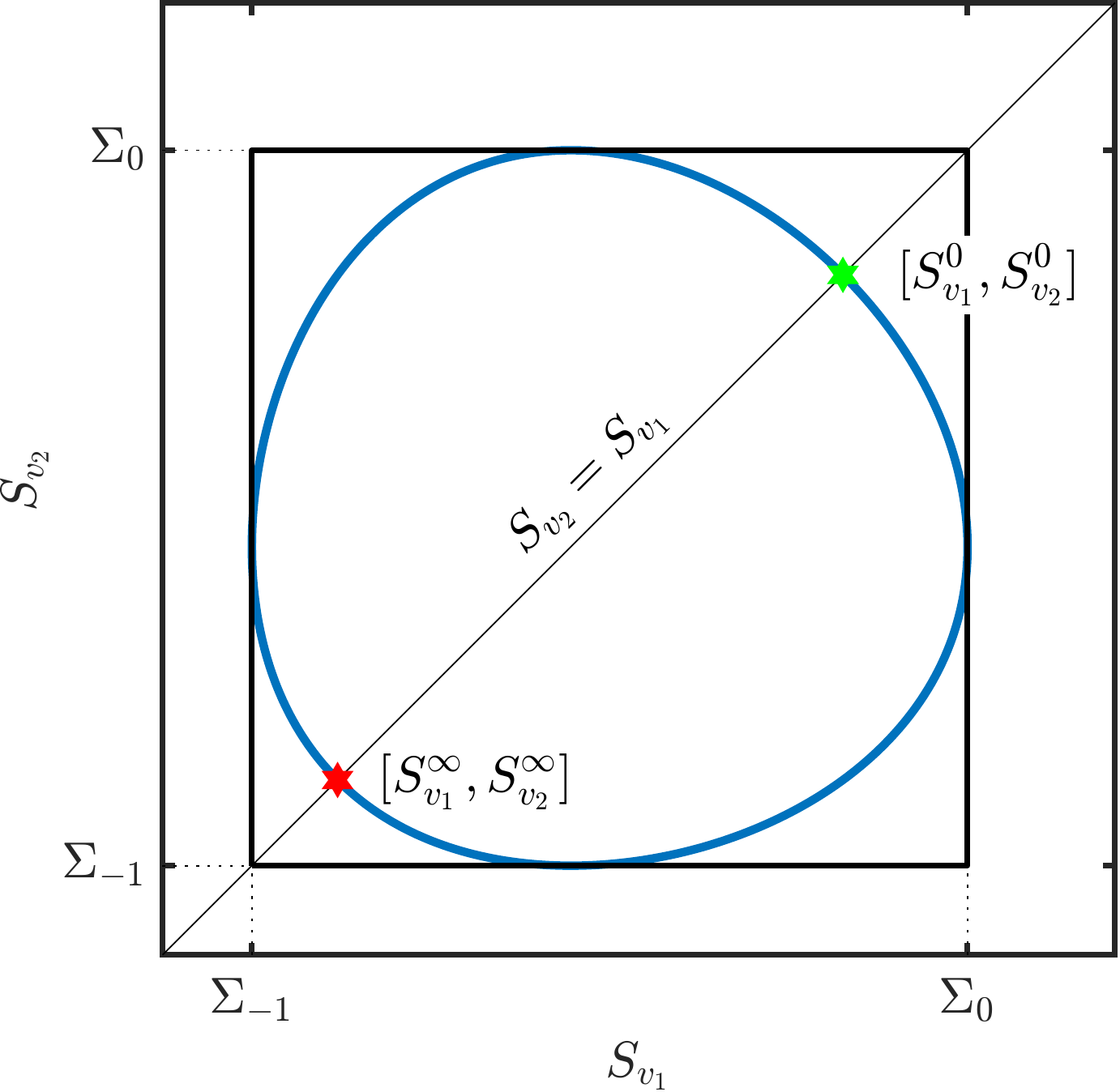}\\
  \end{tabular}
  \caption{Schematic visualisation (red star) of $\mathcal{I}^\infty=\mathcal{I}_{v_1}^\infty=\mathcal{I}_{v_2}^\infty$, resp. $S^\infty=S_{v_1}^\infty=S_{v_2}^\infty$, in the $(\mathcal{I}_{v_1},\mathcal{I}_{v_2})$ plane, resp.  in the $(S_{v_1},S_{v_2})$-plane, in the fully symmetric case. The asymptotic value $\mathcal{I}^\infty$, resp. $S^\infty$, lies at the intersection of the diagonal $\mathcal{I}_{v_1}=\mathcal{I}_{v_2}$, resp. $S_{v_1}=S_{v_2}$, and the implicit curve given by \eqref{manifoldIvinf}, resp. \eqref{manifoldSvinf}.}
  \label{fig:SIinf_sym}
\end{figure}

\begin{lem}[Fully symmetric case.]\label{lem_sym}
Assume that our model is fully symmetric, then the final total population of infected individuals as given by Theorem~\ref{thmFTI} is independent on the vertex that is $\mathcal{I}_v^\infty=\mathcal{I}^\infty$ for each $v\in\V$, and $\mathcal{I}^\infty$ has the following closed form formula
\[
  \mathcal{I}^\infty=\frac{W_{0}\left(-\mathscr{R}_e
    \exp(-\mathscr{R}_0/c_\V)\right)}{\tau}+\frac{\mathscr{R}_0}{c_\V\,\tau},
\]
where $\mathscr{R}_e$ and  $\mathscr{R}_0$
 are respectively the effective and basic 
reproductive number defined in \eqref{basicRN} and $c_\V$ the cardinal of $\V$. See Figure~\ref{fig:SIinf_sym} for an illustration.
\end{lem}

\paragraph{Case of two vertices.}
In this simple case, it is possible to build explicit formulas to deal with the
implicit submanifold equations \eqref{manifoldIvinf}, \eqref{manifoldSvinf}, and
\eqref{manifoldRvinf}. Let $\mathscr{R}_{0,v_k}:=M^0\, \tau_{v_k}/\eta_{v_k}$
  and $\mathscr{R}_{e,v_k}:=S_{v_k}^0\, \tau_{v_k}/\eta_{v_k}$, $k=1,2$ be respectively the
local to vertex $v_k$ basic and effective reproductive number. Then,
\begin{equation}
  S_{v_1}^\infty =
  -\frac{\eta_{v_1}}{\tau_{v_1}}W\left(-\exp\left(-\mathscr{R}_{0,v_1}\right)\mathscr{R}_{e,v_1}
    \left(\mathscr{R}_{e,v_2}\right)^{\frac{\tau_{v_1}\eta_{v_2}}{\tau_{v_2}\eta_{v_1}}}\frac{\exp\left(S_{v_2}^\infty\,
        {\tau_{v_1}}/{\eta_{v_1}}\right)}{\left({S_{v_2}^\infty}\, {\tau_{v_2}}/{\eta_{v_2}}\right)^{\frac{\tau_{v_1}\eta_{v_2}}{\tau_{v_2}\eta_{v_1}}}}
  \right),
\label{eq:Sv1}
\end{equation}
% \begin{equation}
%   S_{v_1}^\infty =
%   -\frac{\eta}{\tau}W\left(-\exp\left(-\mathscr{R}_0\right){\mathscr{R}_{e,v_1}}\mathscr{R}_{e,v_2}
%     \frac{\exp\left({\tau}/{\eta}\, S_{v_2}^\infty\right)}{{\tau}/{\eta}\,
%       S_{v_2}^\infty}\right),\label{eq:Sv1}
% \end{equation}
where the Lambert W function $W$ can be either $W_0$ or $W_{-1}$. Indeed, the argument
of $W$ being negative, two solutions have to be considered. We obviously also have
\begin{equation}
  S_{v_2}^\infty =
  -\frac{\eta_{v_2}}{\tau_{v_2}}W\left(-\exp\left(-\mathscr{R}_{0,v_2}\right)\mathscr{R}_{e,v_2}
    \left(\mathscr{R}_{e,v_1}\right)^{\frac{\tau_{v_2}\eta_{v_1}}{\tau_{v_1}\eta_{v_2}}}\frac{\exp\left(S_{v_1}^\infty\,
        {\tau_{v_2}}/{\eta_{v_2}}\right)}{\left({S_{v_1}^\infty}\, {\tau_{v_1}}/{\eta_{v_1}}\right)^{\frac{\tau_{v_2}\eta_{v_1}}{\tau_{v_1}\eta_{v_2}}}}
  \right),
\label{eq:Sv2}
\end{equation}
% \begin{equation}
%   S_{v_2}^\infty =
%   -\frac{\eta}{\tau}W\left(-\exp\left(-\mathscr{R}_0\right){\mathscr{R}_{e,v_1}}\mathscr{R}_{e,v_2}
%     \frac{\exp\left({\tau}/{\eta}\, S_{v_1}^\infty\right)}{{\tau}/{\eta}\,
%       S_{v_1}^\infty}\right),\label{eq:Sv2}
% \end{equation}
Due to the
definition of the domain of the Lambert W function, the argument
has to be greater than $-\exp(-1)$. So, the following inequality must be
satisfied for $S_{v_2}^\infty$ (respectively of $S_{v_1}^\infty$)
\[
  \frac{-\exp\left(-\mathscr{R}_{0,v_1}\right){\mathscr{R}_{e,v_1}}\left(\mathscr{R}_{e,v_2}\right)^{\frac{\tau_{v_1}\eta_{v_2}}{\tau_{v_2}\eta_{v_1}}}}{\left(S_{v_2}^\infty\,
      {\tau_{v_2}}/{\eta_{v_2}}\right)^{\frac{\tau_{v_1}\eta_{v_2}}{\tau_{v_2}\eta_{v_1}}}\exp\left(-S_{v_2}^\infty\,
      {\tau_{v_1}}/{\eta_{v_1}}\right)}\geq-\exp(-1).
\]
% \[
% -\exp\left(-\mathscr{R}_0\right){\mathscr{R}_{e,v_1}}\mathscr{R}_{e,v_2}
%     \frac{\exp\left({\tau}/{\eta}\, S_{v_2}^\infty\right)}{{\tau}/{\eta}\, S_{v_2}^\infty}\geq -\exp(-1).\]
Solving the equality part of this inequality, we find that
\[
  S_{v_2}^\infty = -\frac{\eta_{v_2}}{\tau_{v_2}}W\left(-\left({\mathscr{R}_{e,v_1}}\right)^{\frac{\tau_{v_2}\eta_{v_1}}{\tau_{v_1}\eta_{v_2}}}\mathscr{R}_{e,v_2}\exp\left(\frac{\tau_{v_2}\eta_{v_1}}{\tau_{v_1}\eta_{v_2}}\left(1-\mathscr{R}_{0,v_1}\right)\right)\right).
\]
% \[
%   S_{v_2}^\infty=-\frac{\eta}{\tau}W\left(-\exp\left(1-\mathscr{R}_0\right){\mathscr{R}_{e,v_1}}\mathscr{R}_{e,v_2}\right).
% \]
This equation has to be verified both for $W_0$ and $W_{-1}$. Let
$\Sigma_{0,-1}^{v_2}$ be defined by
\[
  \Sigma_{0,-1}^{v_2}:=-\frac{\eta_{v_2}}{\tau_{v_2}}W_{0,-1}\left(\mathscr{A}_{v_2}\right),
\]
where
\begin{equation}
  \mathscr{A}_{v_2}=\left({\mathscr{R}_{e,v_1}}\right)^{\frac{\tau_{v_2}\eta_{v_1}}{\tau_{v_1}\eta_{v_2}}}\mathscr{R}_{e,v_2}\exp\left(\frac{\tau_{v_2}\eta_{v_1}}{\tau_{v_1}\eta_{v_2}}\left(1-\mathscr{R}_{0,v_1}\right)\right)\label{eq:Av1}
\end{equation}

\begin{figure}[!htbp]
  \centering
  \begin{tabular}{cc}
    \includegraphics[width=.435\textwidth]{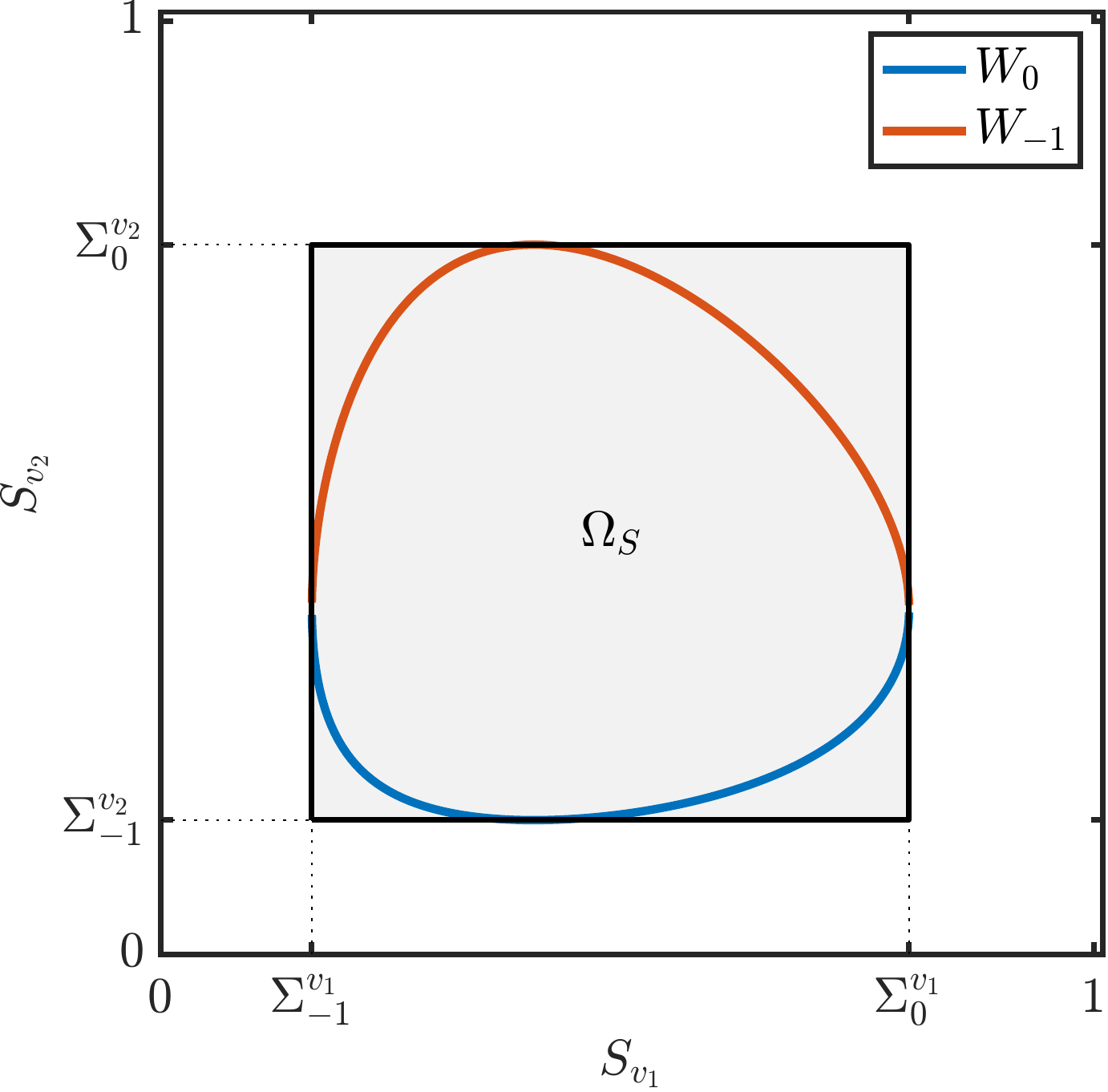} &
   \includegraphics[width=.435\textwidth]{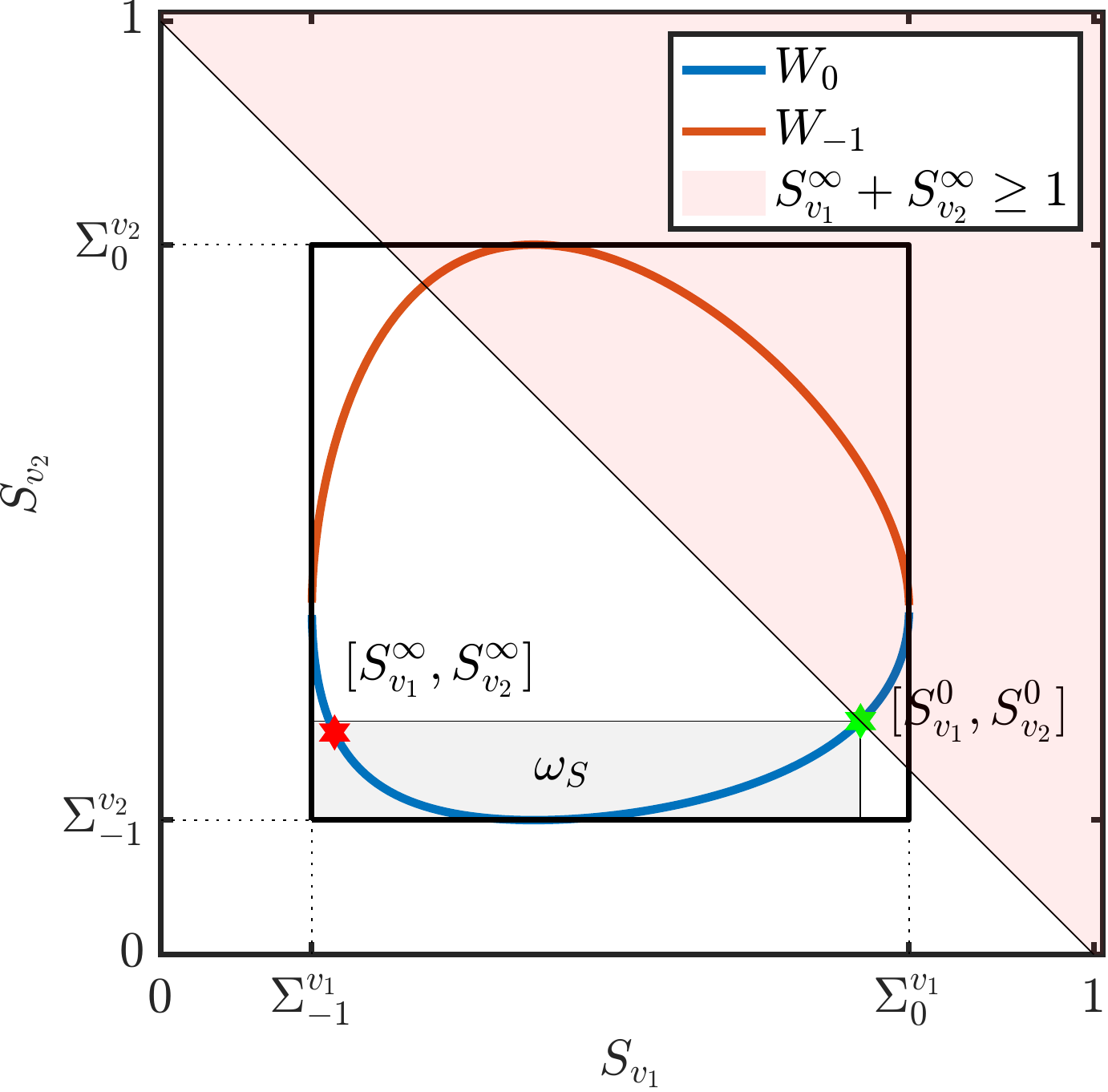}\\
  \end{tabular}
  \caption{Location of $S_{v_1}^\infty$ and $S_{v_2}^\infty$ together with the  visualisation of the domains 
    $\Omega_S$ (left) and  $\omega_S$ (right). The final configuration of
    susceptible individuals $(S_{v_1}^\infty,S_{v_2}^\infty)$ lies on the closed
    curve parametrized by the two branches of the Lambert W function (blue and
    red curve). We note that $\left(S_{v_1}^\infty,S_{v_2}^\infty\right) \in
    \omega_S$ as indicated by the red star on the right figure. Values of the
    parameters are $d=10^{-3}$, $\lambda^1=\lambda^2=6/10$,
    $\alpha^1=\alpha^2=1/8$, $\tau_{v_1}=1$, $\tau_{v_2}=9/10$, $\eta_{v_1}=2/5$, $\eta_{v_2}=2/6$, and initial conditions are set to: 
$I_1^0=I_2^0=10^{-6}$, $S_0^1=3/4-I_1^ 0$, $S_0^2=1/4-I_1^0$ and $u^0(x)=0$. The mass $M^0$ is therefore equal to $1$. }
  \label{fig:boxS}
\end{figure}

% \[
%   \Sigma_{0,-1}:=-\frac{\eta}{\tau}W_{0,-1}\left(-\exp\left(1-\mathscr{R}_0\right){\mathscr{R}_{e,v_1}}\mathscr{R}_{e,v_2}\right).
% \]
Then, the domain of $S_{v_1}^\infty$ as a function of $S_{v_2}^\infty$ is
\[
  S_{v_2}^\infty \in \left[\min\left(\Sigma_{-1}^{v_2},\Sigma_{0}^{v_2}\right),\max\left(\Sigma_{-1}^{v_2},\Sigma_{0}^{v_2}\right)\right].
\]
Concerning $S_{v_2}^\infty$ as a function of $S_{v_1}^\infty$, we have
\[
  S_{v_1}^\infty \in \left[\min\left(\Sigma_{-1}^{v_1},\Sigma_{0}^{v_1}\right),\max\left(\Sigma_{-1}^{v_1},\Sigma_{0}^{v_1}\right)\right],
\]
with
\[
  \Sigma_{0,-1}^{v_1}:=-\frac{\eta_{v_1}}{\tau_{v_1}}W_{0,-1}\left(\mathscr{A}_{v_1}\right),
\]
where
\begin{equation}
  \mathscr{A}_{v_1}={\mathscr{R}_{e,v_1}}\left(\mathscr{R}_{e,v_2}\right)^{\frac{\tau_{v_1}\eta_{v_2}}{\tau_{v_2}\eta_{v_1}}}\exp\left(\frac{\tau_{v_1}\eta_{v_2}}{\tau_{v_2}\eta_{v_1}}\left(1-\mathscr{R}_{0,v_2}\right)\right).\label{eq:Av2}
\end{equation}
Thus,
\[
  \left(S_{v_1}^\infty,S_{v_2}^\infty\right) \in
  \Omega_S:=\left[\min\left(\Sigma_{-1}^{v_1},\Sigma_{0}^{v_1}\right),\max\left(\Sigma_{-1}^{v_1},\Sigma_{0}^{v_1}\right)\right]\times\left[\min\left(\Sigma_{-1}^{v_2},\Sigma_{0}^{v_2}\right),\max\left(\Sigma_{-1}^{v_2},\Sigma_{0}^{v_2}\right)\right]
  .
\]
We present on Figure \ref{fig:boxS} (left) the functions $W_0$ and $W_{-1}$
defining $S_{v_2}^\infty$ as a function of $S_{v_1}^\infty$ and the domain $\Omega_S$ for a given set of the parameters and initial conditions. We refer to Section~\ref{secNumSc} for details regarding the numerical integration of the model and Section~\ref{secSim} for further numerical results on the case of two vertices.

%$I_1^0=I_2^0=10^{-6}$, $S_0^1=3/4-I_1^ 0$, $S_0^2=1/4-I_1^0$ and $u^0(x)=0$. The
%mass $M^0$ is therefore equal to $1$. The
%parameters are $d=10^{-3}$, $\lambda^1=\lambda^2=6/10$, $\alpha^1=\alpha^2=1/8$.
Actually, we can reduce the domain of validity of \eqref{eq:Sv1}-\eqref{eq:Sv2} for $S_{v_1}^\infty$
and $S_{v_2}^\infty$. Indeed, we know that $S_{v_k}$, $k=1,2$, decay with
respect to time, so $S_{v_k}^\infty<S_{v_k}$. Moreover, the sum
$S_{v_1}^\infty+S_{v_2}^\infty<M^0$. Thus, we have
\[
  \left(S_{v_1}^\infty,S_{v_2}^\infty\right) \in \omega_S:=[\min\left(\Sigma_{-1}^{v_1},\Sigma_{0}^{v_1}\right),S_{v_1}^0]\times
  [\min\left(\Sigma_{-1}^{v_2},\Sigma_{0}^{v_2}\right),S_{v_2}^0]\, \cap \left\{ S_{v_1}^\infty+S_{v_2}^\infty<M^0\right\}.
\]
The domain $\omega_S$ is drawn on Figure \ref{fig:boxS} (right).

Concerning $\mathcal{I}_{v_1}^\infty$ and $\mathcal{I}_{v_2}^\infty$, we can
perform the same analysis. Let
\[
  \mathcal{J}_{v_1}^\infty = -\mathscr{R}_{0,v_1}+\frac{\tau_{v_1}\eta_{v_2}}{\tau_{v_2}\eta_{v_1}} \left(\mathscr{R}_{e,v_2}\exp^{-\tau_{v_2} \mathcal{I}_{v_2}^\infty} + \tau_{v_2}\mathcal{I}_{v_2}^\infty\right),
\]
and
\[
  \mathcal{J}_{v_2}^\infty = -\mathscr{R}_{0,v_2}+\frac{\tau_{v_2}\eta_{v_1}}{\tau_{v_1}\eta_{v_2}} \left(\mathscr{R}_{e,v_1}\exp^{-\tau_{v_1} \mathcal{I}_{v_1}^\infty} + \tau_{v_1}\mathcal{I}_{v_1}^\infty\right).
\]
We obtain for $k=1,2$,
\[
  \mathcal{I}_{v_k}^\infty=\frac{1}{\tau_{v_k}}W\left(-\mathscr{R}_{e,v_k}\exp\left(\mathcal{J}_{v_k}^\infty\right)\right)-\frac{\mathcal{J}_{v_k}^\infty}{\tau_{v_k}},
\]
% \[
%   \mathcal{I}_{v_1}^\infty=\frac{1}{\tau}
%   W\left(-\mathscr{R}_{e,v_1}\exp\left(-\mathscr{R}_0+\mathscr{R}_{e,v_2}e^{-\tau
%         \mathcal{I}_{v_2}^\infty}+\tau\mathcal{I}_{v_2}^\infty\right)\right) + \frac{\mathscr{R}_0-\mathscr{R}_{e,v_2}e^{-\tau
%         \mathcal{I}_{v_2}^\infty}-\tau\mathcal{I}_{v_2}^\infty}{\tau},
% \]
% and
% \[
%   \mathcal{I}_{v_2}^\infty=\frac{1}{\tau}
%   W\left(-\mathscr{R}_{e,v_2}\exp\left(-\mathscr{R}_0+\mathscr{R}_{e,v_1}e^{-\tau
%         \mathcal{I}_{v_1}^\infty}+\tau\mathcal{I}_{v_1}^\infty\right)\right) + \frac{\mathscr{R}_0-\mathscr{R}_{e,v_1}e^{-\tau
%         \mathcal{I}_{v_1}^\infty}-\tau\mathcal{I}_{v_1}^\infty}{\tau},
% \]
still with $W$ equal to $W_{-1}$ and $W_0$. Let $\iota^1_{-1,0}$ and $\iota^2_{-1,0}$ be
defined by
\bqq
 \iota_{-1,0}^{v_1}=\frac{W_{-1,0}\left(-\mathscr{A}_{v_1}\right)}{\tau_{v_1}}+\frac{\frac{\tau_{v_1}\eta_{v_2}}{\tau_{v_2}\eta_{v_1}}\left(\mathscr{R}_{0,v_2}-1-\log\left(\mathscr{R}_{e,v_2}\right)\right)}{\tau_{v_1}},
 \label{eq:iota1}
\eqq
% \bqq
%  \iota_{-1,0}^1=\frac{W_{-1,0}\left(-\mathscr{R}_{e,v_1} \exp\left(-\left(\log(\exp(-1)/\mathscr{R}_{e,v_2})+\mathscr{R}_0\right)\right)\right)}{\tau}+\frac{\log(\exp(-1)/\mathscr{R}_{e,v_2})+\mathscr{R}_0}{\tau},
%  \label{eq:iota1}
% \eqq
and
\bqq
 \iota_{-1,0}^{v_2}=\frac{W_{-1,0}\left(-\mathscr{A}_{v_2}\right)}{\tau_{v_2}}
 +\frac{\frac{\tau_{v_2}\eta_{v_1}}{\tau_{v_1}\eta_{v_2}}\left(\mathscr{R}_{0,v_1}-1-\log\left(\mathscr{R}_{e,v_1}\right)\right)}{\tau_{v_2}},
\label{eq:iota2}
\eqq
with $\mathscr{A}_{v_1}$ and $\mathscr{A}_{v_2}$ given by \eqref{eq:Av1} and \eqref{eq:Av2}.
% \bqq
%  \iota_{-1,0}^2=\frac{W_{-1,0}\left(-\mathscr{R}_{e,v_2} \exp\left(-\left(\log(\exp(-1)/\mathscr{R}_{e,v_1})+\mathscr{R}_0\right)\right)\right)}{\tau}+\frac{\log(\exp(-1)/\mathscr{R}_{e,v_1})+\mathscr{R}_0}{\tau}.
% \label{eq:iota2}
% \eqq
Then,
\[
  (\mathcal{I}_{v_1}^\infty,\mathcal{I}_{v_2}^\infty)\in [\min\left(\iota_{-1}^{v_1},\iota_{0}^{v_1}\right),\max\left(\iota_{-1}^{v_1},\iota_{0}^{v_1}\right)]\times[\min\left(\iota_{-1}^{v_2},\iota_{0}^{v_2}\right),\max\left(\iota_{-1}^{v_2},\iota_{0}^{v_2}\right)].
\]
We can show that $\min\left(\iota_{-1}^{v_k},\iota_{0}^{v_k}\right)<0$ for $k=1,2$. So,
we can reduce this domain since $\mathcal{I}_{v_k}^\infty>0$. So, we define the
domain $\omega_{\mathcal{I}}$
\[
  (\mathcal{I}_{v_1}^\infty,\mathcal{I}_{v_2}^\infty)\in \omega_{\mathcal{I}}:=[0,\max\left(\iota_{-1}^{v_1},\iota_{0}^{v_1}\right)]\times[0,\max\left(\iota_{-1}^{v_2},\iota_{0}^{v_2}\right)].
\]
As a consequence, we have proved the following lemma.

\begin{figure}[!htbp]
  \centering
  \begin{tabular}{cc}
    \includegraphics[height=.29\textheight]{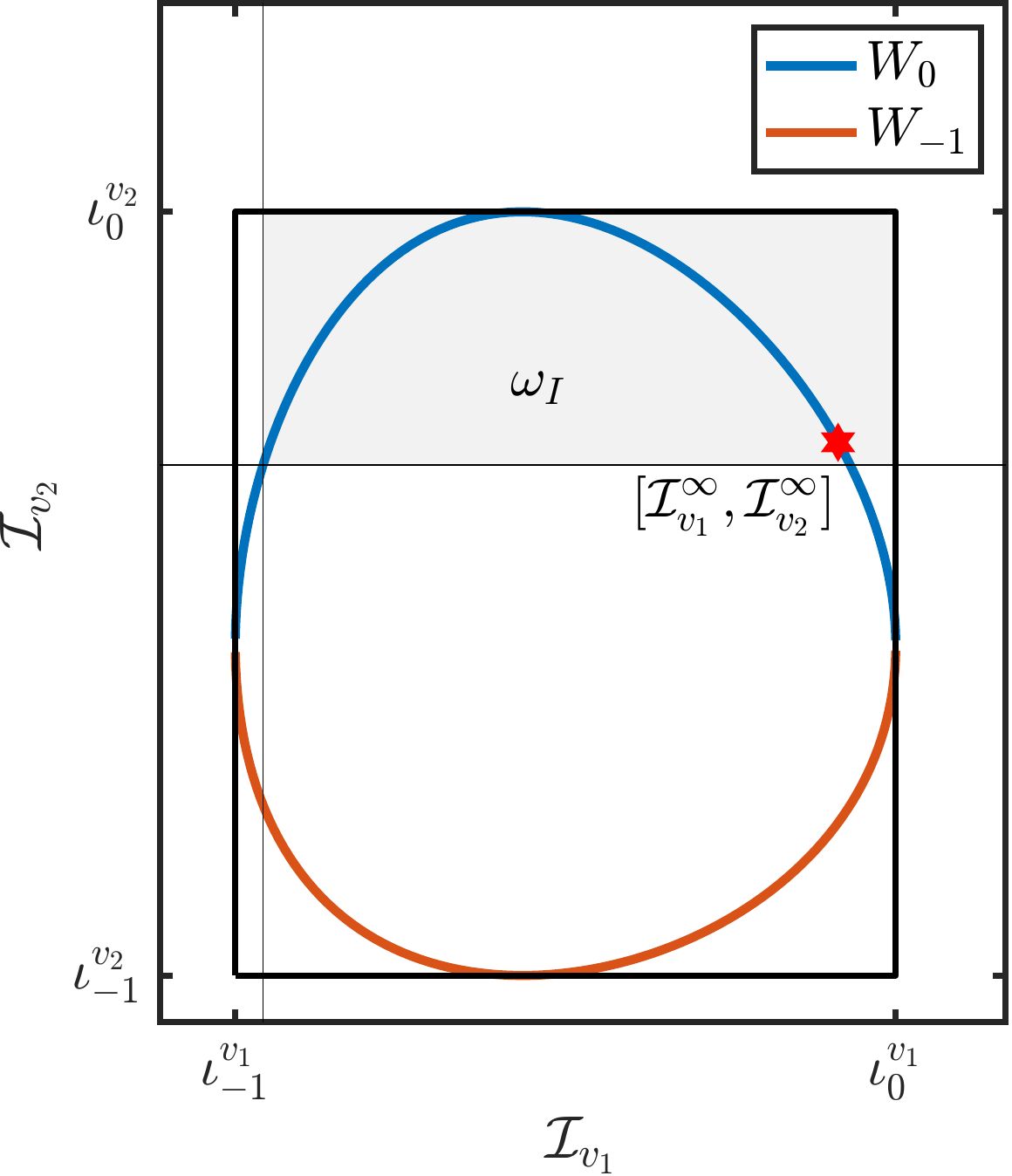} &
   \includegraphics[height=.29\textheight]{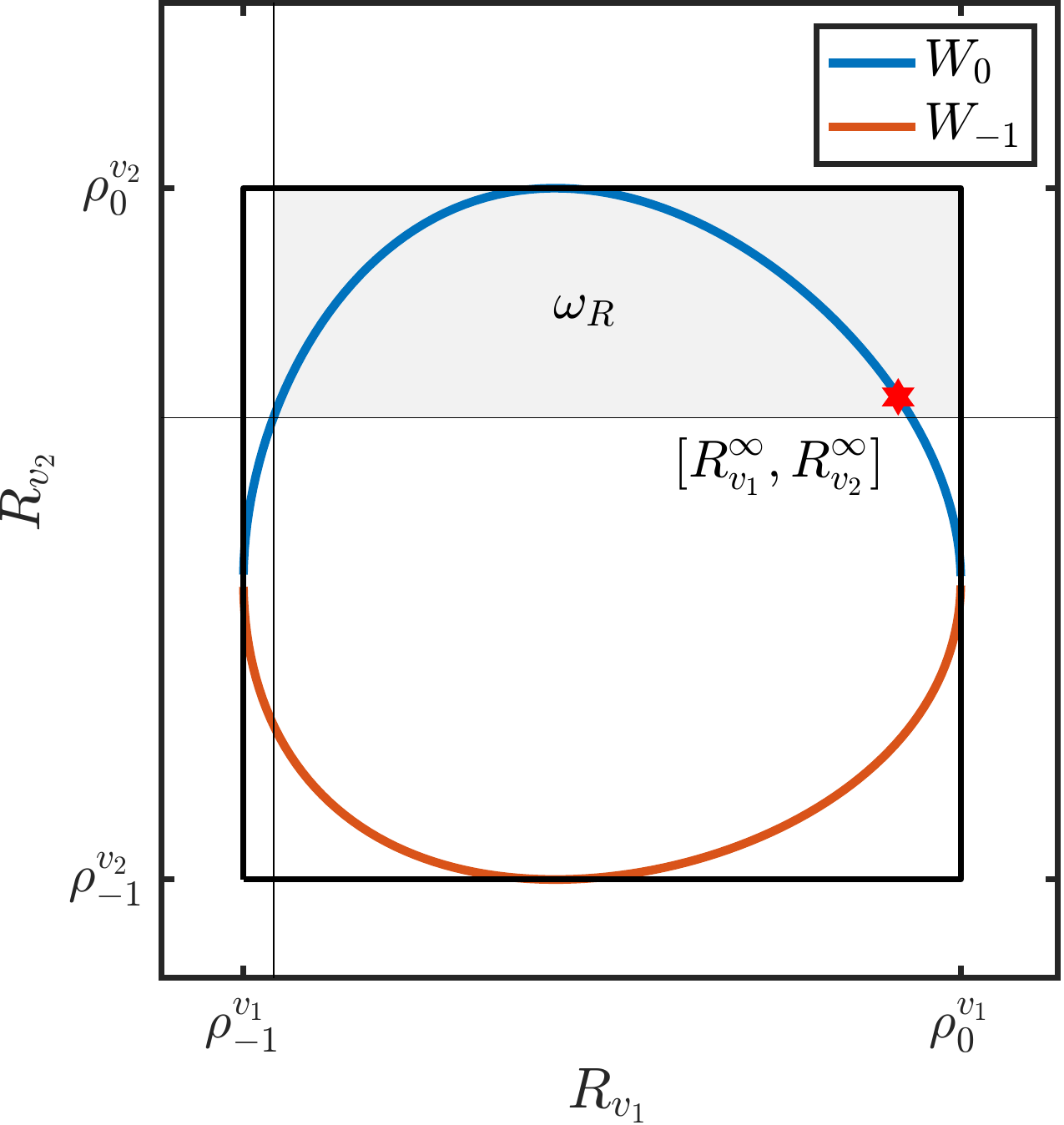}\\
  \end{tabular}
  \caption{Location of $\mathcal{I}_k^\infty$ and $R_{v_k}^\infty$, $k=1,2$,  and visualisation of the domains
    $\omega_{\mathcal{I}}$ (left), and domain $\omega_R$ (right). In both cases,  $(\mathcal{I}_1^\infty,\mathcal{I}_2^\infty)\in\omega_{\mathcal{I}}$ and $(R_1^\infty,R_2^\infty)\in\omega_{R}$ are represented by a red star. Values of the parameters and initial conditions are similar to Figure~\ref{fig:boxS}.}
  \label{fig:omega_I_R}
\end{figure}

\begin{lem}{Case of two vertices.}\label{lem_2vert}
Assume that $|\V|=2$ and $|\E|=1$. The final total population of infected
individuals at each vertex $\mathcal{I}_{v_k}^\infty$, $k=1,2$ can be expressed
as
\[
  \mathcal{I}_{v_k}^\infty=\frac{1}{\tau_{v_k}}W\left(-\mathscr{R}_{e,v_k}\exp\left(\mathcal{J}_{v_k}^\infty\right)\right)-\frac{\mathcal{J}_{v_k}^\infty}{\tau_{v_k}},
\]
with
\[
  \mathcal{J}_{v_k}^\infty = -\mathscr{R}_{0,v_k}+\frac{\tau_{v_k}\eta_{v_j}}{\tau_{v_j}\eta_{v_k}} \left(\mathscr{R}_{e,v_j}\exp^{-\tau_{v_j} \mathcal{I}_{v_j}^\infty} + \tau_{v_j}\mathcal{I}_{v_j}^\infty\right), \quad k\neq j \in\left\{1,2\right\},
\]
% \[
%   \mathcal{I}_{v_k}^\infty=\frac{1}{\tau}
%   W\left(-\mathscr{R}_{e,v_k}\exp\left(-\mathscr{R}_0+\mathscr{R}_{e,v_j}e^{-\tau
%         \mathcal{I}_{v_j}^\infty}+\tau\mathcal{I}_{v_j}^\infty\right)\right) + \frac{\mathscr{R}_0-\mathscr{R}_{e,v_j}e^{-\tau
%         \mathcal{I}_{v_j}^\infty}-\tau\mathcal{I}_{v_j}^\infty}{\tau}, \quad k\neq j \in\left\{1,2\right\},
% \]
where $\mathscr{R}_{0,v_k}:=M^0\, \tau_{v_k}/\eta_{v_k}$ and
$\mathscr{R}_{e,v_k}:=S_{v_k}^0\, {\tau_{v_k}}/{\eta_{v_k}}$, $k=1,2$. Furthermore, we have the sharp bound
\[
  (\mathcal{I}_{v_1}^\infty,\mathcal{I}_{v_2}^\infty)\in \omega_{\mathcal{I}}:=[0,\max\left(\iota_{-1}^{v_1},\iota_{0}^{v_1}\right)]\times[0,\max\left(\iota_{-1}^{v_2},\iota_{0}^{v_2}\right)],
\]
with $\iota_{-1,0}^{v_k}$, $k=1,2$ defined in \eqref{eq:iota1}-\eqref{eq:iota2}. See Figure~\ref{fig:omega_I_R} for an illustration.
\end{lem}

\begin{rmk}
As the solutions $R_{v_k}^\infty$, $k=1,2$, are simply given
by $R_{v_k}^\infty=\eta_{v_k}\mathcal{I}_{v_k}^\infty$, if we let
$\rho_{-1,0}^{v_k}:=\eta\iota_{-1,0}^{v_k}$ then we have
\[
  (R_{v_1}^\infty,R_{v_2}^\infty)\in \omega_{R}:=[0,\max\left(\rho_{-1}^{v_1},\rho_{0}^{v_1}\right)]\times[0,\max\left(\rho_{-1}^{v_2},\rho_{0}^{v_2}\right)].
\]
We represent on Figure~\ref{fig:omega_I_R} the domain $\omega_R$.
\end{rmk}

\section{A semi-implicit numerical scheme which preserves total mass}\label{secNumSc}

In this section, we propose a semi-implicit in time numerical scheme based on finite differences in space which has the property to preserve the {\it discrete} total mass.

\subsection{Notations}

For each $e\in\E$, we denote $\delta x_e>0$ the space discretization of each edge, and $J_e\in \N$ the number of points of the corresponding discretization. For each $i=1,\cdots,J_e$, the space grid on each edge is given by $x_i=(i-1)\delta x_e$ with $\ell_e=(J_e-1)\delta x_e$. And we let ${\bf J}:=\sum_{e\in\E} J_e \in\N$. Let $\delta t>0$ be the time discretization and denote $t_m=m\delta t$ for $m\geq 0$. 

 For a given function $u \in \mathscr{C}^{1,2}(\R_+\times \cG,\R^+)$, its space-time discretization is given by some sequence of vectors 
\bqs
u \sim (U^m)_{m\geq 0}, \text{ with } U^m=\left(U_1^m,\cdots,U_{\bf J}^m\right)^{\bf t} \in \R^{\bf J}.
\eqs
For each $e\in\E$, there exists an integer $j_e \in \N$ such that
\bqs
u_e(t_m, x_i) \sim U_{j_e+i}^m, \quad i=1,\cdots,J_e, \quad m\geq 0.
\eqs
We approximate the laplacian on each edge via finite differences. That is, for each $e\in\E$,
\bqs
\partial_x^2 u_e(t_m,x_i) \sim \frac{U_{j_e+i-1}^m-2U_{j_e+i}^m+U_{j_e+i+1}^m}{\delta x_e^2}, \quad i=2,\cdots,J_e-1, \quad m\geq 1,
\eqs
where we have only considered the interior points of the discretized domain. Let us now precise how we approximate the laplacian at a given vertex $v\in\V$ of the graph. So let $v\in\V$ such that there are $\delta_v$ edges incident to the vertex. We locally label $e\sim v = (e_1,\cdots,e_{\delta_v})$ all these incident edges. For each $v\in\V$, we introduce the map $\sigma_v:\left\{e_1,\cdots,e_{\delta_v}\right\}\rightarrow \left\{1,\cdots,{\bf J}\right\}$ such that $\sigma_v(e_k)$ corresponds to the global index of the grid discretization associated to the vertex $v$ on edge $e_k$. Finally, we denote by $n(\sigma_v(e_k))$ the global index of the nearest neighbor on edge $e_k$ to the vertex $v$. Note that either $n(\sigma_v(e_k))=\sigma_v(e_k)-1$ or $n(\sigma_v(e_k))=\sigma_v(e_k)+1$. To approximate the laplacian at a given vertex $v\in\V$ on edge $e_k$, we use the following formula
\bqs
\partial_x^2 u_{e_k}(t_m,v) \sim \frac{U_{\sigma_v(e_k)}^{*,m}-2U_{\sigma_v(e_k)}^m+U_{n(\sigma_v(e_k))}^m}{\delta x_{e_k}^2}:=\mathcal{Z}_{v,k}^{m}, \quad k=1,\cdots,\delta_v, \quad m\geq 1.
\eqs
The unknown $U_{\sigma_v(e_k)}^{*,m}$ can be expressed by discretization of the boundary condition as follows. For each $v\in\V$ with $e\sim v = (e_1,\cdots,e_{\delta_v})$, we approximate the normal derivative $\partial_n u_{e_k}(t_m,v)$ as
\bqs
\partial_n u_{e_k}(t_m,v) \sim \frac{U_{\sigma_v(e_k)}^{*,m}-U_{n(\sigma_v(e_k))}^m}{2\delta x_{e_k}}, \quad k=1,\cdots,\delta_v, \quad m\geq 1.
\eqs
Using \eqref{BCgraph}, and denoting $I_v^m$ the time approximation of $I_v(t_m)$, we obtain the following expression for $U_{\sigma_v(e_k)}^{*,m}$
\bqs
U_{\sigma_v(e_k)}^{*,m}=U_{n(\sigma_v(e_k))}^m-\frac{2\delta x_{e_k}}{d_{e_k}}\left( \alpha_{e_k}^vU_{\sigma_v(e_k)}^m + \sum_{l=1}^{\delta_v} (N_v)_{kl} U_{\sigma_v(e_l)}^m -\lambda_{e_k}^v I_v^m \right), \quad k=1,\cdots,\delta_v, \quad m\geq 1.
\eqs
As a consequence, we obtain that for each $k=1,\cdots,\delta_v$ and $m\geq1$
\bqs
\mathcal{Z}_{v,k}^{m}=\frac{2U_{n(\sigma_v(e_k))}^m-2U_{\sigma_v(e_k)}^m}{\delta x_{e_k}^2}-\frac{2}{d_{e_k}\delta x_{e_k}}\left( \alpha_{e_k}^vU_{\sigma_v(e_k)}^m + \sum_{l=1}^{\delta_v} (N_v)_{kl} U_{\sigma_v(e_l)}^m -\lambda_{e_k}^v I_v^m \right).
\eqs

\subsection{The semi-implicit numerical scheme}

We introduce the following scheme for each $m\geq 0$
\bqq
\left\{
\begin{split}
U^{m+1}_{j_e+i}&=U^m_{j_e+i}+ d_e\delta t \left( \frac{U_{j_e+i-1}^{m+1}-2U_{j_e+i}^{m+1}+U_{j_e+i+1}^{m+1}}{\delta x_e^2}\right), \quad i=2,\cdots,J_e-1, \quad e\in\E,\\
U_{\sigma_v(e_k)}^{m+1}&=U_{\sigma_v(e_k)}^m+d_{e_k}\delta t \mathcal{Z}_{v,k}^{m+1},\quad  k=1,\cdots,\delta_v,\quad v\in\V,\\
S_v^{m+1}&=S_v^m-\delta t \tau_v S_v^{m+1}I_v^m,\\
I_v^{m+1}&=I_v^m+\delta t \left(\tau_v S_v^{m+1}I_v^m-\eta_v I_v^{m+1}\right)+\delta t \left(\sum_{k=1}^{\delta_v} \alpha^v_{e_k} U^{m+1}_{\sigma_v(e_k)}-\overline{\lambda}_vI_v^{m+1}\right),\\
R_v^{m+1}&=R_v^m+\delta t\eta_v I_v^{m+1},
\end{split}
\right.
\label{numsch}
\eqq
initialized with $U^0\in \R^{\bf J}$ and some $(S_v^0,I_v^0,R_v^0)_{v\in\V}$. One can find similar semi-implicit discretization for the SIR part of the model in \cite{SE11}.

\paragraph{Well-posedness and positivity.} We prove that the numerical scheme defined through \eqref{numsch} is well defined and preserves positivity under some condition on $\delta t$. Indeed, we first remark that the equation for $S_v^{m+1}$ and $I_v^{m+1}$ in \eqref{numsch} can be used to obtain that
\begin{align*}
S_v^{m+1}&=\frac{S_v^m}{1+\delta t \tau_v I_v^m},\\
I_v^{m+1}&=\frac{I_v^m+\delta t \tau_v I_v^m(S_v^m+ I_v^m)}{\left(1+\delta t(\eta_v+\overline{\lambda}_v)\right)(1+\delta t \tau_v I_v^m)}+\frac{\delta t}{1+\delta t (\eta_v+\overline{\lambda}_v)}\sum_{l=1}^{\delta_v} \alpha^v_{e_l} U^{m+1}_{\sigma_v(e_l)},
\end{align*}
such that $\mathcal{Z}_{v,k}^{m+1}$ can be expressed only in terms of elements of $U^{m+1}$ as
\begin{align*}
\mathcal{Z}_{v,k}^{m+1}&=\frac{2U_{n(\sigma_v(e_k))}^{m+1}-2U_{\sigma_v(e_k)}^{m+1}}{\delta x_{e_k}^2}-\frac{2}{d_{e_k}\delta x_{e_k}}\left( \alpha_{e_k}^vU_{\sigma_v(e_k)}^{m+1} + \sum_{l=1}^{\delta_v} (N_v)_{kl} U_{\sigma_v(e_l)}^{m+1} -\frac{\delta t \lambda_{e_k}^v}{1+\delta t (\eta_v+\overline{\lambda}_v)}\sum_{l=1}^{\delta_v} \alpha^v_{e_l} U^{m+1}_{\sigma_v(e_l)}  \right)\\
&~~~+\frac{2\lambda_{e_k}^v}{d_{e_k}\delta x_{e_k}}\left(\frac{I_v^m+\delta t \tau_v I_v^m(S_v^m+ I_v^m)}{\left(1+\delta t(\eta_v+\overline{\lambda}_v)\right)(1+\delta t \tau_v I_v^m)}\right).
\end{align*}
As a consequence, there exists a matrix $\mathcal{A}\in \mathscr{M}_{\bf J}(\R)$ such that
\bqs
\left(\mathrm{I}_{\bf J}+\mathcal{A}\right)U^{m+1}=U^m+\mathcal{Y}^{m}, 
\eqs
where $\mathcal{Y}^{m}\in\R^{\bf J}$ is such that
\bqs
\mathcal{Y}^{m}_j=\left\{
\begin{array}{cl}
 \dfrac{2\delta t\lambda_{e_k}^v}{\delta x_{e_k}}\left(\dfrac{I_v^m+\delta t \tau_v I_v^m(S_v^m+ I_v^m)}{\left(1+\delta t(\eta_v+\overline{\lambda}_v)\right)(1+\delta t \tau_v I_v^m)}\right)&, \text{ if } j=\sigma_v(e_k),\\
 0&, \text{ otherwise}.
\end{array}
\right.
\eqs

\begin{lem}\label{lem_wellpnum} There exists a constant $C_0>0$, which only depends on the parameters of the system, such that if $0<\delta t <C_0$ then we have
\begin{itemize}
\item $\mathrm{I}_{\bf J}+\mathcal{A}$ is invertible;
\item if $N_v$ is symmetric for each $v\in\V$, then given $V\in\R^{\bf J}$ with $V\geq {\bf 0}$, the unique solution $U\in\R^{\bf J}$ of $(\mathrm{I}_{\bf J}+\mathcal{A})U=V$ also satisfies $U\geq {\bf 0}$.
\end{itemize}
\end{lem}

\begin{Proof}
Let $U\in\R^{\bf J}\neq {\bf 0} $ be such that $(\mathrm{I}_{\bf J}+\mathcal{A})U=0$. Without loss of generality, assume that $U_{j_0}=\max_{j=1,\cdots,\mathbf{J}} U_j>0$. If there exists $e\in\E$ such that $j_0=j_e+i_0$ for some $i_0\in\left\{2,\cdots,J_e-1\right\}$, then we have
\bqs
U_{j_0}+\frac{d_e \delta t}{\delta x_e^2}(2U_{j_0}-U_{j_0-1}-U_{j_0+1})=0,
\eqs
which is a contradiction by definition of $U_{j_0}$. Next if $j_0$ is such that there is $v\in\V$ and $k\in\{1,\cdots,\delta_v\}$ such that $j_0=\sigma_v(e_k)$, then we have
\bqs
U_{\sigma_v(e_k)}+\frac{2d_{e_k}\delta t}{\delta x_{e_k}^2}(U_{\sigma_v(e_k)}-U_{n(\sigma_v(e_k))})=-\frac{2\delta t}{\delta x_{e_k}}\left( \alpha_{e_k}^vU_{\sigma_v(e_k)} + \sum_{l=1}^{\delta_v} (N_v)_{kl} U_{\sigma_v(e_l)} -\frac{\delta t \lambda_{e_k}^v}{1+\delta t (\eta_v+\overline{\lambda}_v)}\sum_{l=1}^{\delta_v} \alpha^v_{e_l} U_{\sigma_v(e_l)}  \right).
\eqs
The left-hand side of the above equality is strictly positive and we claim that the right-hand side is negative. We use the fact that $(N_v)_{kl}=-\nu^v_{e_l,e_k}$ when $k\neq l$ and $(N_v)_{kk}=\sum_{j\neq k}\nu^v_{e_k,e_j}$
\begin{align*}
\sum_{l=1}^{\delta_v} (N_v)_{kl} U_{\sigma_v(e_l)}&=\left(\sum_{j\neq k}\nu^v_{e_k,e_j}\right) U_{\sigma_v(e_k)}-\sum_{l\neq k} \nu^v_{e_l,e_k} U_{\sigma_v(e_l)}\\
&=\left[\sum_{j\neq k}\nu^v_{e_k,e_j} - \sum_{l\neq k} \nu^v_{e_l,e_k} \right] U_{\sigma_v(e_k)}+\sum_{l\neq k} \nu^v_{e_l,e_k}\left(U_{\sigma_v(e_k)}- U_{\sigma_v(e_l)}\right).
\end{align*}
As a consequence, we deduce that
\begin{align*}
\mathcal{U}_k^v&:=\alpha_{e_k}^vU_{\sigma_v(e_k)} + \sum_{l=1}^{\delta_v} (N_v)_{kl} U_{\sigma_v(e_l)} -\frac{\delta t \lambda_{e_k}^v}{1+\delta t (\eta_v+\overline{\lambda}_v)}\sum_{l=1}^{\delta_v} \alpha^v_{e_l} U_{\sigma_v(e_l)}\\
&=\left[\alpha_{e_k}^v+\sum_{j\neq k}\nu^v_{e_k,e_j} - \sum_{l\neq k} \nu^v_{e_l,e_k} -\frac{\delta t \lambda_{e_k}^v}{1+\delta t(\eta_v+ \overline{\lambda}_v)}\sum_{l=1}^{\delta_v} \alpha^v_{e_l} \right] U_{\sigma_v(e_k)}+\sum_{l\neq k} \nu^v_{e_l,e_k}\left(U_{\sigma_v(e_k)}- U_{\sigma_v(e_l)}\right)\\
&~~~+\frac{\delta t \lambda_{e_k}^v}{1+\delta t(\eta_v+ \overline{\lambda}_v)}\sum_{l=1}^{\delta_v} \alpha^v_{e_l}\left(U_{\sigma_v(e_k)}- U_{\sigma_v(e_l)}\right).
\end{align*}
The last two terms are positive by definition of $U_{j_0}=U_{\sigma_v(e_k)}=\max_{j=1,\cdots,\mathbf{J}} U_j>0$. Now using Hypothesis~\ref{hypDD}, we have that
\bqs
\alpha_{e_k}^v+\sum_{j\neq k}\nu^v_{e_k,e_j} - \sum_{l\neq k} \nu^v_{e_l,e_k}>0,
\eqs
such that the term in bracket is positive provided that
\bqs
\frac{\delta t \lambda_{e_k}^v}{1+\delta t(\eta_v+ \overline{\lambda}_v)}\sum_{l=1}^{\delta_v} \alpha^v_{e_l}<\alpha_{e_k}^v+\sum_{j\neq k}\nu^v_{e_k,e_j} - \sum_{l\neq k} \nu^v_{e_l,e_k},
\eqs
or equivalently
\bqs
\delta t \left[ \lambda_{e_k}^v\sum_{l=1}^{\delta_v} \alpha^v_{e_l}-(\eta_v+\overline{\lambda}_v) \left(\alpha_{e_k}^v+\sum_{j\neq k}\nu^v_{e_k,e_j} - \sum_{l\neq k} \nu^v_{e_l,e_k} \right)  \right]<\alpha_{e_k}^v+\sum_{j\neq k}\nu^v_{e_k,e_j} - \sum_{l\neq k} \nu^v_{e_l,e_k}.
\eqs
As a consequence, we impose that
\bqs
0<\delta t < \underset{v\in\V}{\min}~  \underset{k=1,\cdots,\delta_v}{\min} \frac{\alpha_{e_k}^v+\sum_{j\neq k}\nu^v_{e_k,e_j} - \sum_{l\neq k} \nu^v_{e_l,e_k}}{\left[ \lambda_{e_k}^v\sum_{l=1}^{\delta_v} \alpha^v_{e_l}-(\eta_v+\overline{\lambda}_v) \left(\alpha_{e_k}^v+\sum_{j\neq k}\nu^v_{e_k,e_j} - \sum_{l\neq k} \nu^v_{e_l,e_k} \right)  \right]_+},
\eqs
where it is understood that when the positive part is zero there is no condition on $\delta t$. And we have reached a contradiction since
\bqs
0<
U_{\sigma_v(e_k)}+\frac{2d_{e_k}\delta t}{\delta x_{e_k}^2}(U_{\sigma_v(e_k)}-U_{n(\sigma_v(e_k))})=-\frac{2\delta t}{\delta x_{e_k}}\mathcal{U}_k^v<0.
\eqs
This shows that $\mathrm{I}_{\bf J}+\mathcal{A}$ is invertible.

Next let $U\in\R^{\bf J}$ be the unique solution of $(\mathrm{I}_{\bf J}+\mathcal{A})U=V$ with $V\geq {\bf 0}$. We denote by $U^- \in\R^{\bf J}$ the vector with components given by
\bqs
U^-_j = \min(0,U_j), \quad j=1,\cdots,{\bf J}.
\eqs
Our aim is to evaluate $\langle (\mathrm{I}_{\bf J}+\mathcal{A})U, U^-\rangle_{\bf J}$ where $\langle \cdot,\cdot\rangle_{\bf J}$ is the following  scalar product on $\R^{\bf J}$:
\bqs
\langle U,V \rangle_{\bf J}:=\sum_{e\in\E}  \left( \sum_{i=2}^{J_e-1}  U_{j_e+i}V_{j_e+i}\right)+\frac{1}{2}\sum_{v\in\V}  \left(\sum_{k=1}^{\delta_v}U_{\sigma_v(e_k)}V_{\sigma_v(e_k)}\right).
\eqs
We divide $\langle (\mathrm{I}_{\bf J}+\mathcal{A})U, U^-\rangle_{\bf J}$ into three parts:
\bqs
\langle (\mathrm{I}_{\bf J}+\mathcal{A})U, U^-\rangle_{\bf J} = \mathcal{Q}_1+\mathcal{Q}_2+\mathcal{Q}_3,
\eqs
where
\begin{align*}
 \mathcal{Q}_1&:=\sum_{e\in\E} \sum_{i=2}^{J_e-1} \left(U_{j_e+i}+\frac{d_e \delta t}{\delta x_e^2}\left(2U_{j_e+i}-U_{j_e+i-1}-U_{j_e+i+1}\right) \right)U^-_{j_e+i},\\
  \mathcal{Q}_2&:=\frac{1}{2}\sum_{v\in\V} \sum_{k=1}^{\delta_v} \left(U_{\sigma_v(e_k)}+\frac{2d_{e_k}\delta t}{\delta x_{e_k}^2} \left(U_{\sigma_v(e_k)}-U_{n(\sigma_v(e_k))}\right)\right) U_{\sigma_v(e_k)}^-,\\
 \mathcal{Q}_3&:=\delta t\sum_{v\in\V} \sum_{k=1}^{\delta_v} \frac{1}{\delta x_{e_k}}\left( \alpha_{e_k}^vU_{\sigma_v(e_k)} + \sum_{l=1}^{\delta_v} (N_v)_{kl} U_{\sigma_v(e_l)} -\frac{\delta t \lambda_{e_k}^v}{1+\delta t (\eta_v+\overline{\lambda}_v)}\sum_{l=1}^{\delta_v} \alpha^v_{e_l} U_{\sigma_v(e_l)}  \right) U_{\sigma_v(e_k)}^-.
\end{align*}
The first and second terms are handled as follows
\bqs
 \mathcal{Q}_1+\mathcal{Q}_2=\langle U, U^-\rangle_{\bf J}+\sum_{e\in\E} \frac{d_e \delta t}{\delta x_e^2}\sum_{i=1}^{J_e-1}\left(U_{j_e+i+1}-U_{j_e+i}\right)\left( U^-_{j_e+i+1}-U^-_{j_e+i}\right)\geq 0.
\eqs
For the third term $\mathcal{Q}_3$, if we further assume that $N_v$ is symmetric, then the matrix $K_v=A_v+N_v$ is symmetric positive definite, and thus for each $v\in\V$ there exists some $\beta_v>0$ such that
\bqs
\sum_{k=1}^{\delta_v} \frac{1}{\delta x_{e_k}}\left( \alpha_{e_k}^vU_{\sigma_v(e_k)} + \sum_{l=1}^{\delta_v} (N_v)_{kl} U_{\sigma_v(e_l)}  \right) U_{\sigma_v(e_k)}^- \geq \beta_v \sum_{k=1}^{\delta_v}  \frac{1}{\delta x_{e_k}}U_{\sigma_v(e_k)} U_{\sigma_v(e_k)}^-,
\eqs
while there exists $\omega_v>0$ such that
\bqs
\left( \sum_{k=1}^{\delta_v}  \frac{\lambda_{e_k}}{\delta x_{e_k}}U_{\sigma_v(e_k)}^-\right)\left(\sum_{l=1}^{\delta_v} \alpha^v_{e_l} U_{\sigma_v(e_l)}  \right)\leq \omega_v \sum_{k=1}^{\delta_v}  \frac{1}{\delta x_{e_k}}U_{\sigma_v(e_k)} U_{\sigma_v(e_k)}^-.
\eqs
And thus, we get an estimate for $\mathcal{Q}_3$ of the form
\bqs
\mathcal{Q}_3\geq \delta t \sum_{v\in\V}\left[ \left( \beta_v-\frac{\delta t \omega_v}{1+\delta t (\eta_v+ \overline{\lambda}_v)}\right)\sum_{k=1}^{\delta_v}  \frac{1}{\delta x_{e_k}}U_{\sigma_v(e_k)} U_{\sigma_v(e_k)}^-\right],
\eqs
which is positive provided that $\delta t$ is small enough. As a consequence, we have proved that
\bqs
0\leq \langle (\mathrm{I}_{\bf J}+\mathcal{A})U, U^-\rangle_{\bf J} = \langle V, U^-\rangle_{\bf J}\leq 0,
\eqs
which implies that $U^-={\bf 0}$ and thus $U\geq {\bf 0}$.
\end{Proof}

The previous lemma demonstrates the well-posedness of our numerical scheme \eqref{numsch}. It also ensures that if we start with positive initial conditions $U^0\geq {\bf 0}$ and $S_v^0>0$, $I_v^0\geq0$ with $\sum_{v\in\V}I_v^0>0$ and $R_v^0\geq 0$, then for all $m\geq 1$ we also have that $U^m\geq {\bf 0}$, $S_v^m>0$, $I_v^{m}\geq0$ and $R_v^m\geq 0$, provided $\delta t>0$ is small enough and $N_v$ is symmetric for each $v\in\V$.

\paragraph{Preservation of total discrete mass.} For any $U\in\R^{\bf J}$, we define the following quantity 
\bqs
\mathrm{trap}_{\bf J}(U):=\sum_{e\in\E} \delta x_e \left( \sum_{i=2}^{J_e-1}  U_{j_e+i}\right)+\frac{1}{2}\sum_{v\in\V} \delta x_{e_k} \left(\sum_{k=1}^{\delta_v}U_{\sigma_v(e_k)}\right).
\eqs
The expression $\mathrm{trap}_{\bf J}(U)$ is simply the trapezoidal rule applied to the elements of $U$ adapted to our graph $\cG$.
From \eqref{numsch}, we get that
\begin{align*}
\mathrm{trap}_{\bf J}(U^{m+1})&=\mathrm{trap}_{\bf J}(U^m)+\delta t \sum_{e\in\E}  \frac{d_e}{\delta x_e} \left( \sum_{i=2}^{J_e-1}  \left( U_{j_e+i-1}^{m+1}-2U_{j_e+i}^{m+1}+U_{j_e+i+1}^{m+1}\right) \right)\\
&~~~+\frac{\delta t}{2}\sum_{v\in\V} \delta x_{e_k} d_{e_k} \left(\sum_{k=1}^{\delta_v}\mathcal{Z}^{m+1}_{v,k}\right).
\end{align*}
Upon denoting $\mathcal{Z}_{v,k}^{1,m}$ the following quantity
\bqs
\mathcal{Z}_{v,k}^{1,m}:=\frac{2U_{n(\sigma_v(e_k))}^m-2U_{\sigma_v(e_k)}^m}{\delta x_{e_k}^2}
\eqs
we get that
\bqs
\frac{1}{2}\sum_{v\in\V} \delta x_{e_k} d_{e_k} \left(\sum_{k=1}^{\delta_v}\mathcal{Z}^{1,m+1}_{v,k}\right)=\sum_{v\in\V} \frac{ d_{e_k}}{\delta x_{e_k}} \left(\sum_{k=1}^{\delta_v}U_{n(\sigma_v(e_k))}^{m+1}-U_{\sigma_v(e_k)}^{m+1}\right).
\eqs
Next, we observe that
\bqs
\sum_{e\in\E}  \frac{d_e}{\delta x_e} \left( \sum_{i=2}^{J_e-1}  \left( U_{j_e+i-1}^{m+1}-2U_{j_e+i}^{m+1}+U_{j_e+i+1}^{m+1}\right) \right)+\sum_{v\in\V} \frac{ d_{e_k}}{\delta x_{e_k}} \left(\sum_{k=1}^{\delta_v}U_{n(\sigma_v(e_k))}^{m+1}-U_{\sigma_v(e_k)}^{m+1}\right)=0,
\eqs
where the cancellation comes from the specific structure of the discretized laplacian through finite differences. As a consequence, we have that
\bqs
\mathrm{trap}_{\bf J}(U^{m+1})=\mathrm{trap}_{\bf J}(U^m)-\delta t \sum_{v\in\V}  \sum_{k=1}^{\delta_v}\left( \alpha_{e_k}^vU_{\sigma_v(e_k)}^{m+1} + \sum_{l=1}^{\delta_v} (N_v)_{kl} U_{\sigma_v(e_l)}^{m+1} -\lambda_{e_k}^v I_v^{m+1} \right).
\eqs
We also have that
\bqs
\sum_{k=1}^{\delta_v}  \sum_{l=1}^{\delta_v} (N_v)_{kl} U_{\sigma_v(e_l)}^{m+1} =  \sum_{l=1}^{\delta_v} \left(\sum_{k=1}^{\delta_v}  (N_v)_{kl} \right)U_{\sigma_v(e_l)}^{m+1}=0,
\eqs
as the sum over the lines of $N_v$ vanishes. And thus we get
\bqs
\mathrm{trap}_{\bf J}(U^{m+1})=\mathrm{trap}_{\bf J}(U^m)-\delta t \sum_{v\in\V} \left(  \sum_{k=1}^{\delta_v}\alpha_{e_k}^vU_{\sigma_v(e_k)}^{m+1} -\overline{\lambda}_v I_v^{m+1} \right).
\eqs
On the other hand, from \eqref{numsch} we also have
\bqs
\sum_{v\in\V}\left( S_v^{m+1}+I_v^{m+1}+R_v^{m+1}\right)= \sum_{v\in\V}\left( S_v^{m}+I_v^{m}+R_v^{m}\right)+\delta t \sum_{v\in\V} \left(  \sum_{k=1}^{\delta_v}\alpha_{e_k}^vU_{\sigma_v(e_k)}^{m+1} -\overline{\lambda}_v I_v^{m+1} \right).
\eqs
As a conclusion, we have proved the following result.

 \begin{lem}\label{lem_discmass}
 Let $(U^m,S^m_v,I^m_v,R^m_v)$ a solution of  \eqref{numsch}, then we have for each $m\geq 0$
\bqs
\mathrm{trap}_{\bf J}(U^{m+1})+\sum_{v\in\V}\left( S_v^{m+1}+I_v^{m+1}+R_v^{m+1}\right)=\mathrm{trap}_{\bf J}(U^m)+\sum_{v\in\V}\left( S_v^{m}+I_v^{m}+R_v^{m}\right).
\eqs
\end{lem}

This is the discrete conter part of conservation of mass for the continuous model.

Now, combining Lemma~\ref{lem_wellpnum}-\ref{lem_discmass}, we have proved the following theorem.

\begin{thm}\label{thmNumSc}
There exists a constant $C_0>0$, which only depends on the parameters of the system, such that if $0<\delta t <C_0$, then the numerical scheme \eqref{numsch} defines a unique  sequence $(U^m,S^m_v,I^m_v,R^m_v)_{m\geq0}$. If we further assume that $N_v$ is symmetric for each $v\in\V$, then the numerical scheme \eqref{numsch} preserves the positivity of the initial condition. Finally, for each solution of \eqref{numsch}, the total discrete mass is preserved, namely for each $m\geq 0$, we have
\bqs
\mathrm{trap}_{\bf J}(U^{m+1})+\sum_{v\in\V}\left( S_v^{m+1}+I_v^{m+1}+R_v^{m+1}\right)=\mathrm{trap}_{\bf J}(U^m)+\sum_{v\in\V}\left( S_v^{m}+I_v^{m}+R_v^{m}\right).
\eqs
\end{thm}

\section{Numerical results for a selection of graphs}\label{secSim}

In the present section, we illustrate our theoretical results with a collection of numerical simulations for various types of graphs. Throughout this section the time discretization is set to $\delta t=0.01$ while the space discretization to $\delta_{x_e}=0.01$ for each $e\in\E$.

\subsection{Case of 2 vertices and 1 edge}
\begin{center}
\begin{tikzpicture}
    \coordinate (A1) at (0,0);
    \coordinate (A2) at (2,0);

    \node[below] at (A1) {$v_1$};
    \node[below] at (A2) {$v_2$};

    \node at (A1) {$\bullet$};
    \node at (A2) {$\bullet$};

    \draw (A1) -- (A2);
\end{tikzpicture}
\end{center}
We first consider the case where $c_{\V}=2$ and $c_{\E}=1$, where $c_{\E}$
denotes the cardinal of $\E$. In this setting, we recall that our model reads as follows
\bqs
\partial_t u(t,x) = d \partial_x^2 u(t,x), \quad t>0, \quad x\in (0,\ell),
\eqs
with boundary conditions
\bqs
\left\{
\begin{split}
-d\partial_x u(t,0)+\alpha_1 u(t,0) &=\lambda_1 I_1(t),\\
d\partial_x u(t,\ell)+\alpha_2 u(t,\ell) &=\lambda_2 I_2(t),
\end{split}
\right. \quad t>0,
\eqs
where $(S_i(t),I_i(t),R_i(t))$, for $i=1,2$, solution of
\bqs
\left\{
\begin{split}
S'_i(t)&=-\tau_i S_i(t)I_i(t),\\
I_i'(t)&=\tau_i S_i(t)I_i(t)-\eta_i I_v(t)+\alpha_iu(t,v_i)-\lambda_iI_i(t),\\
R_i'(t)&=\eta_i I_i(t),
\end{split}
\right.  \quad t>0,
\eqs
where $v_1=0$ and $v_2=\ell$. This system is complemented by some initial condition $(u^0,S_i^0,I_i^0,R_i^0)$ with $S_i>0$, $I_1^0+I_2^0>0$, $R_i^0=0$ and $u^0\geq 0$ such that the boundary condition is satisfied initially. Finally, we normalize the total mass as follows
\bqs
M^0 = \int_0^\ell u^0(x)\md x +\sum_{i=1}^2\left(S_i^0+I_i^0 \right)=1.
\eqs

\begin{figure}[!t]
\centering
\begin{tabular}{ccc}
\includegraphics[height=0.16\textheight]{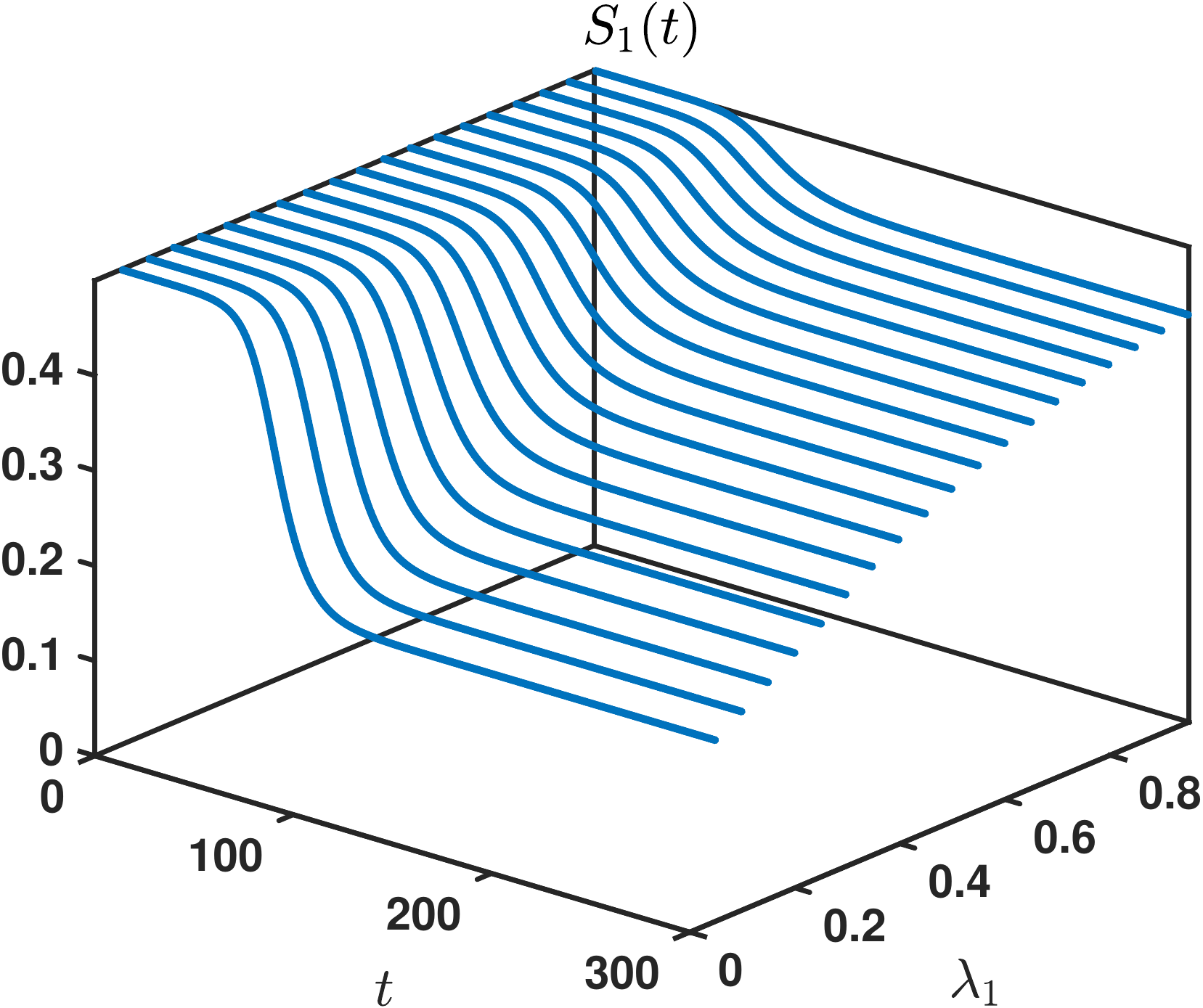}&
\includegraphics[height=0.16\textheight]{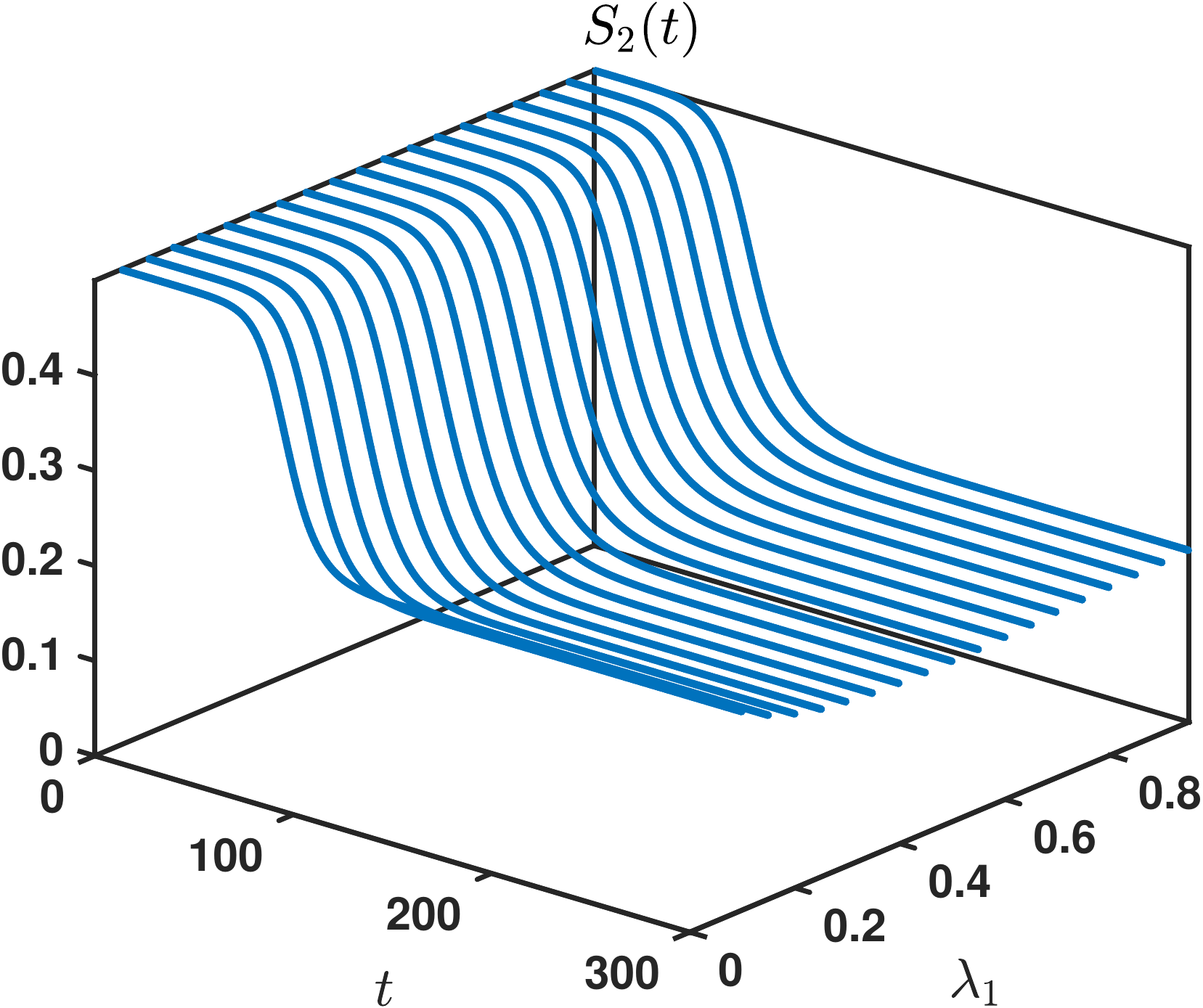}&
\includegraphics[height=0.16\textheight]{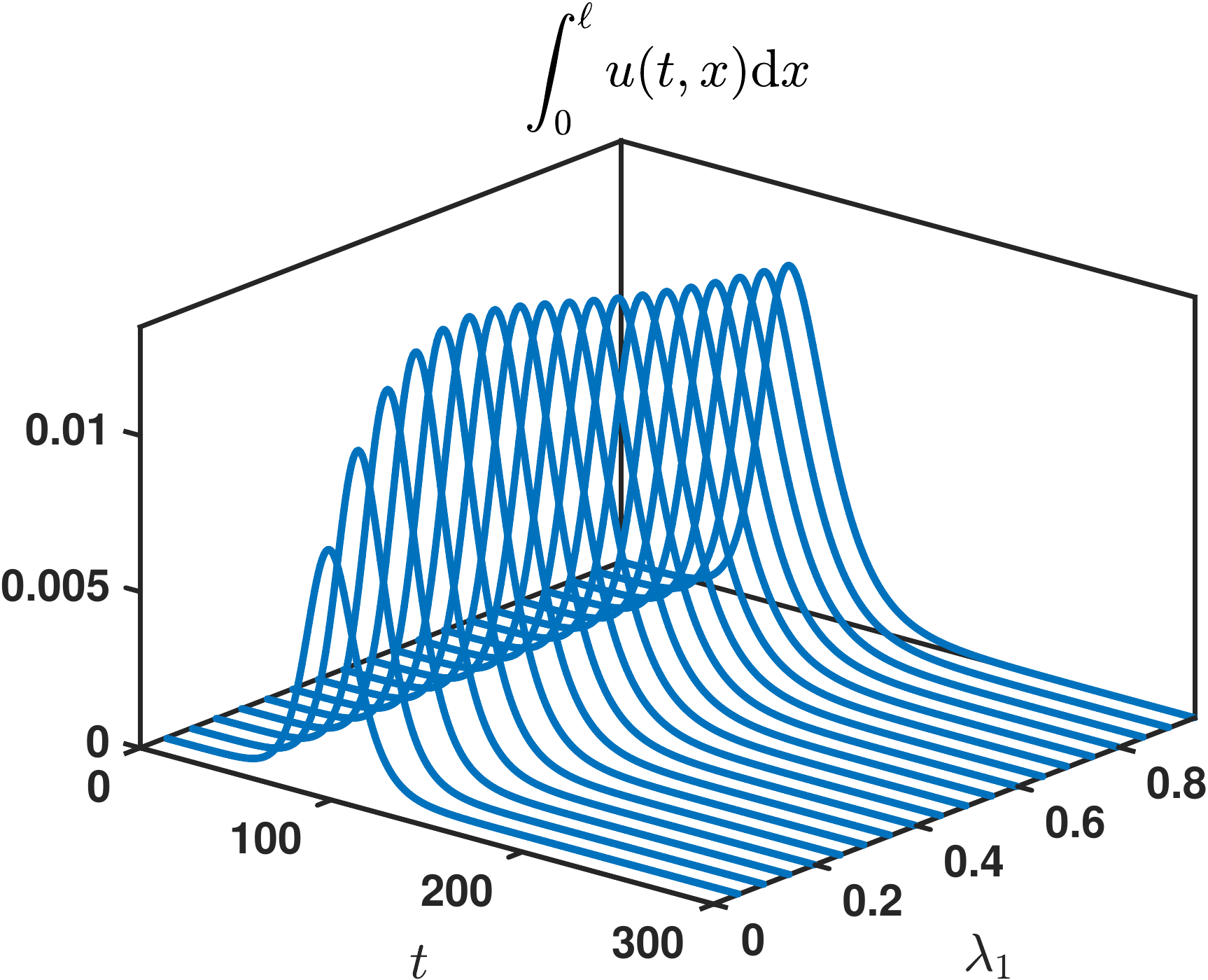}\\
\includegraphics[height=0.16\textheight]{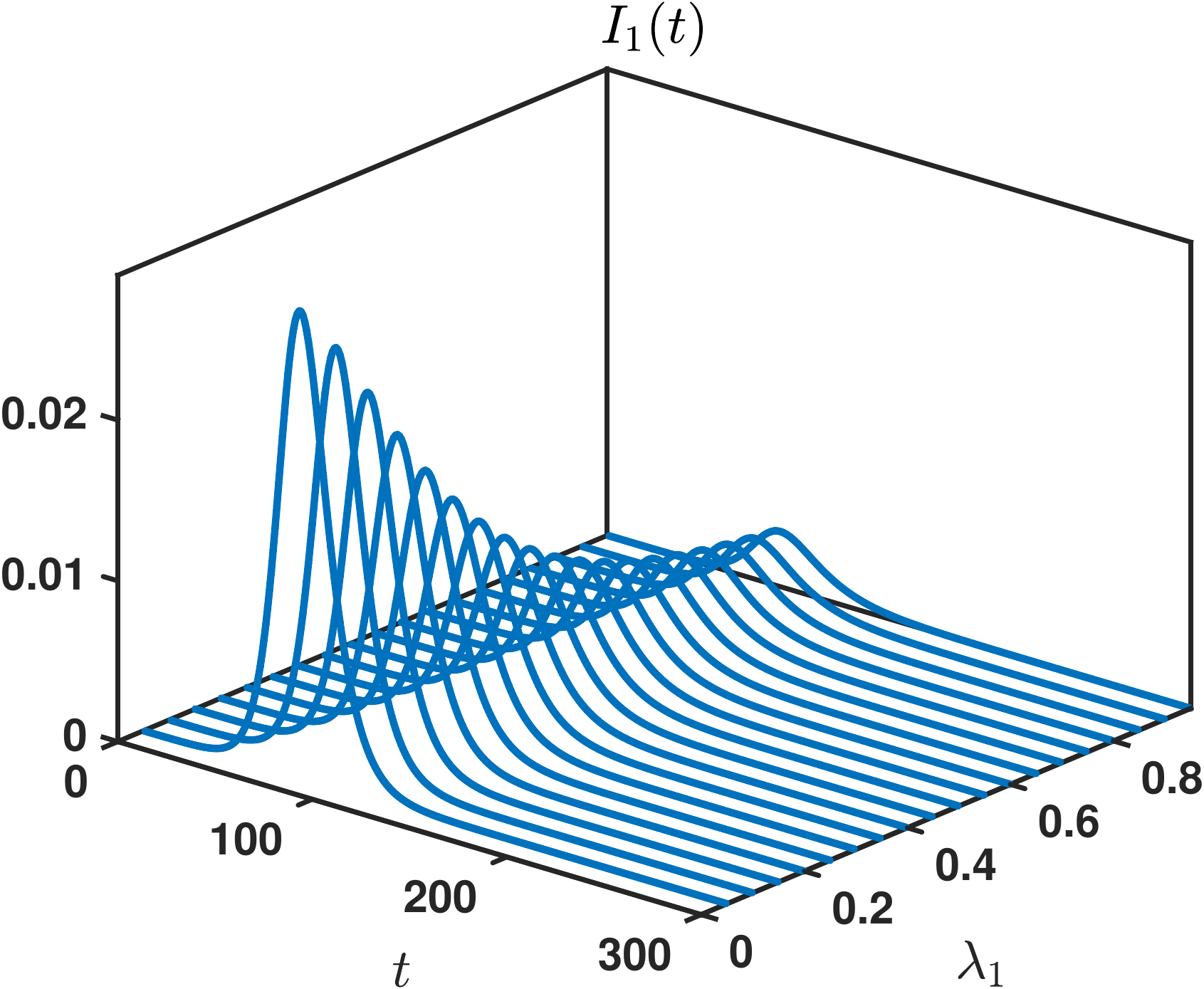}&
\includegraphics[height=0.16\textheight]{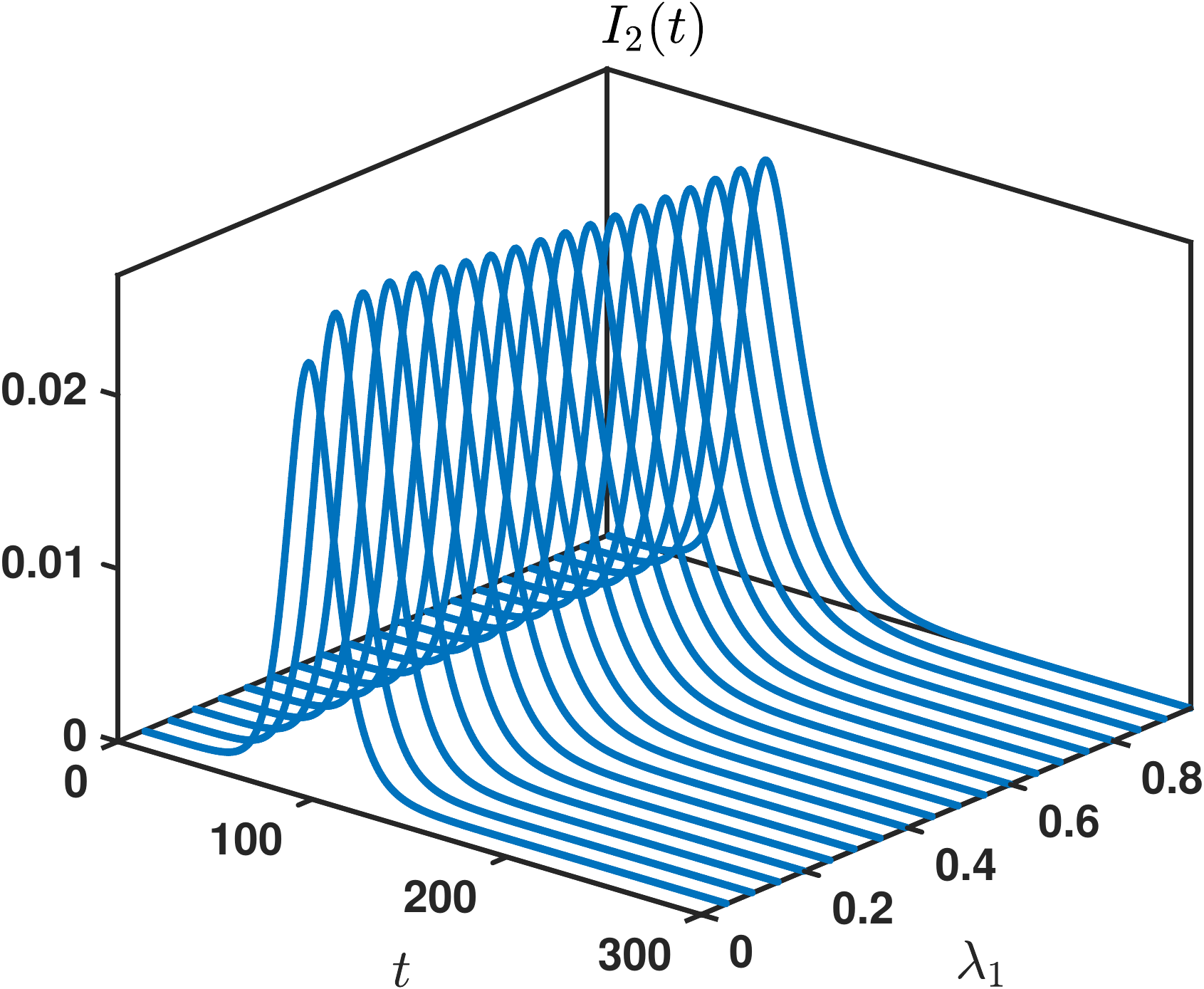}&
\includegraphics[height=0.16\textheight]{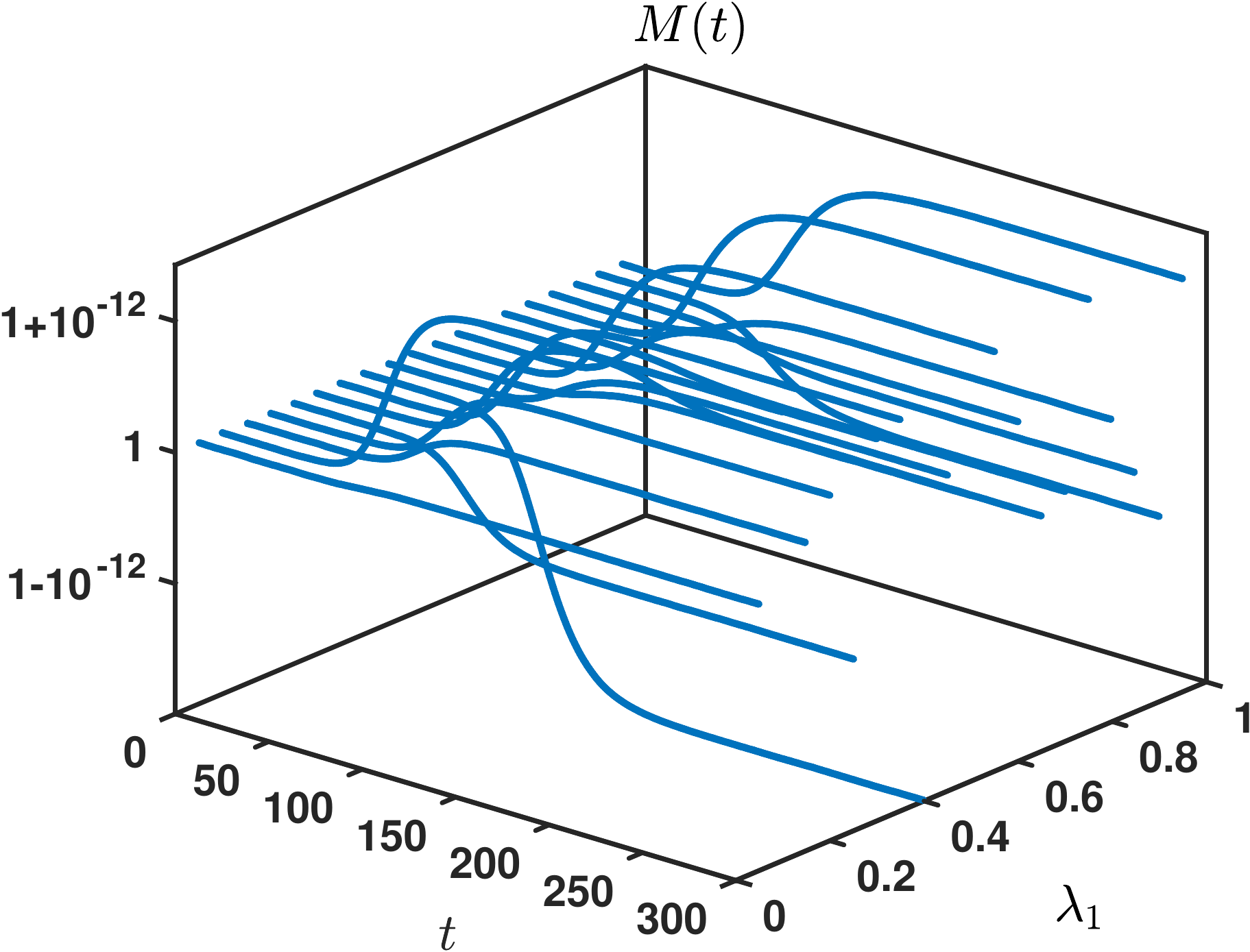}\\
\end{tabular}
\caption{Profiles of the solutions $(S_i(t),I_i(t))$ together with the total population on the edge $\int_0^\ell u(t,x)\md x$ and the total mass of the system $M(t)$ as the parameter $\lambda_1$ is varied from $0.05$ to $0.95$. All other parameters are fixed and set to $d=\ell=1$, $\lambda_2=1/10$, $\alpha_1=\alpha_2=1/4$, and $\tau_1=\tau_2=1$ with $\eta_1=\eta_2=1/3$. For the initial condition we have $(S_0,I_0)=(1/2,10^{-6})$.}
\label{fig:2c1r}
\end{figure}

\begin{figure}[!t]
  \centering
  \begin{tabular}{c}
\begin{tikzpicture}[scale=.8]
  % % Dimensions du repere
  % \def\xmin{-1} \def\xmax{6} \def\ymin{-1} \def\ymax{5}
  % % Grilles
  % \draw [step=0.1,gray!30,very thin]  (\xmin,\ymin) grid (\xmax,\ymax);
  % \draw [step=1,gray,thin] (\xmin,\ymin) grid (\xmax,\ymax);
  % \draw [step=5,thin,color=blue!50]  (\xmin,\ymin) grid (\xmax,\ymax); 
  % \foreach \x in {\xmin,...,\xmax}
  %     \node[color=gray!80] at (\x,-0.25) {\x};
  % \foreach \y in {\ymin,...,\ymax}
  %     \node[color=gray!80] at (-0.25,\y) {\y};
   \node[above right] at (0,0){\includegraphics[height=0.165\textheight]{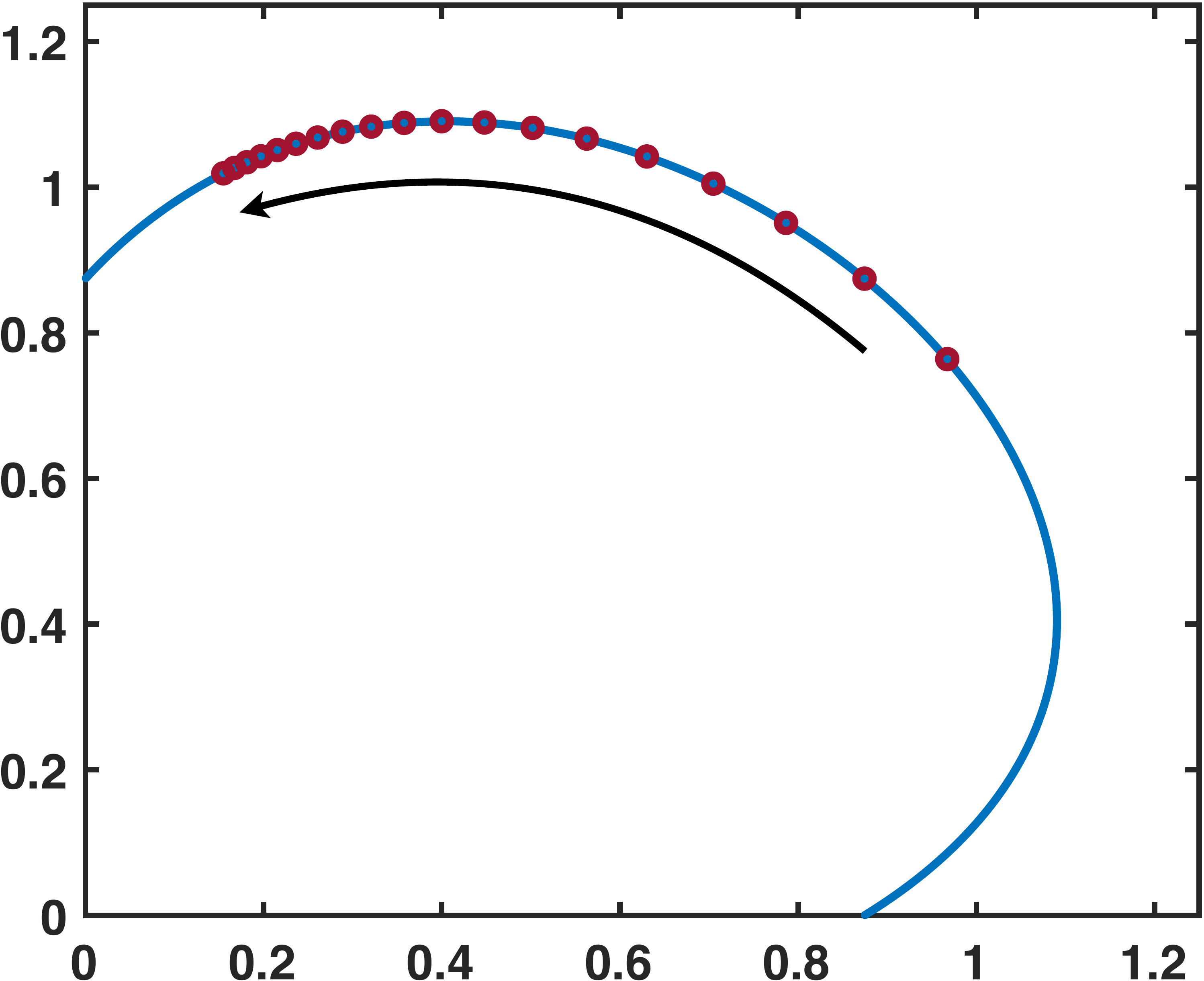}};
   \node[font=\fontsize{9pt}{0}\selectfont] at (2.9,-0.2) {$\mathcal{I}_1^\infty$};
   \node[font=\fontsize{9pt}{0}\selectfont] at (-0.2,2.5) {$\mathcal{I}_2^\infty$};
   \node[font=\fontsize{9pt}{0}\selectfont] at (2.4,3.5) {$\lambda_1$};
   % \node[right,font=\fontsize{4pt}{0}\selectfont] at (4.5,3) {$\lambda_1=0.05$};
   % \node[left,font=\fontsize{4pt}{0}\selectfont] at (1.5,4.1) {$\lambda_1=0.95$};
\end{tikzpicture}
\begin{tikzpicture}[scale=.8]
   \node[above right] at (0,0){\includegraphics[height=0.165\textheight]{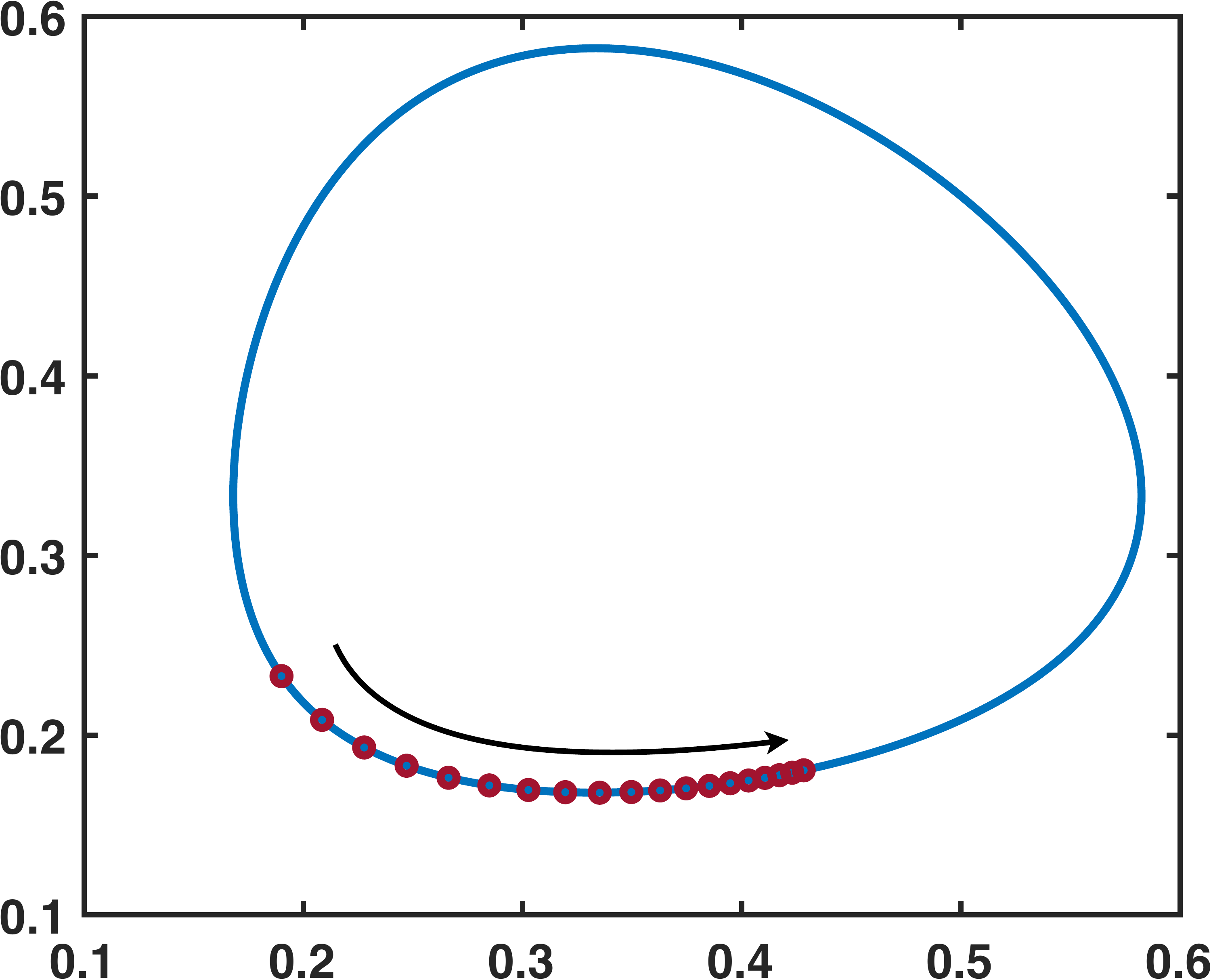}};
   \node[font=\fontsize{9pt}{0}\selectfont] at (2.9,-0.2) {$S_1^\infty$};
   \node[font=\fontsize{9pt}{0}\selectfont] at (-0.2,2.5) {$S_2^\infty$};
   \node[font=\fontsize{9pt}{0}\selectfont] at (2.6,1.6) {$\lambda_1$};
\end{tikzpicture} 
\begin{tikzpicture}[scale=.8]
   \node[above right] at (0,0){\includegraphics[height=0.165\textheight]{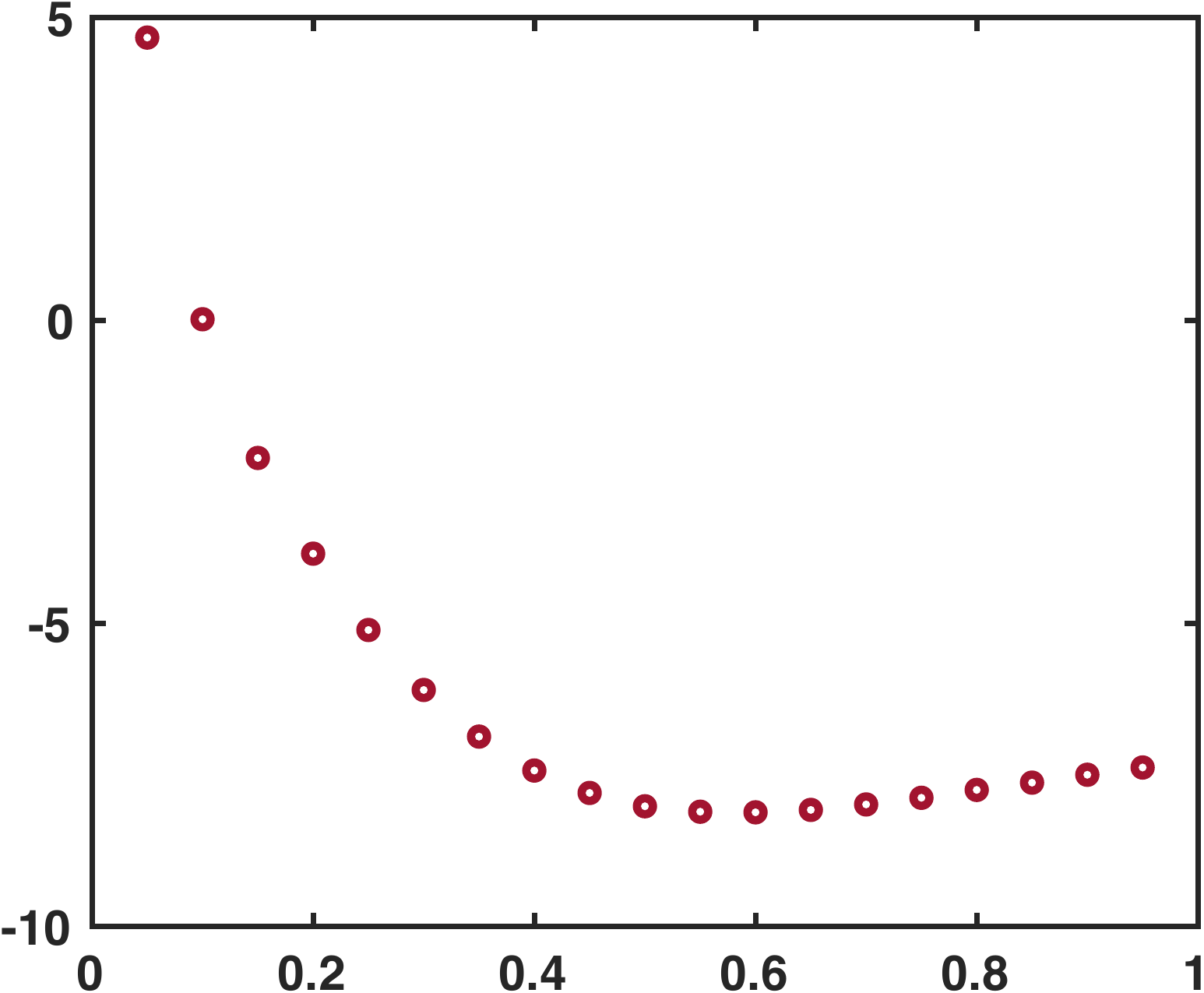}};
   \node[font=\fontsize{9pt}{0}\selectfont] at (2.9,-0.2) {$\lambda_1$};
   \node[font=\fontsize{9pt}{0}\selectfont] at (0,2.5) {$\Delta \mathcal{T}$};
\end{tikzpicture}\\[-2mm]
\begin{tikzpicture}[scale=.8]
   \node[above right] at (0,0){\includegraphics[height=0.165\textheight]{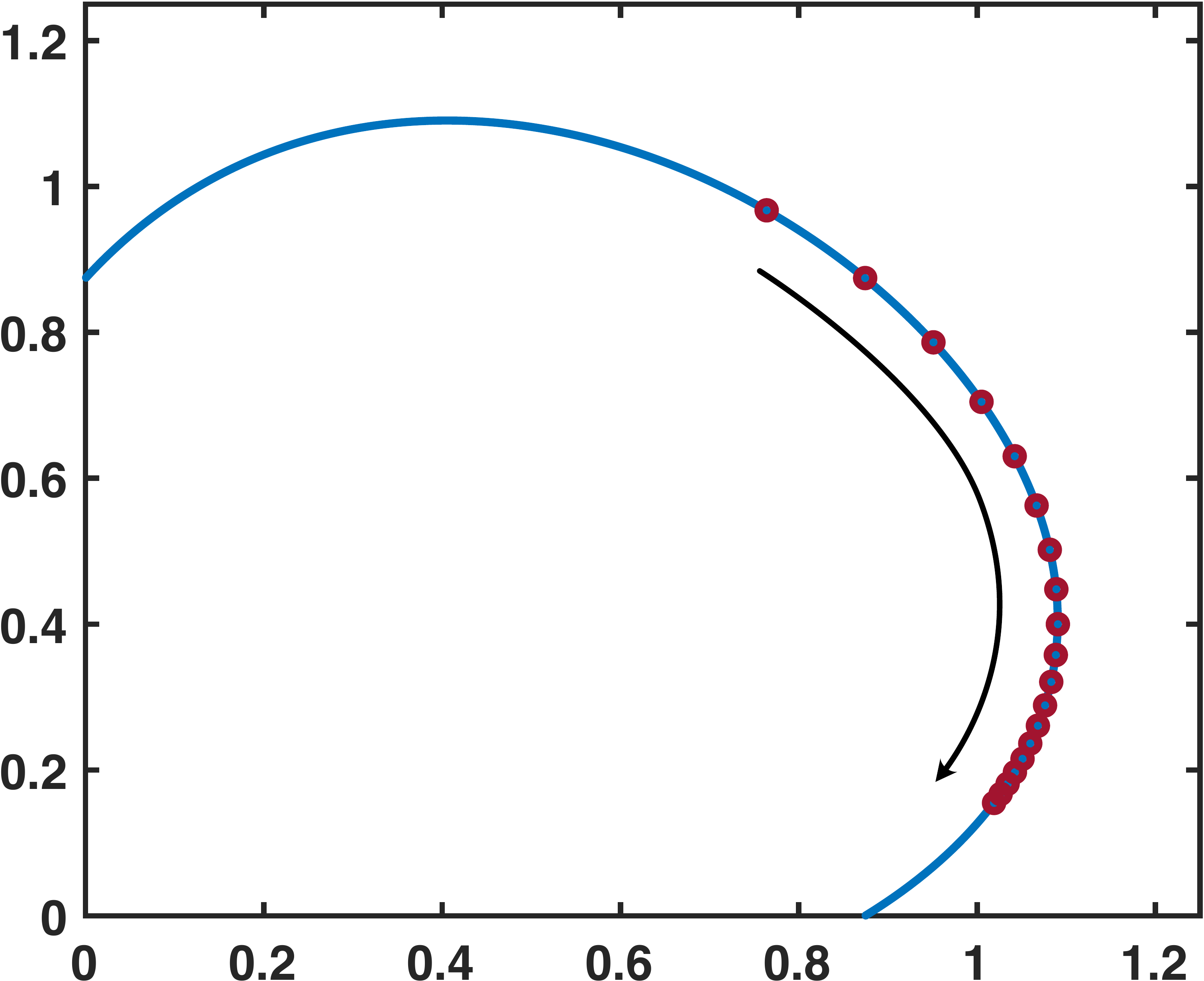}};
   \node[font=\fontsize{9pt}{0}\selectfont] at (2.9,-0.2) {$\mathcal{I}_1^\infty$};
   \node[font=\fontsize{9pt}{0}\selectfont] at (-0.2,2.5) {$\mathcal{I}_2^\infty$};
   \node[font=\fontsize{9pt}{0}\selectfont] at (4.4,2.3) {$\lambda_2$};
\end{tikzpicture}
\begin{tikzpicture}[scale=.8]
   \node[above right] at (0,0){\includegraphics[height=0.165\textheight]{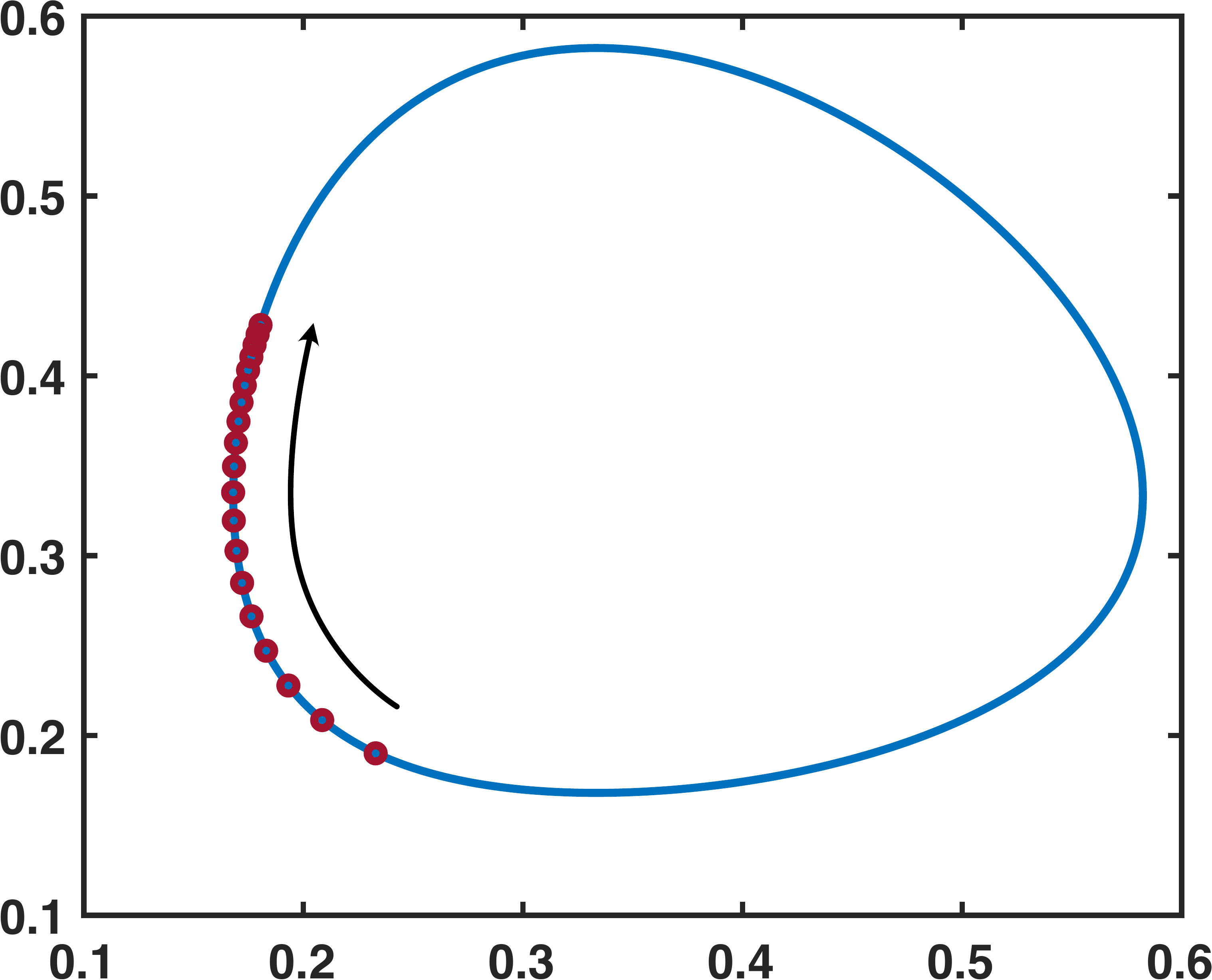}};
   \node[font=\fontsize{9pt}{0}\selectfont] at (2.9,-0.2) {$S_1^\infty$};
   \node[font=\fontsize{9pt}{0}\selectfont] at (-0.2,2.5) {$S_2^\infty$};
   \node[font=\fontsize{9pt}{0}\selectfont] at (2.1,2.3) {$\lambda_2$};
\end{tikzpicture} 
\begin{tikzpicture}[scale=.8]
   \node[above right] at (0,0){\includegraphics[height=0.165\textheight]{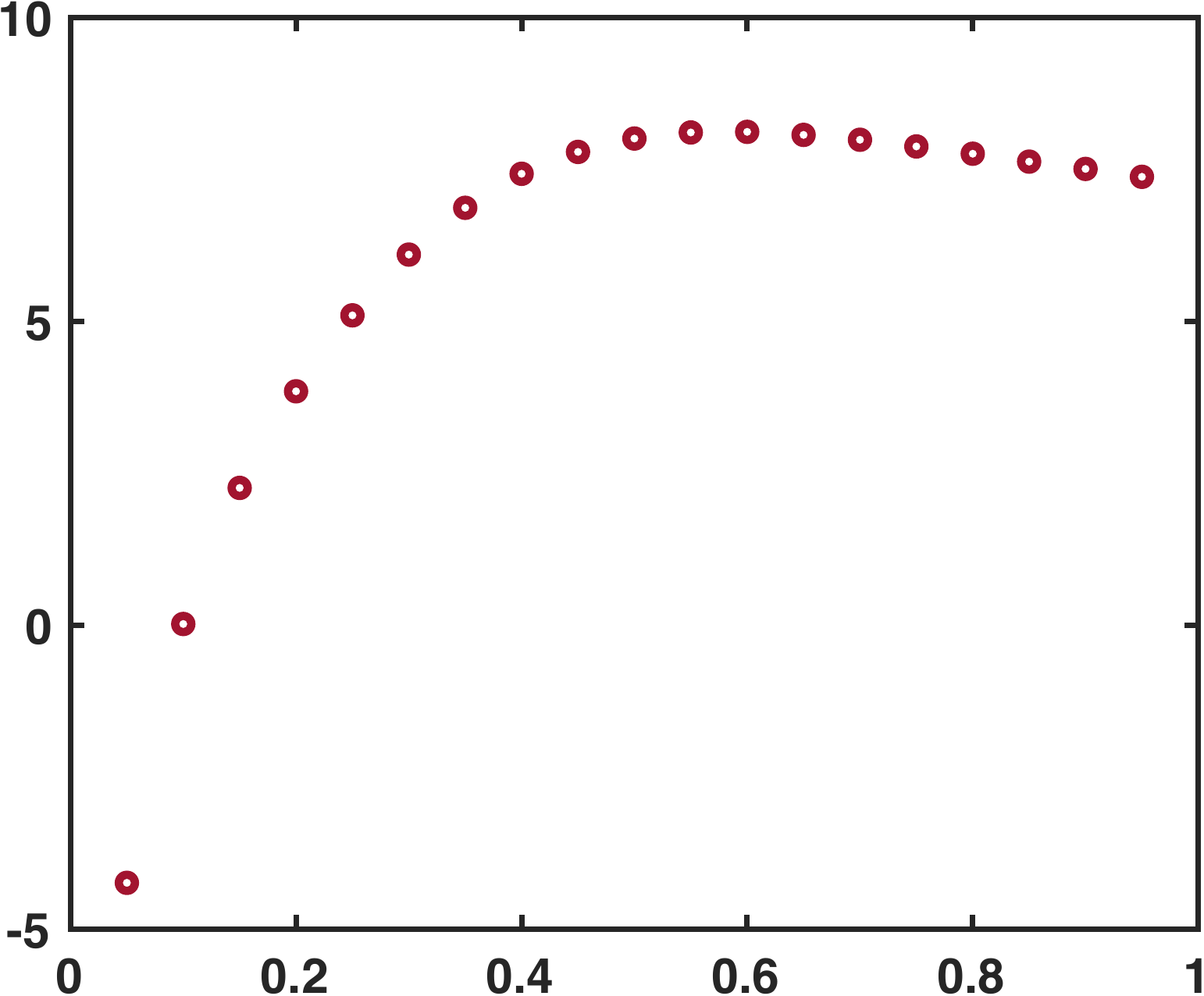}};
   \node[font=\fontsize{9pt}{0}\selectfont] at (2.9,-0.2) {$\lambda_2$};
   \node[font=\fontsize{9pt}{0}\selectfont] at (0,2.5) {$\Delta \mathcal{T}$};
\end{tikzpicture}\\[-2mm]
\begin{tikzpicture}[scale=.8]
   \node[above right] at (0,0){\includegraphics[height=0.165\textheight]{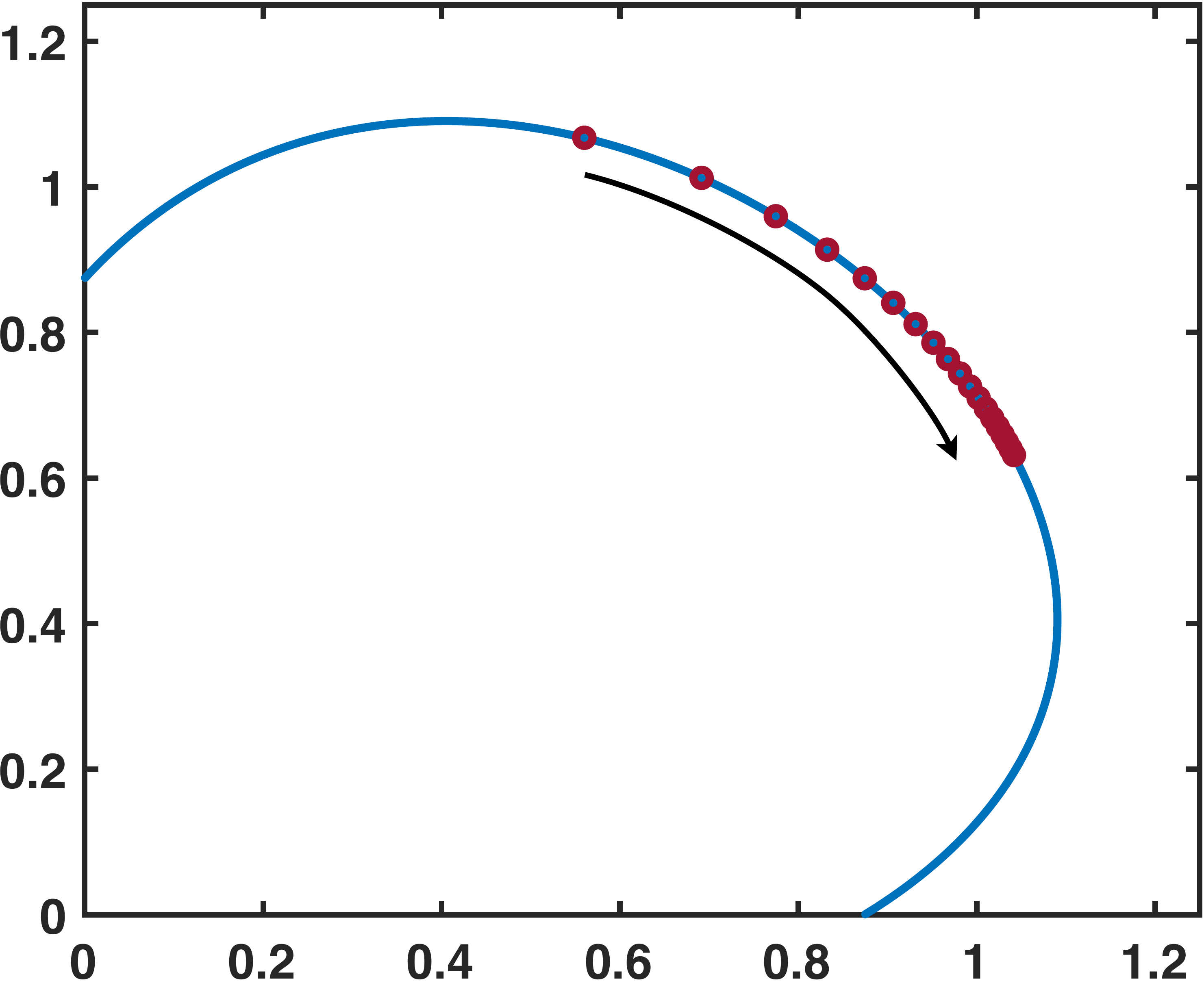}};
   \node[font=\fontsize{9pt}{0}\selectfont] at (2.9,-0.2) {$\mathcal{I}_1^\infty$};
   \node[font=\fontsize{9pt}{0}\selectfont] at (-0.2,2.5) {$\mathcal{I}_2^\infty$};
   \node[font=\fontsize{9pt}{0}\selectfont] at (3.5,3.2) {$\alpha_1$};
\end{tikzpicture}
\begin{tikzpicture}[scale=.8]
   \node[above right] at (0,0){\includegraphics[height=0.165\textheight]{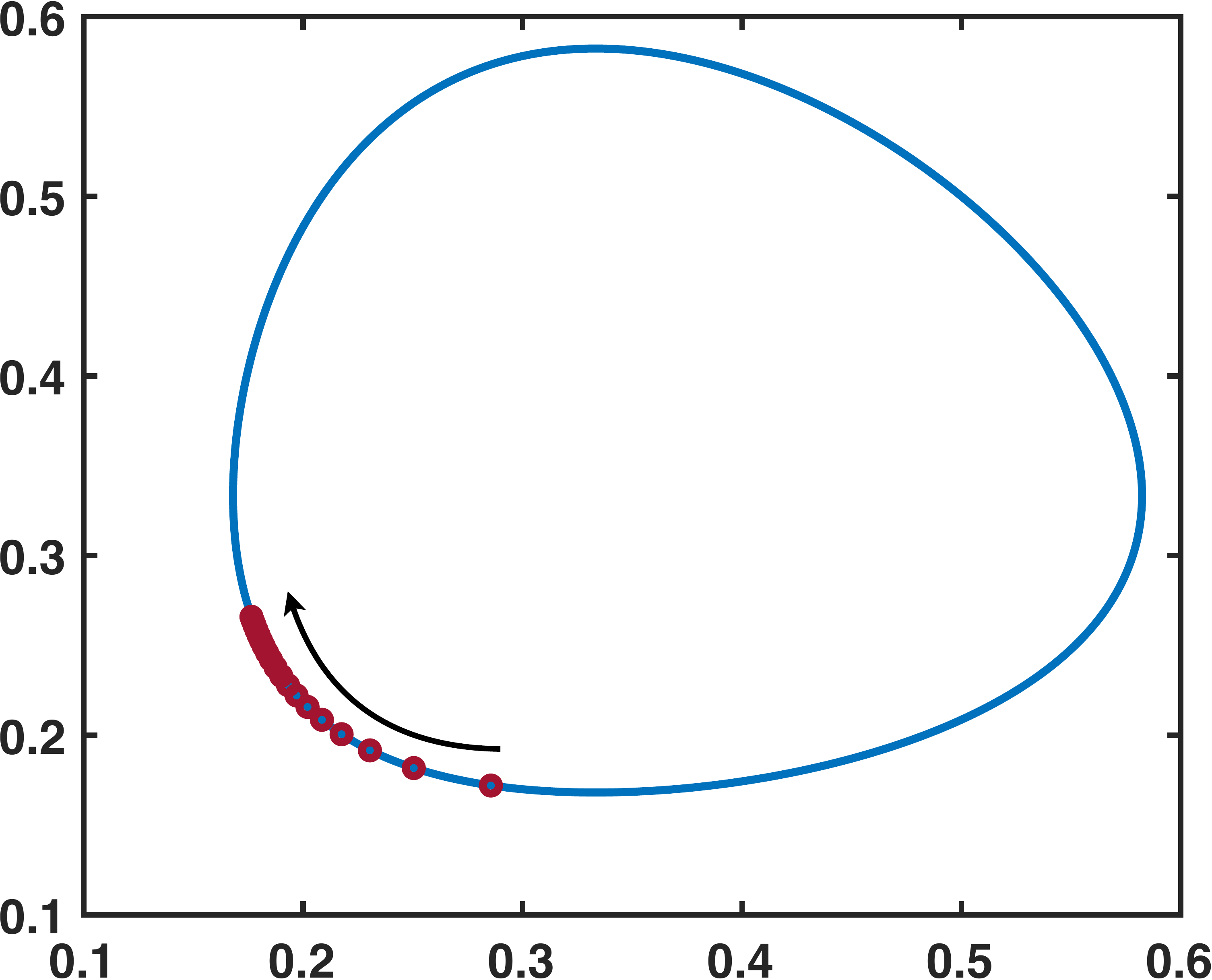}};
   \node[font=\fontsize{9pt}{0}\selectfont] at (2.9,-0.2) {$S_1^\infty$};
   \node[font=\fontsize{9pt}{0}\selectfont] at (-0.2,2.5) {$S_2^\infty$};
   \node[font=\fontsize{9pt}{0}\selectfont] at (2.1,2.0) {$\alpha_1$};
\end{tikzpicture}
\hspace*{-3mm}    
\begin{tikzpicture}[scale=.8]
   \node[above right] at (0,0){\includegraphics[height=0.165\textheight]{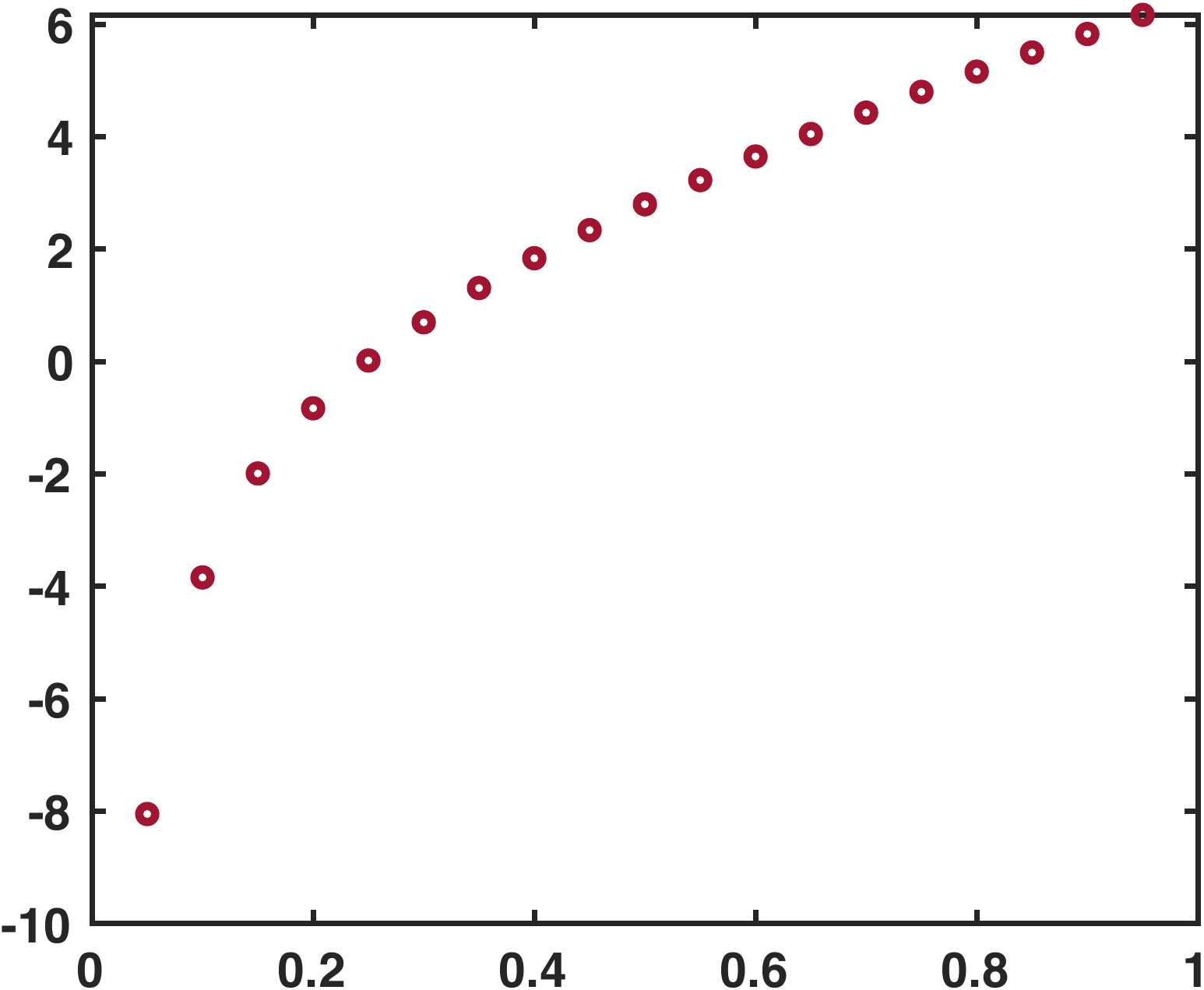}};
   \node[font=\fontsize{9pt}{0}\selectfont] at (2.9,-0.2) {$\alpha_1$};
   \node[font=\fontsize{9pt}{0}\selectfont] at (-0.2,2.5) {$\Delta \mathcal{T}$};
\end{tikzpicture}\\[-2mm]
\begin{tikzpicture}[scale=.8]
   \node[above right] at (0,0){\includegraphics[height=0.165\textheight]{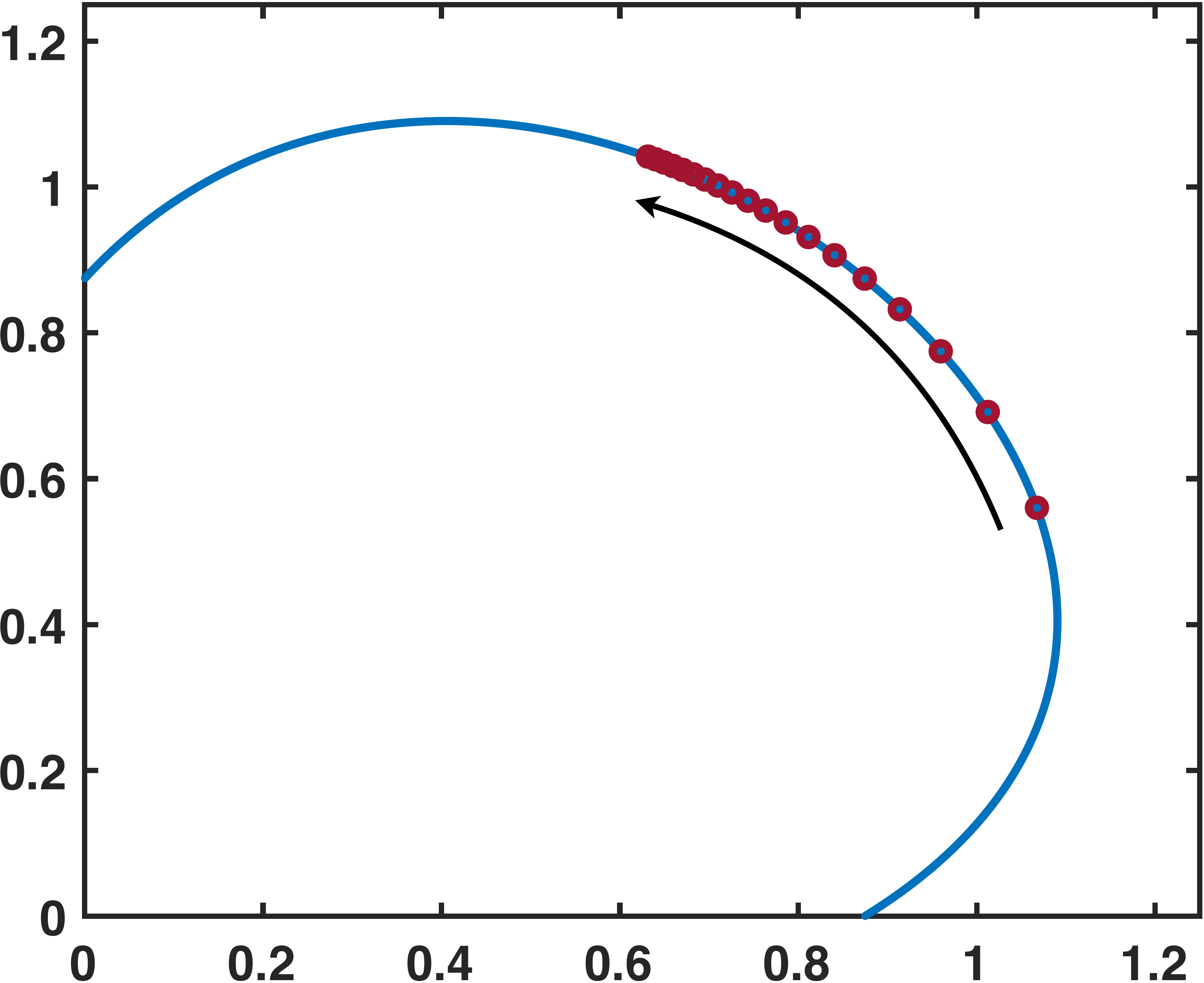}};
   \node[font=\fontsize{9pt}{0}\selectfont] at (2.9,-0.2) {$\mathcal{I}_1^\infty$};
   \node[font=\fontsize{9pt}{0}\selectfont] at (-0.2,2.5) {$\mathcal{I}_2^\infty$};
   \node[font=\fontsize{9pt}{0}\selectfont] at (3.9,2.9) {$\alpha_2$};
\end{tikzpicture}
\begin{tikzpicture}[scale=.8]
   \node[above right] at (0,0){\includegraphics[height=0.165\textheight]{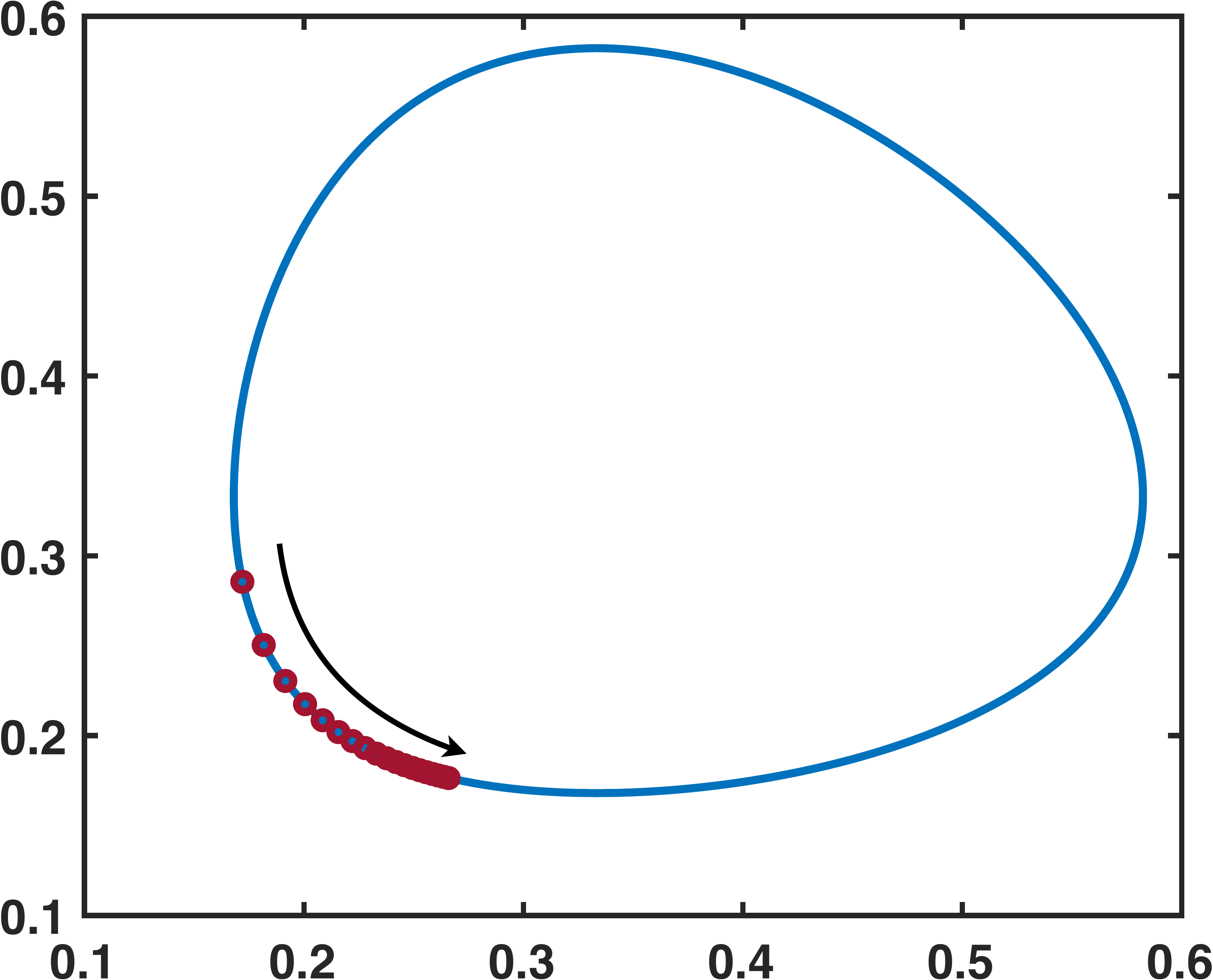}};
   \node[font=\fontsize{9pt}{0}\selectfont] at (2.9,-0.2) {$S_1^\infty$};
   \node[font=\fontsize{9pt}{0}\selectfont] at (-0.2,2.5) {$S_2^\infty$};
   \node[font=\fontsize{9pt}{0}\selectfont] at (2.1,1.9) {$\alpha_2$};
\end{tikzpicture} 
\begin{tikzpicture}[scale=.8]
   \node[above right] at (0,0){\includegraphics[height=0.165\textheight]{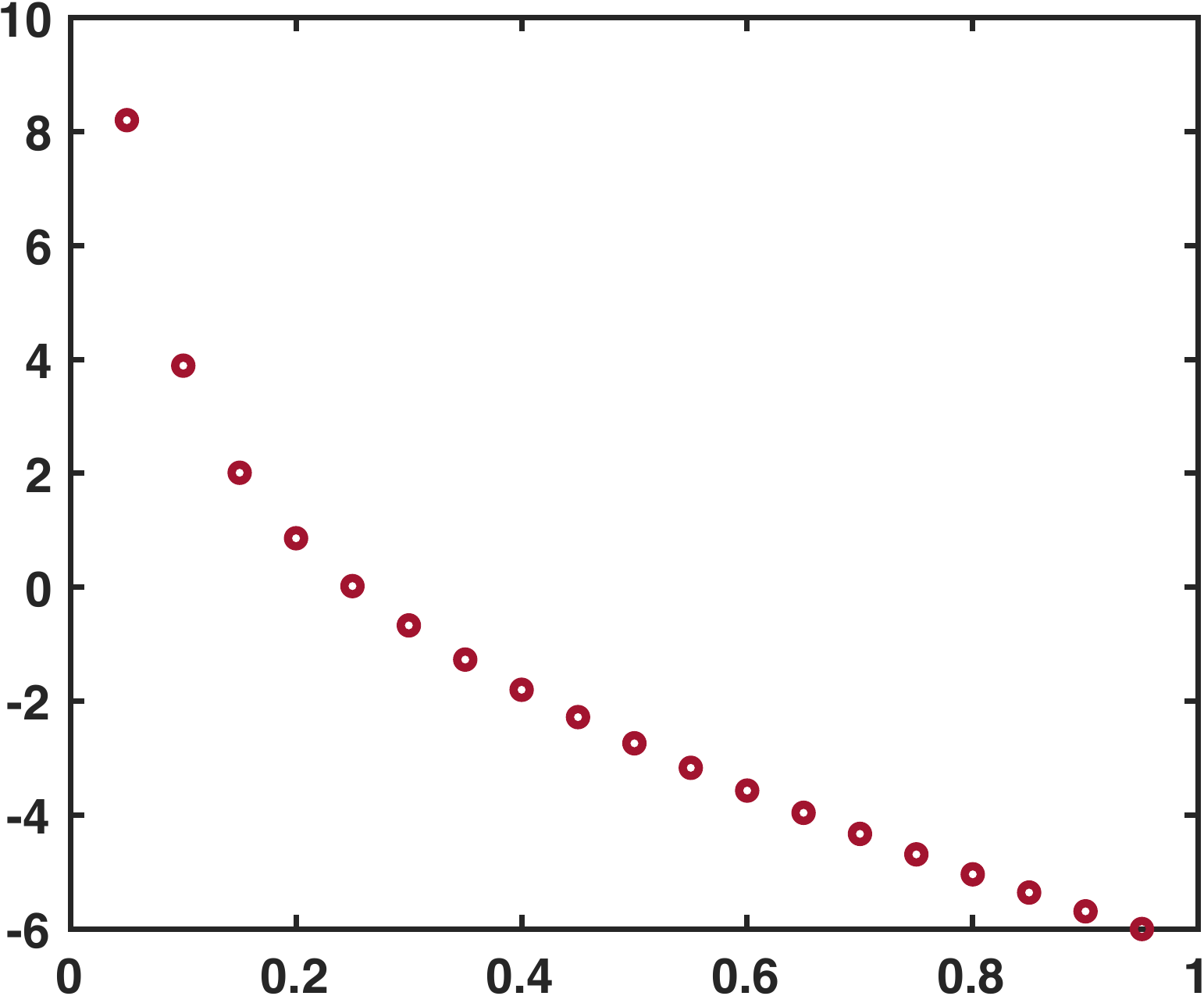}};
   \node[font=\fontsize{9pt}{0}\selectfont] at (2.9,-0.2) {$\alpha_2$};
   \node[font=\fontsize{9pt}{0}\selectfont] at (-0.2,2.5) {$\Delta \mathcal{T}$};
\end{tikzpicture}\\    
\end{tabular}

\caption{Final total populations of infected individuals (left) and corresponding final population of susceptible individuals (middle) as one parameter is varied from $0.05$ to $0.95$ while all parameters are fixed. The blue curve is the location of $(\mathcal{I}_1^\infty,\mathcal{I}_2^\infty)$ respectively $(S_1^\infty,S_2^\infty)$ while the dark red circles indicate the numerically computed values. Right: relative distance $\Delta \mathcal{T} := \mathcal{T}_2-\mathcal{T}_1$ between time of maximal infection $\mathcal{T}_j$ in each population, indicated by dark red circles, as the parameter is varied from $0.05$ to $0.95$. Varying parameters: $\lambda_1$ (top panel), $\lambda_2$ (second panel), $\alpha_1$ (third panel) and $\alpha_2$ (bottom panel). }
\label{fig:2c1rfinal}
\end{figure}
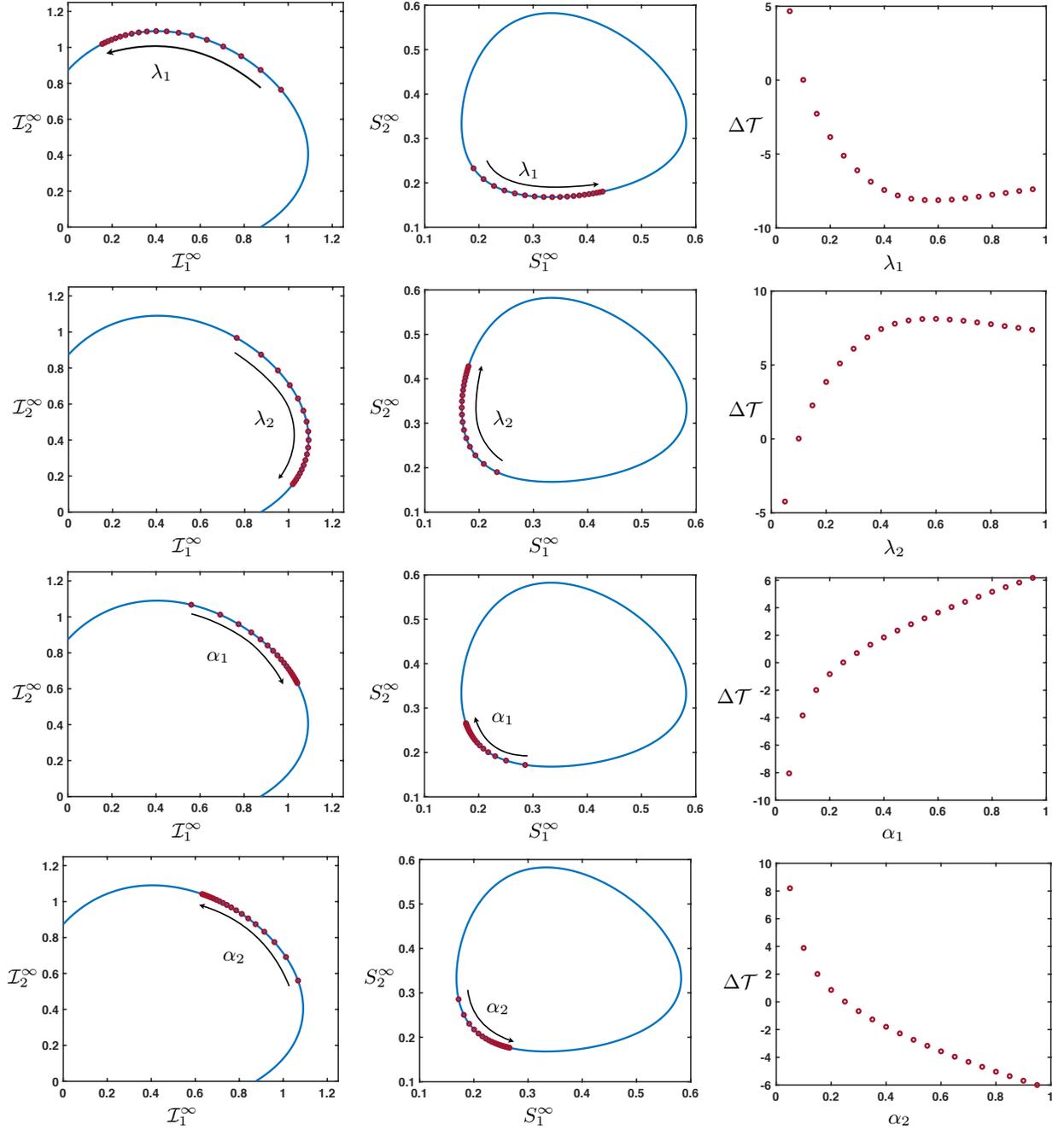

\begin{figure}[!t]
\centering
\includegraphics[width=0.5\textwidth]{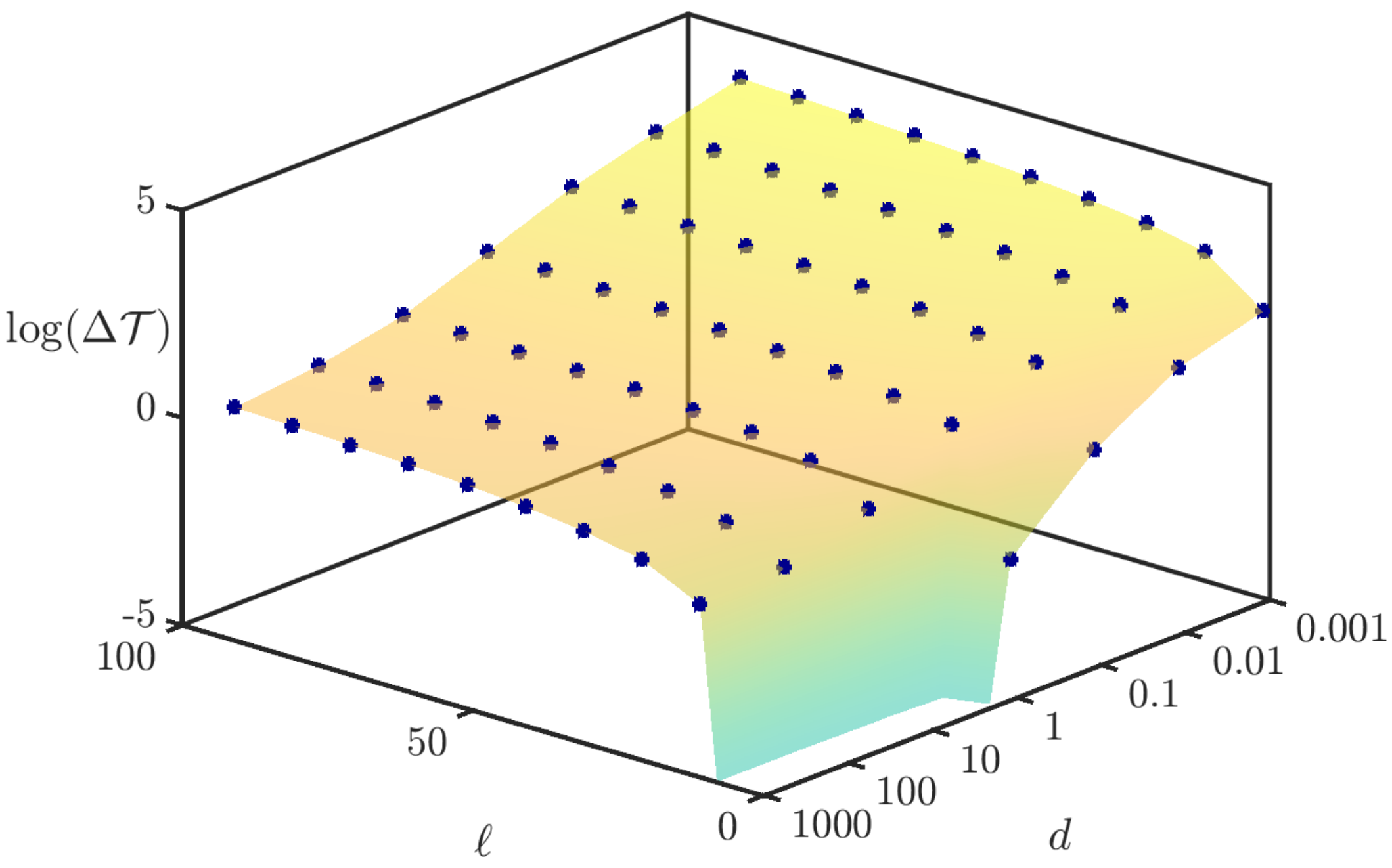}
\caption{Log-plot of the relative distance $\Delta \mathcal{T} = \mathcal{T}_2-\mathcal{T}_1$ between time of maximal infection $\mathcal{T}_j$ in each population $I_j(t)$ as the diffusion coefficient  $d$ and the length of the edge $\ell$ are varied while all other parameters are fixed to $\lambda_1=\lambda_2=1/10$, $\alpha_1=\alpha_2=1/4$, and $\tau_1=\tau_2=1$ with $\eta_1=\eta_2=1/3$. For the initial condition we have $(S_0,I_0)=(1/2,10^{-6})$. We note that as $d$ becomes smaller $\Delta \mathcal{T}$ rapidly increases as $\ell$ increases.}
\label{fig:2c1rdell}
\end{figure}

\begin{figure}[!t]
\centering
\includegraphics[width=0.45\textwidth]{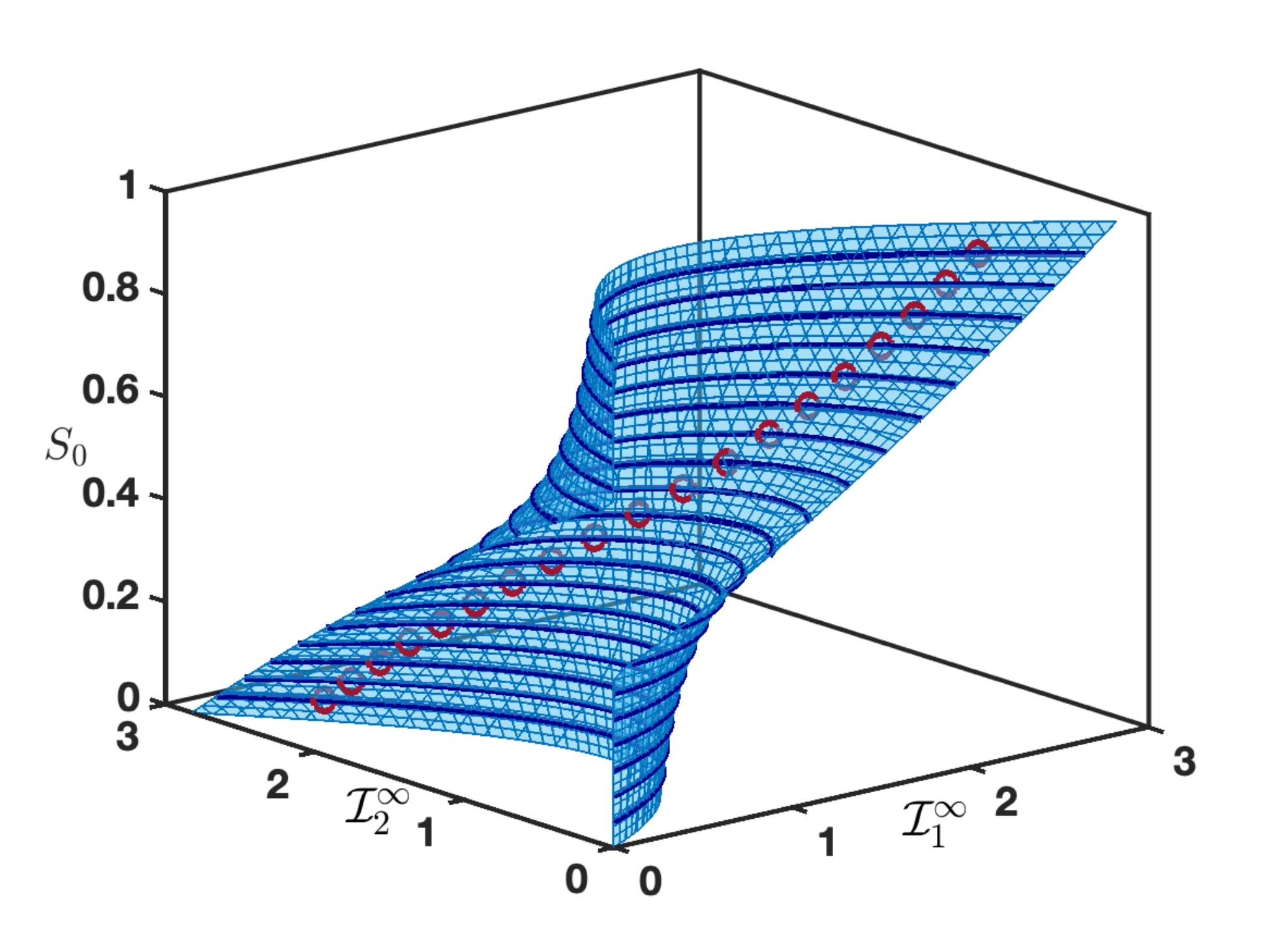}
\includegraphics[width=0.4\textwidth]{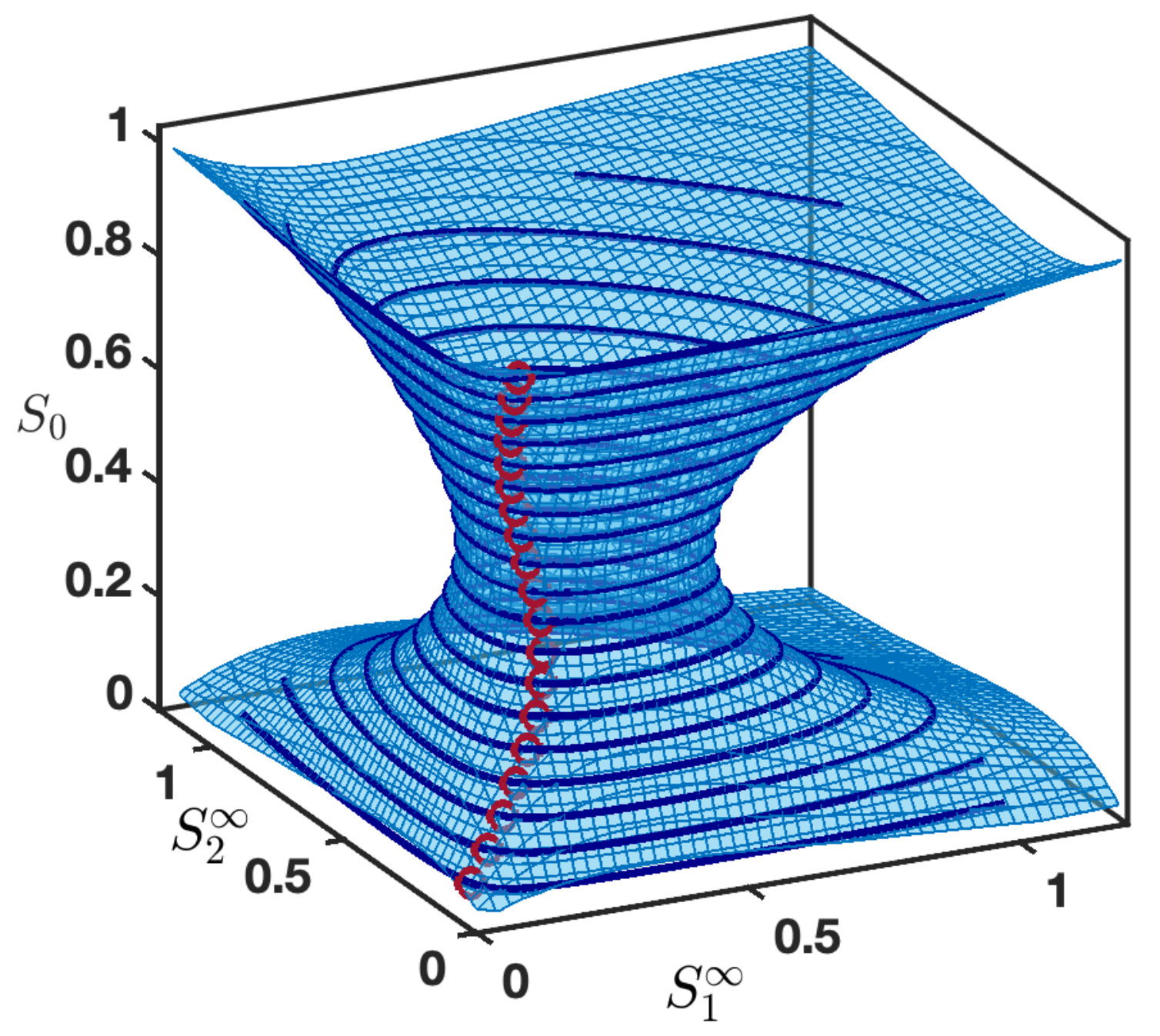}
\caption{Final total populations of infected individuals (left) and corresponding final population of susceptible individuals (right) as the initial population of susceptible individuals $S_0$ is varied from $0.05$ to $0.95$ while $I_0=10^{-6}$ is fixed. The dark blue curves are the location of $(\mathcal{I}_1^\infty,\mathcal{I}_2^\infty)$ respectively $(S_1^\infty,S_2^\infty)$ for each value of $S_0$, while the dark red circles indicate the numerically computed values. Each dark blue curve is a level set of the parameterized surface given by the conservation of total mass \eqref{manifoldIvinf}. All other parameters are fixed to $d=\ell=1$, $\lambda_1=\lambda_2=1/10$, $\alpha_1=\alpha_2=1/4$, and $\tau_1=\tau_2=1$ with $\eta_1=\eta_2=1/3$.  }
\label{fig:2c1rS0}
\end{figure}

For the numerical simulations, we have fixed initial conditions to be of the form
\bqs
u^0(x)=\frac{\lambda_1 I_0}{\alpha_1}\exp\left({-\dfrac{\alpha_2x^2}{2d\ell}}\right), x \in[0,\ell],
\eqs
with 
\bqs
(S_1^0,I_1^0,S_2^0,I_2^0)=\left(S_0-I_0-\int_0^\ell u^0(x)\md x,I_0,1-S_0,0\right), 
\eqs
where $S_0$ and $I_0$ may vary.
In Figures~\ref{fig:2c1r}-\ref{fig:2c1rfinal}-\ref{fig:2c1rdell}, $S_0$ and $I_0$ are fixed to $(S_0,I_0)=(1/2,10^{-6})$, while in Figure~\ref{fig:2c1rS0}, $S_0$ is allowed to vary and $I_0$ is fixed to $I_0=10^{-6}$.

In Figure~\ref{fig:2c1r}, we report the profiles of the solutions $(S_i(t),I_i(t))$ together with the total population on the edge $\int_0^\ell u(t,x)\md x$ and the total mass of the system $M(t)$ as the parameter $\lambda_1$ is varied from $0.05$ to $0.95$, while all other parameters are being kept fixed. We observe that the dynamics of the epidemic at the second vertex is almost independent of the parameter $\lambda_1$ while it has a significant impact on the dynamics at the first vertex. Indeed, as $\lambda_1$ is increased, the maximum of infected individuals $\max_{t\geq 0} I_1(t)$ is decreased. In the last panel of the figure, we also illustrate the conservation of total population where the fluctuations around $M^0=1$ is of order $10^{-12}$. In the top panel of Figure~\ref{fig:2c1rfinal}, we present the final total populations of infected individuals  and corresponding final population of susceptible individuals as $\lambda_1$ is varied. The blue curve is the location of $(\mathcal{I}_1^\infty,\mathcal{I}_2^\infty)$ respectively $(S_1^\infty,S_2^\infty)$ while the dark red circles indicate the numerically computed values. We recover the fact that $\lambda_1$ has a more significant impact on the final total populations at the first vertex than it has at the second vertex. The get a better understanding of the intricate dynamics between the epidemic at the two vertices, we also present the relative distance $\Delta \mathcal{T} := \mathcal{T}_2-\mathcal{T}_1$ between time of maximal infection $\mathcal{T}_j$ in each population as $\lambda_1$ is varied. We observe that $\Delta \mathcal{T}$ is not monotone in $\lambda_1$, as it  first decreases and then increases. But we also note that $\Delta \mathcal{T}<0$ for $\lambda_1\geq 0.1$ traducing the fact that the pick of the epidemic occurs at the second vertex before it does at the first vertex, although initially $I_2^0=0$. This illustrates the effect of the diffusion of infected individuals along the edge.

Similarly, in Figure~\ref{fig:2c1rfinal}, we report  the final total populations of infected individuals and corresponding final population of susceptible individuals as $\lambda_2$ (second panel), $\alpha_1$ (third panel) and $\alpha_2$ (bottom panel) are varied from $0.05$ to $0.95$. As expected, the final total population of infected individuals at the second vertex decreases as $\lambda_2$ increases while at the first vertex it varies less significantly.  As $\alpha_1$ increases, the final total population of infected individuals at the first vertex increases while it decreases at the second vertex. This time the relative distance $\Delta \mathcal{T} := \mathcal{T}_2-\mathcal{T}_1$ between time of maximal infection is monotonically increasing with $\alpha_1$. We get the opposite monotonicity properties as $\alpha_2$ is varied.

In Figure~\ref{fig:2c1rdell}, we investigate the joint effect of the diffusion coefficient $d$ and the length of the edge $\ell$ on the dynamics of the epidemic at the vertices. Here, we focus on the delay between time of maximal infection $\mathcal{T}_j$ in each infected population $I_j(t)$. As expected, when the diffusion coefficient is really small while the length is being kept at order one, $\Delta \mathcal{T}$ takes large value: $\Delta \mathcal{T} \sim 10^4$ when $d=10^{-3}$ and $\ell=1$. Biologically, this means that when the diffusion coefficient is really small  it takes more time for infected individuals from vertex one to reach the second vertex and start an epidemic. We also note that at fixed $\ell$, $\Delta \mathcal{T}$ monotonically decreases as $d$ increases, while at fixed $d$, $\Delta \mathcal{T}$ monotonically increases as $\ell$ increases.

\begin{figure}[!t]
\centering
\includegraphics[width=0.45\textwidth]{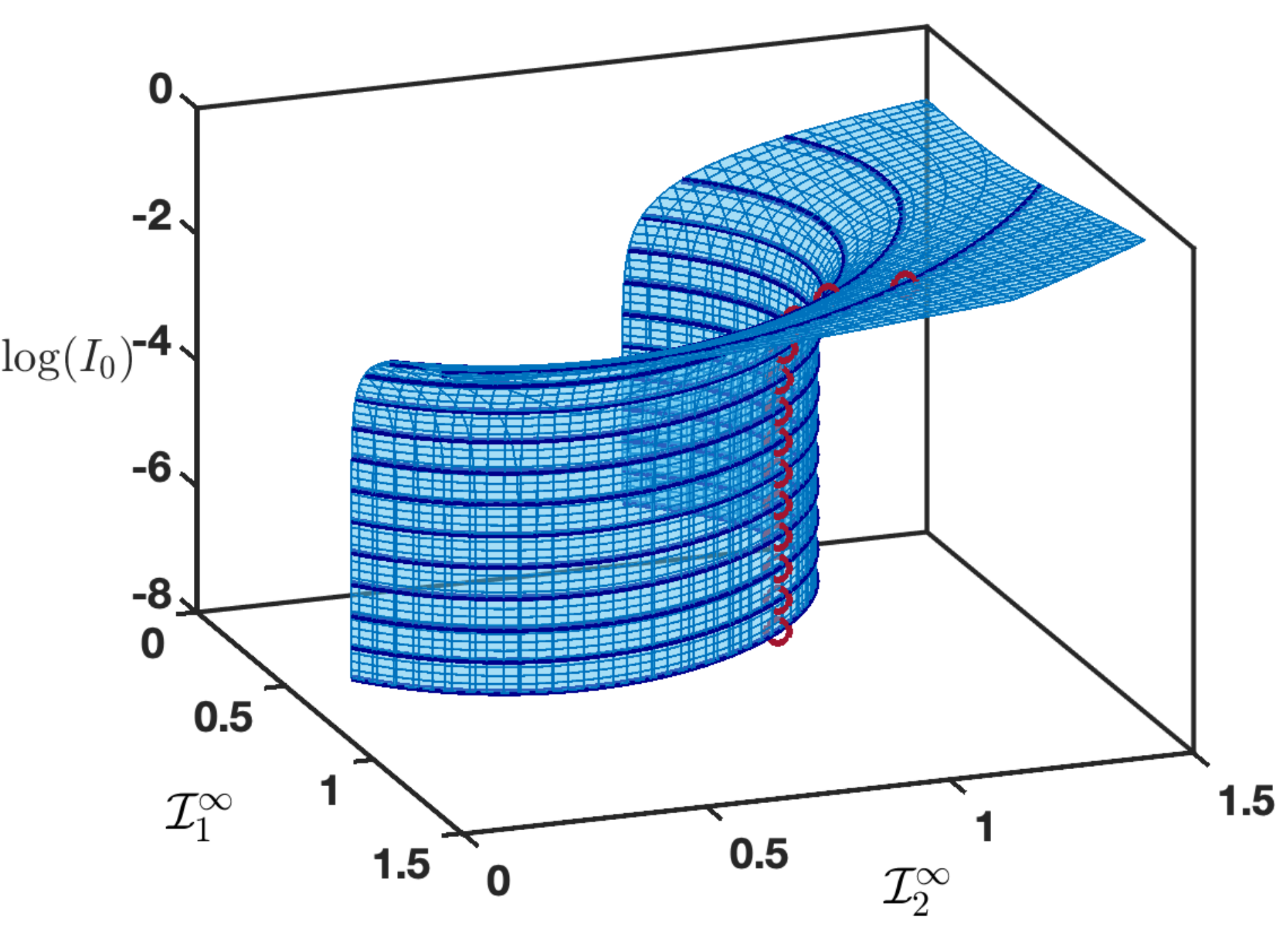}
\includegraphics[width=0.45\textwidth]{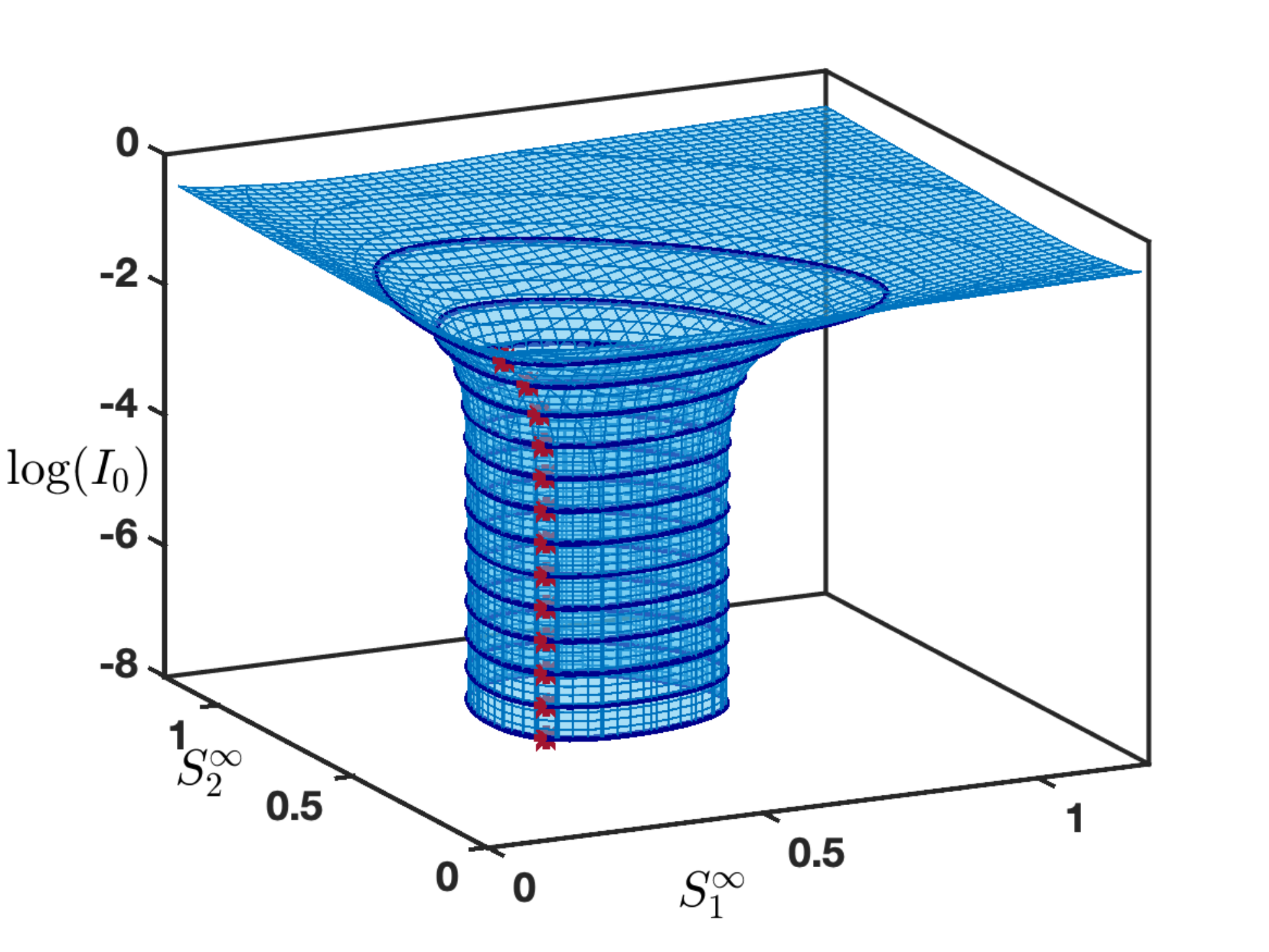}
\caption{Final total populations of infected individuals (left) and corresponding final population of susceptible individuals (right) as the initial population of susceptible individuals $I_0$ is varied from $10^{-7}$ to $10^{-1}$ in log-scale while $S_0=1/2$ is fixed. The dark blue curves are the location of $(\mathcal{I}_1^\infty,\mathcal{I}_2^\infty)$ respectively $(S_1^\infty,S_2^\infty)$ for each value of $S_0$, while the dark red circles indicate the numerically computed values. Each dark blue curve is a level set of the parameterized surface given by the conservation of total mass \eqref{manifoldIvinf}. All other parameters are fixed to $d=\ell=1$, $\lambda_1=\lambda_2=1/10$, $\alpha_1=\alpha_2=1/4$, and $\tau_1=\tau_2=1$ with $\eta_1=\eta_2=1/3$.  }
\label{fig:2c1rI0}
\end{figure}

\begin{figure}[!htbp]
  \centering
  \begin{tabular}{cc}
    \includegraphics[height=0.27\textheight]{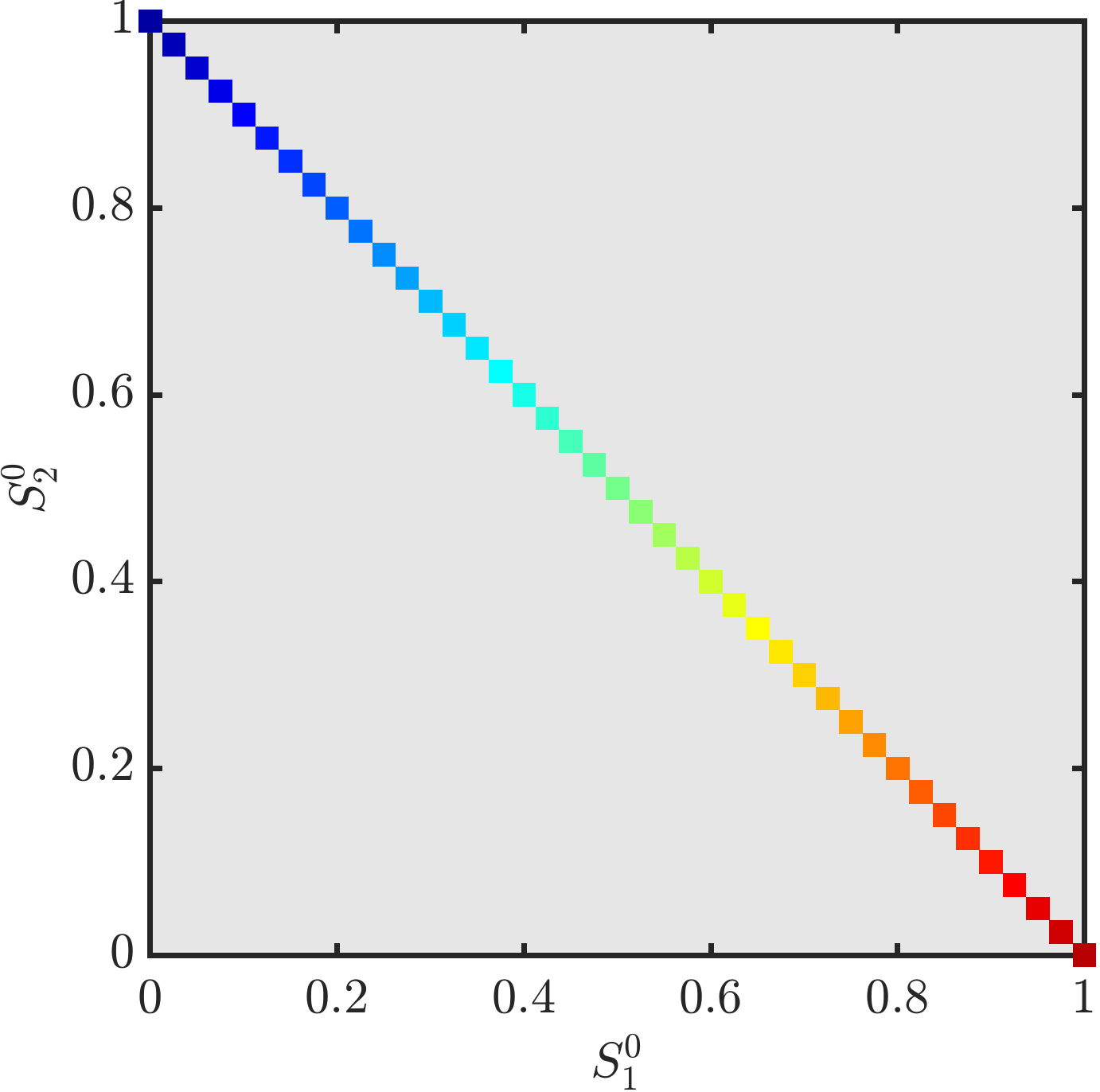}&\includegraphics[height=0.27\textheight]{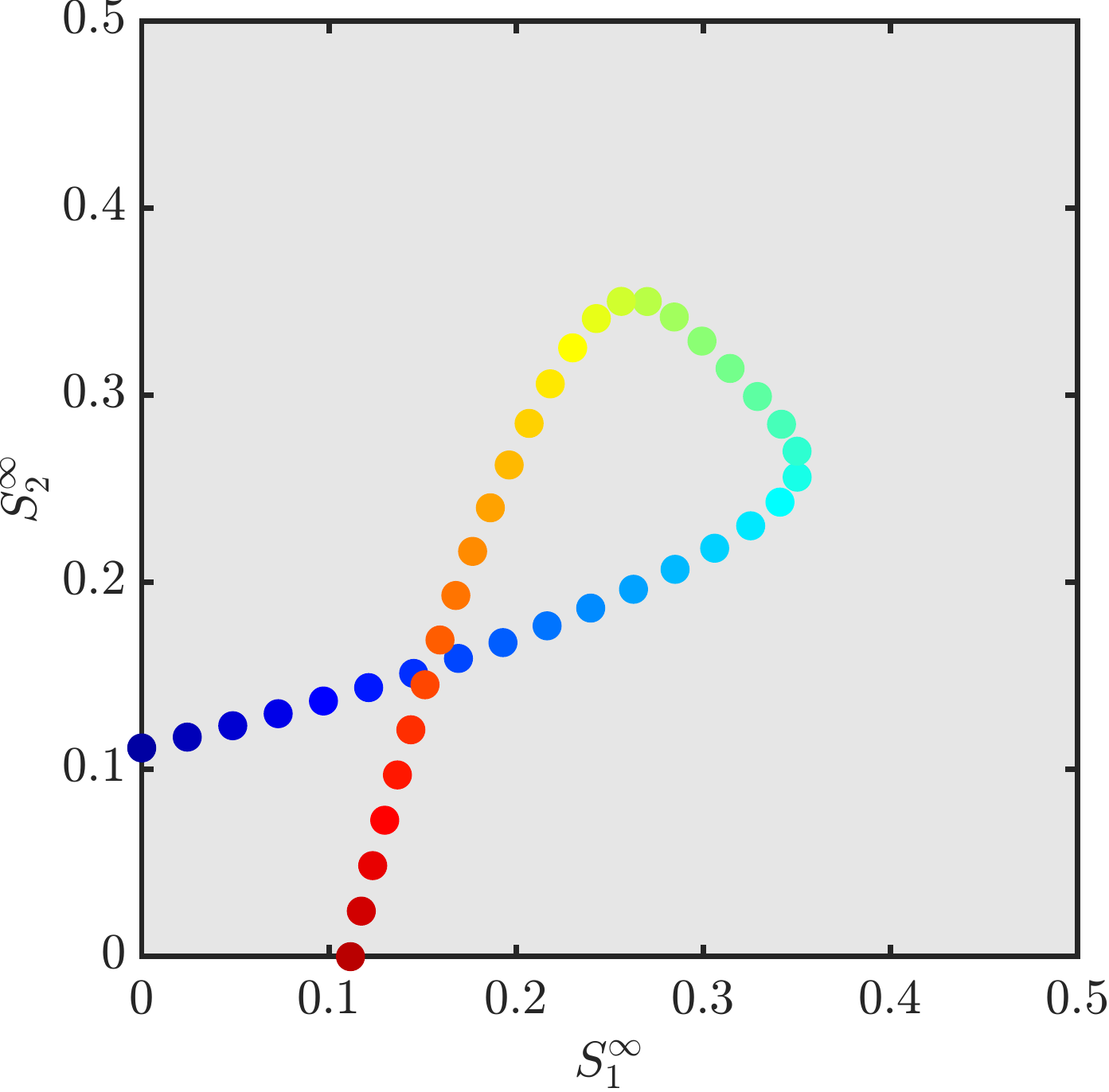}\\%&\includegraphics[height=0.22\textheight]{eval_Rinfini.pdf}\\
  \end{tabular}
  \begin{tabular}{cc}
    \includegraphics[height=0.25\textheight]{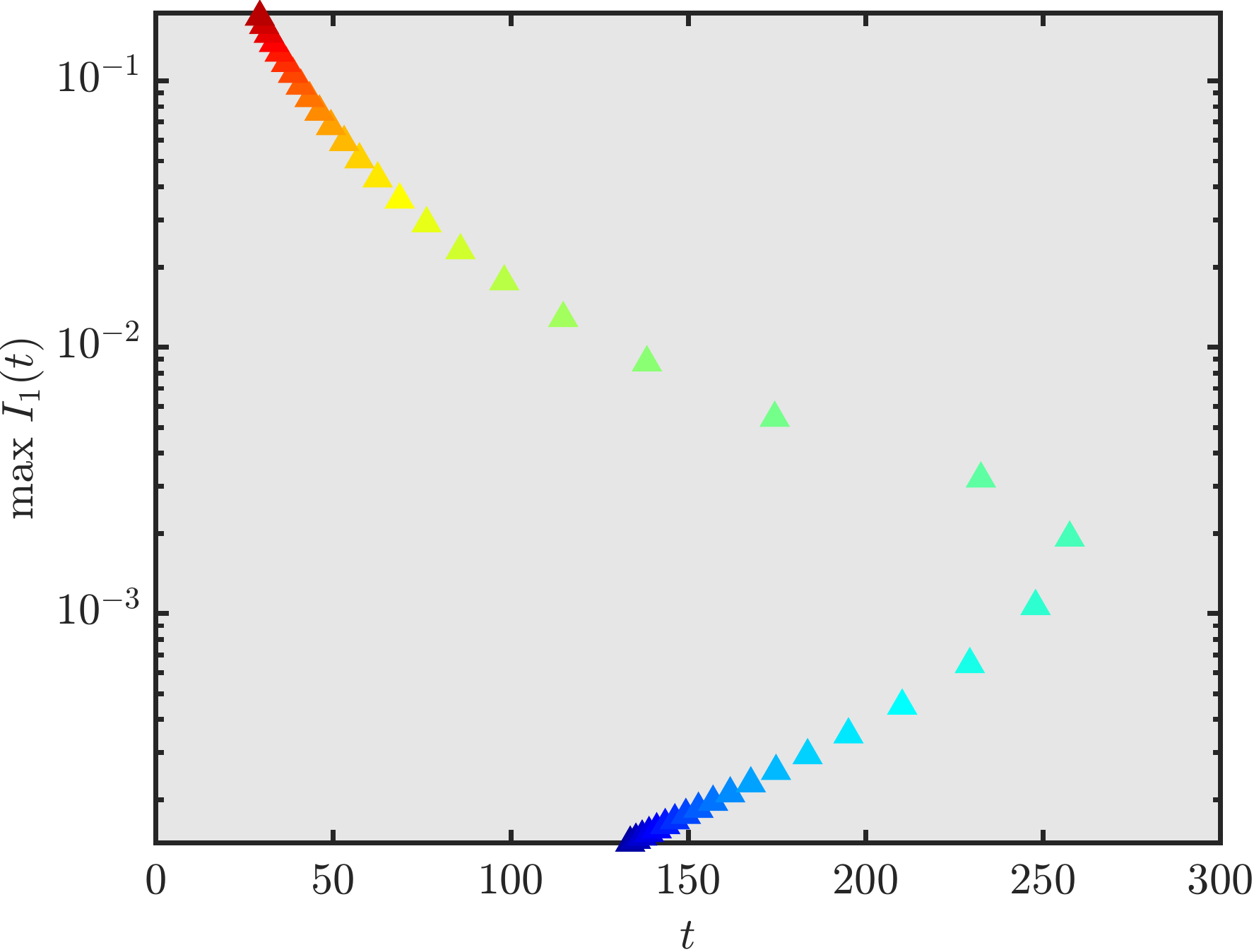}&\includegraphics[height=0.25\textheight]{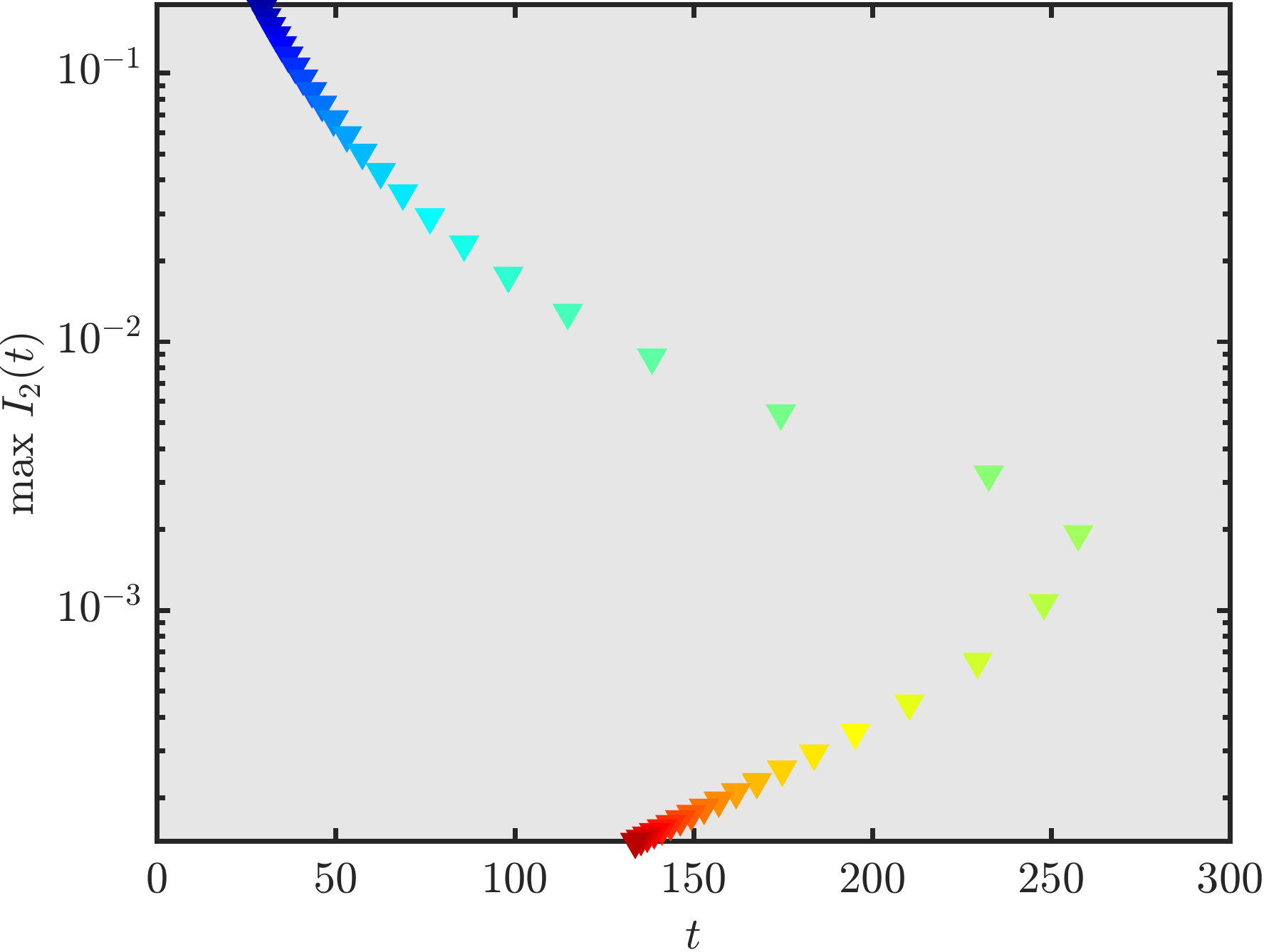}\\
  \end{tabular}
  \caption{Locations of $S_{1,2}^\infty$ (top) and $\max I_{1,2}(t)$ (bottom left and right) as functions of $S_{1,2}^0$. Values of all other parameters are similar to Figure~\ref{fig:2c1rS0}. The initial condition is of the form $(S_0,1-S_0)$ with $S_0\in[0,1]$ and to each initial configuration is associated a color code from blue to red. The curve in the top right panel is a projection on the $(S_{1}^\infty,S_{2}^\infty)$-plane of the parametrized a curve from Figure~\ref{fig:2c1rS0}, right panel.} 
  \label{fig:toto}
\end{figure}

In Figures~\ref{fig:2c1rI0}-\ref{fig:2c1rS0}, we vary respectively the initial population of susceptible individuals $S_0$ and infected individuals $I_0$. We visualize the final total populations of infected individuals  and corresponding final population of susceptible individuals on the parameterized surfaces $(\mathcal{I}_1^\infty,\mathcal{I}_2^\infty,S_0)$ and $(S_1^\infty,S_2^\infty,S_0)$, respectively $(\mathcal{I}_1^\infty,\mathcal{I}_2^\infty,I_0)$ and $(S_1^\infty,S_2^\infty,I_0)$, where the level sets of the parameterized surface are given by the conservation of total mass \eqref{manifoldIvinf}. We note that $(\mathcal{I}_1^\infty,\mathcal{I}_2^\infty)$ and $(S_1^\infty,S_2^\infty)$ are almost independent of $I_0$ when $I_0\leq 10^{-3}$ with sensible variations only occurring for larger values of $I_0$. On the other hand, we observe that as $S_0$ is increased the final total population of infected individuals increases at the first vertex while it decreases at the second one. The dependence of $(S_1^\infty,S_2^\infty)$ as a function of $S_0$ is more subtile and is presented in Figure~\ref{fig:toto}. In the same figure, we also show the location of $\max I_j(t)$ and its amplitude. We observe a strong nonlinear dependence with respect to $S_0$. As $S_0$ increases, we first see that the time at which $I_1(t)$ is maximal increases and then decreases, while $\max I_1(t)$ is monotonically increasing. The converse is observed at the second vertex.

\subsection{Case of 3 vertices and 3 edges}
\begin{center}
\begin{tikzpicture}
    \coordinate (A) at (-2,0);
    \coordinate (B) at (2,-2);
    \coordinate (C) at (1.5,1.5);
    
    \node[below] at (A) {$v_1$};
    \node[below] at (B) {$v_2$};
    \node[above] at (C) {$v_3$};

    % \node[label=left:$A$]  at (A) {$\bullet$};
    % \node[label=right:$B$] at (B)  {$\bullet$};
    % \node[label=right:$C$] at (C) {$\bullet$};

    \node at (A) {$\bullet$};
    \node at (B) {$\bullet$};
    \node at (C) {$\bullet$};

    % \draw (A) -- (B) -- (C) -- (A);
    \draw (A) -- (B) node[midway,below] {$A$};
    \draw (B) -- (C) node[midway,right] {$B$};
    \draw (C) -- (A) node[midway,above] {$C$};
\end{tikzpicture}
\end{center}

\begin{figure}[!t]
\centering
\includegraphics[width=0.45\textwidth]{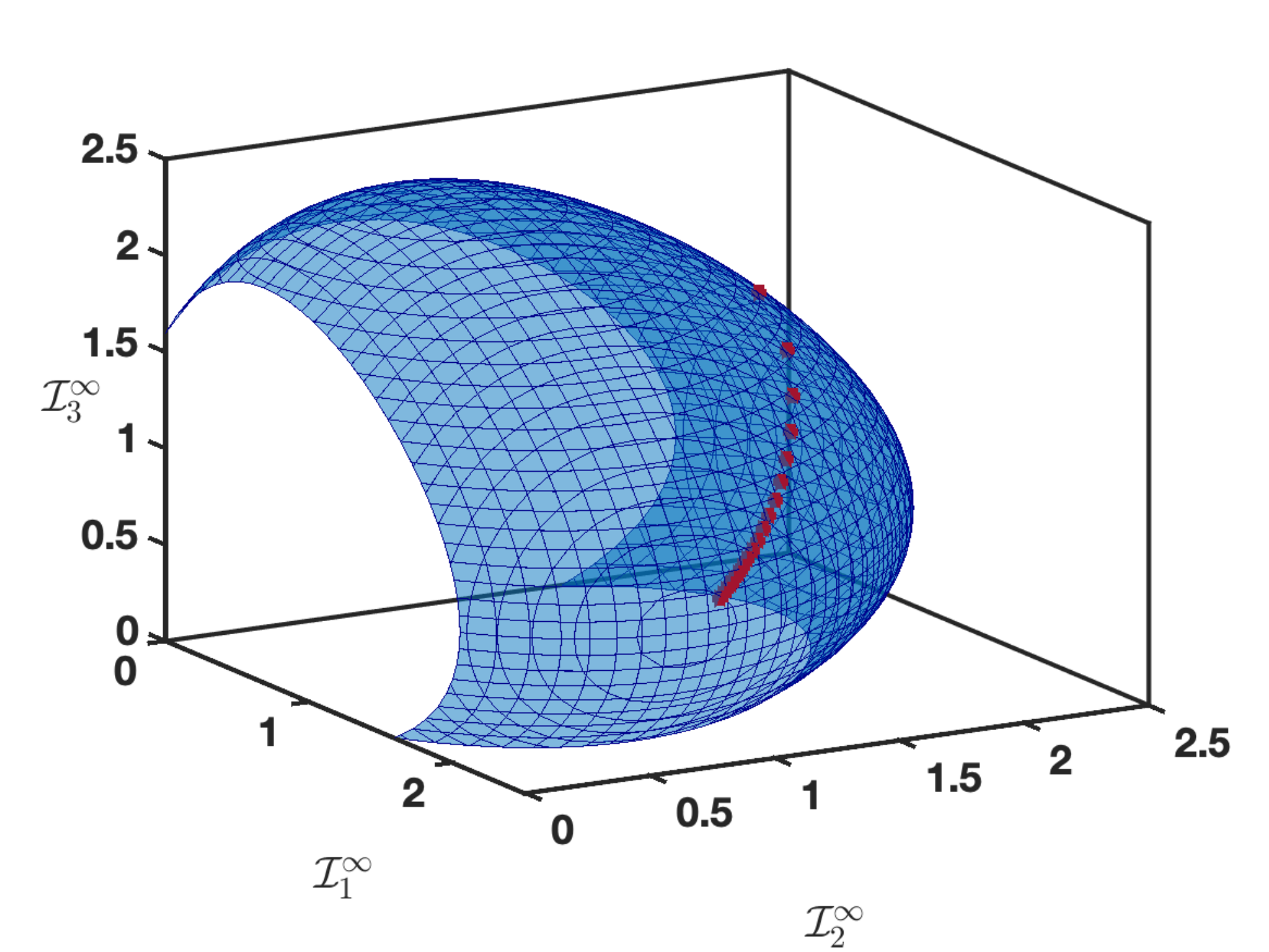}
\includegraphics[width=0.45\textwidth]{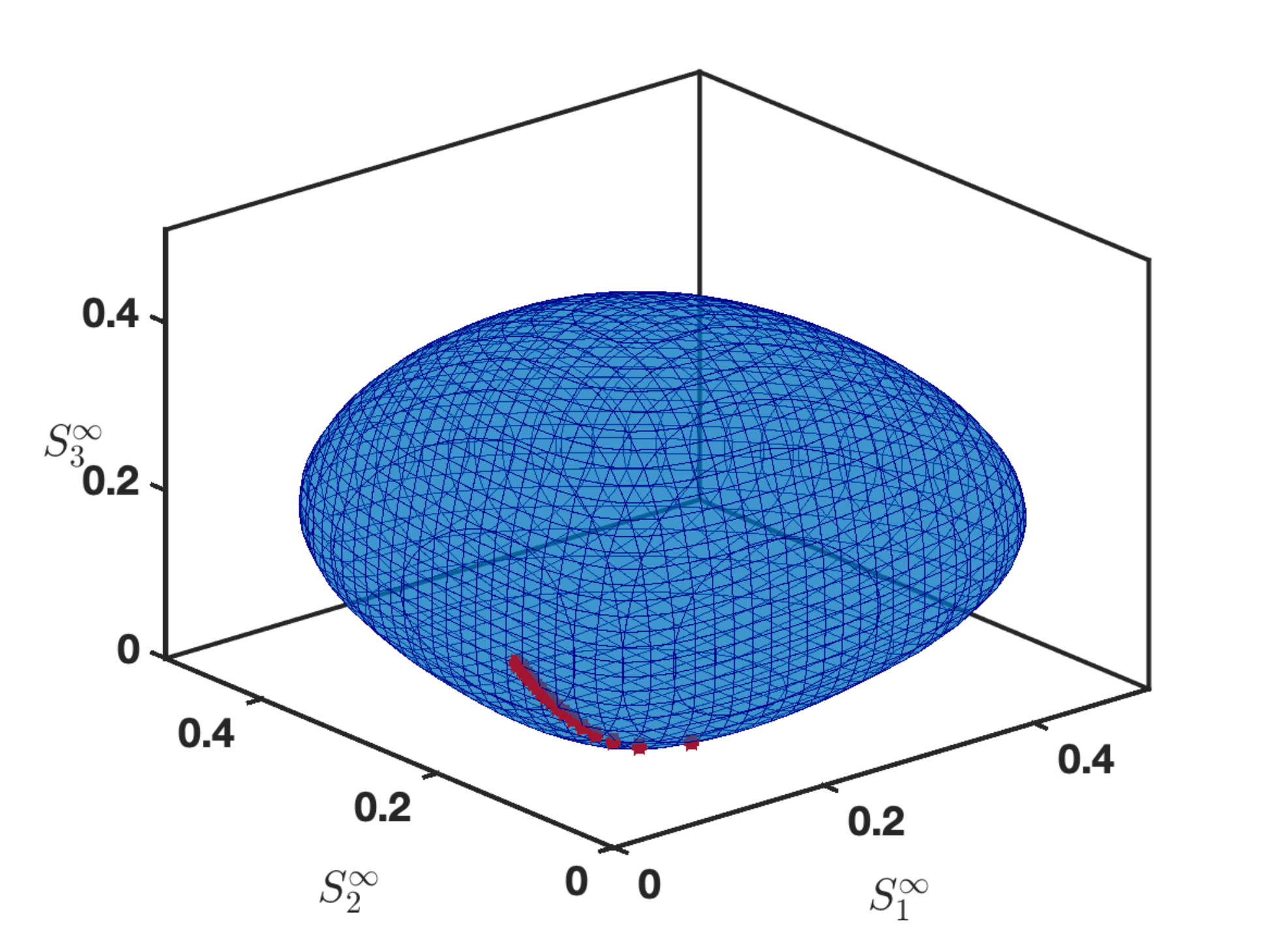}\hspace{0.2cm}
\caption{Final total populations of infected individuals (left) and corresponding final population of susceptible individuals (right) as $\nu$ is varied from $0.05$ to $0.95$. The dark blue surfaces are the location of $(\mathcal{I}_1^\infty,\mathcal{I}_2^\infty,\mathcal{I}_3^\infty)$ respectively $(S_1^\infty,S_2^\infty,S_3^\infty)$ for each value of $S_0$, while the dark red circles indicate the numerically computed values. Parameters were set to $\ell=d=1$, $(\tau,\eta)=(1,1/6)$, and $(\alpha,\lambda)=(1/8,1/10)$, while the initial condition is $(S_0,I_0)=(1,10^{-6})$.}
\label{fig:3c3r}
\end{figure}

Next, we consider the case of $3$ vertices and $3$ edges arranged in a triangular configuration. For the numerical simulations presented in Figure~\ref{fig:3c3r}, we have assumed full symmetry in the parameters that is
\begin{align*}
& (\ell_e,d_e)=(\ell,d), \quad e\in\E,\quad (\tau_v,\eta_v)=(\tau,\eta), \quad v\in\V,\\
& (\alpha_e^v,\lambda_e^v)=(\alpha,\lambda), \quad (e,v)\in\E\times\V, \quad \nu_{e,e'}^v=\nu, \quad (e,e',v)\in\E\times\E\times\V.
\end{align*}
 Regarding the initial condition, we have chosen
 \bqs
 (S^0_{1},I^0_1,S^0_2,I^0_2,S^0_3,I^0_3)=\left(S^0-I_0,I_0,S^0,0,S^0,0\right) \quad v\in\V,
 \eqs
 for a given $(S^0,I_0)$, while for each $e\in\E$ we have set $u_e^0(x)=0$ on $\Omega_e$. Note that, we have initially a boundary layer as our initial condition does not satisfy \eqref{BCgraph} for small times. We remark that the final total populations of infected individuals  and corresponding final population of susceptible individuals belong to a surface as provided by \eqref{manifoldIvinf}-\eqref{manifoldSvinf} from Theorem~\ref{thmFTI}.

\begin{figure}[!t]
\centering
\includegraphics[width=0.45\textwidth]{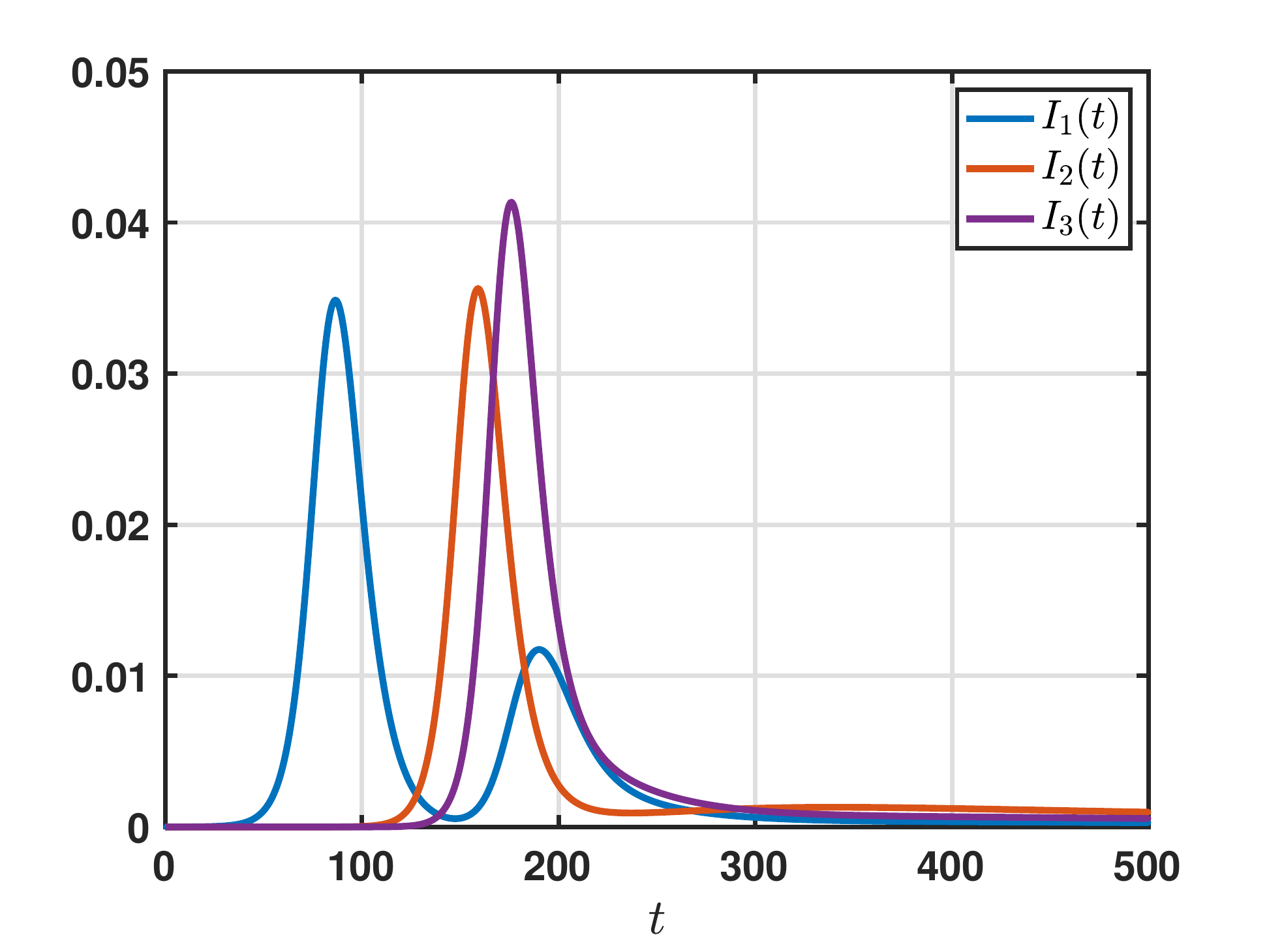}
\includegraphics[width=0.45\textwidth]{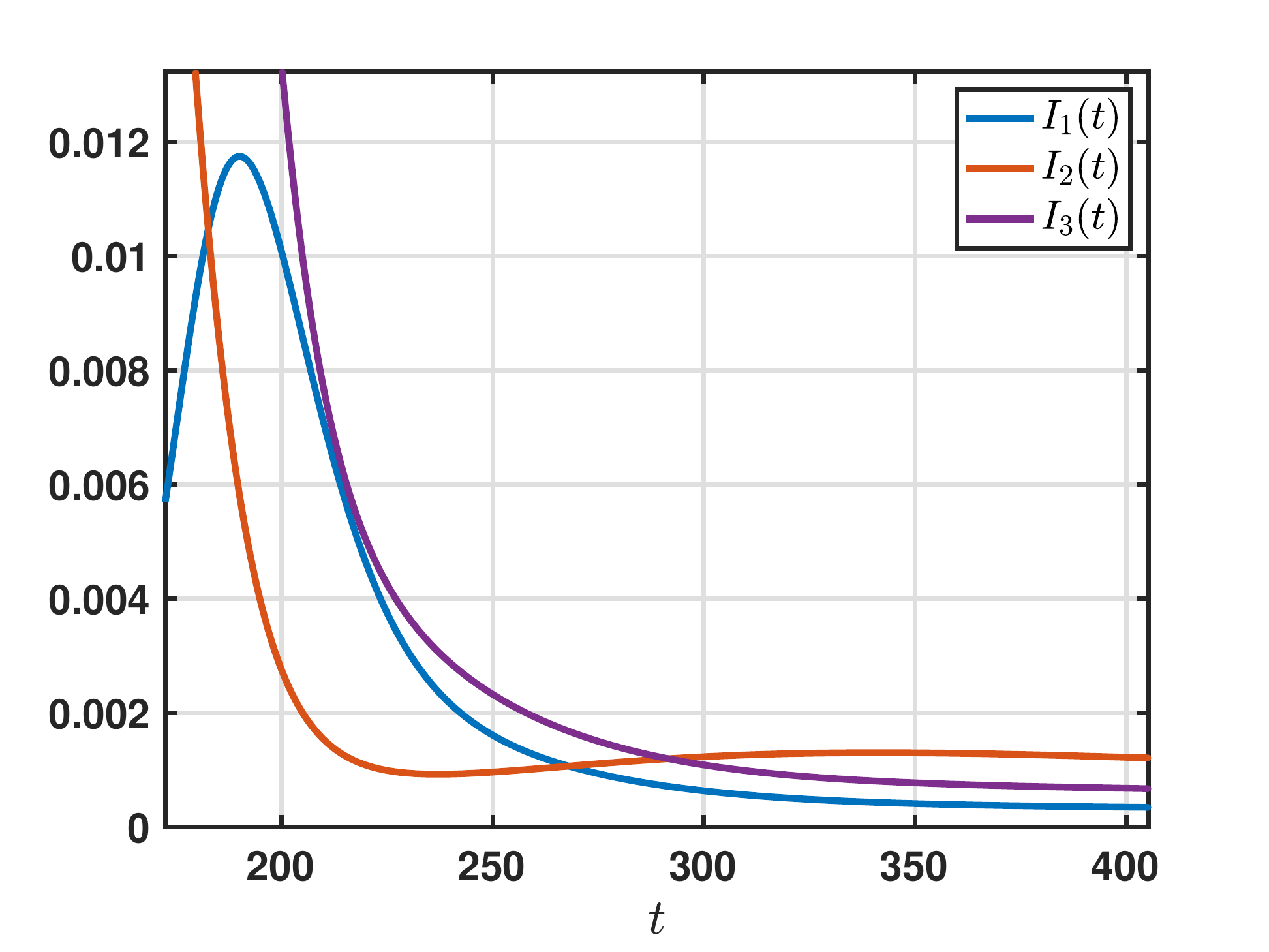}
\caption{Time plot of infected populations in the case of $3$ vertices and $3$ edges in a triangular configuration between times $[0,500]$ (left) and a zoom for times between $[150,400]$ (right). We observe a second wave of infection at the first vertex resulting from incoming infected individuals that have successively passed through the two other vertices. This second wave is also present at the second vertex with a slight increase of $I_2(t)$ after the second wave has reached the first vertex. Parameters values are set in the text.}
\label{fig:3c3rsec}
\end{figure}

In Figure~\ref{fig:3c3rsec}, we tested a different configuration. Upon labeling by $A$ the edge between vertices $v_1$ and $v_2$, $B$ the edge between vertices $v_2$ and $v_3$ and $C$ the edge between vertices $v_1$ and $v_3$, we have set the parameters to
\bqs
\alpha_A^1= \alpha_B^2=\alpha_C^3=0, \text{ and } \alpha_A^2=\alpha_B^3=\alpha_C^2= 1/10,
\eqs
while
\bqs
\lambda_A^1= \lambda_B^2=\lambda_C^3=1/20, \text{ and } \lambda_A^2=\lambda_B^3=\lambda_C^2= 0,
\eqs
and
\bqs
\nu_{A,C}^1=\nu_{B,A}^2=\nu_{C,B}^3=0, \text{ and } \nu_{C,A}^1=\nu_{A,B}^2=\nu_{B,C}^3=1/30.
\eqs
The length of each edge is fixed $\ell_e=\ell=1$ and $(\tau_v,\eta_v)=(1,1/7)$ at each vertex $v\in\V$. Finally, we have set different coefficients on each edge, namely $d_A=1$, $d_B=10^{-2}$ and $d_C=10^{-3}$. Initially, we assume that infected individuals are only present at vertex $v_1$ and each vertex has the same number of susceptible individuals fixed to $1/3$. Finally, for each $e\in\E$ we have set $u_e^0(x)=0$ on $\Omega_e$. We see in Figure~\ref{fig:3c3rsec} that such a configuration can generate a second wave of infection at the first and second vertices showing that transient dynamics can be complex with multiple bumps of infection.

\subsection{Case of 4 vertices and 3 edges}
\begin{center}
\begin{tikzpicture}
    \coordinate (A) at (-2,0);
    \coordinate (O) at (0,0);
    \coordinate (B) at (2,-2);
    \coordinate (C) at (2,2);
    
    \node[below] at (A) {$v_1$};
    \node[below] at (O) {$v_2$};
    \node[below] at (B) {$v_3$};
    \node[above] at (C) {$v_4$};

    % \node[label=left:$A$]  at (A) {$\bullet$};
    % \node[label={[shift={(-0.2,-0.9)}]$O$}] at (O) {$\bullet$};
    % \node[label=right:$B$] at (B)  {$\bullet$};
    % \node[label=right:$C$] at (C) {$\bullet$};

    \node  at (A) {$\bullet$};
    \node at (O) {$\bullet$};
    \node at (B)  {$\bullet$};
    \node at (C) {$\bullet$};

    \draw (A) -- (O);
    \draw (B) -- (O);
    \draw (C) -- (O);
\end{tikzpicture}
\end{center}

Next, we consider a star-shape graph with 4 vertices and 3 edges where one vertex is connected to the three others. In this configuration, we assume that our parameters may vary with respect to time, modeling locked down strategies for example \cite{GMS20,LMSW20}. More precisely, we will assume that there exists $T_{lock}$ and $\mu_{lock}$ such that the transmission rates can be written as
\bqs
\tau_v(t)=\left\{
\begin{array}{cl}
\tau\,, & t\in[0,T_{lock}],\\
\frac{\tau \exp(-\mu_{lock}(t-T_{lock}))+\tau_{lock}}{1+\exp(-\mu_{lock}(t-T_{lock}))}\,, & t>T_{lock},
\end{array}
\right.
\eqs
for each $v\in\V$ and for a given $0<\tau_{lock}<\tau$. We will assume that the four vertices are at equal distance such that $\ell_e=\ell$ for each $e\in\E$ and that the coefficient diffusion are equal on each edge, $d_e=d$, $e\in\E$. We further assume that at the central vertex $v_2$ exchanges are no longer allowed. That is, we impose that
\bqs
\alpha_e^2=\left\{
\begin{array}{cl}
\alpha\,, & t\in[0,T_{lock}],\\
\alpha \exp(-\mu_{lock}(t-T_{lock}))\,, & t>T_{lock},
\end{array}
\right. \quad e\in\E,
\eqs
while $\alpha_e^j=\alpha$ for $j\neq2$ and $e\in\E$, together with
\bqs
\lambda_e^2=\left\{
\begin{array}{cl}
\lambda\,, & t\in[0,T_{lock}],\\
\lambda \exp(-\mu_{lock}(t-T_{lock}))\,, & t>T_{lock},
\end{array}
\right. \quad e\in\E,
\eqs
while $\lambda_e^j=\lambda$ for $j\neq2$ and $e\in\E$, and also
\bqs
\nu_{e,e'}^2=\left\{
\begin{array}{cl}
\nu\,, & t\in[0,T_{lock}],\\
\nu \exp(-\mu_{lock}(t-T_{lock}))\,, & t>T_{lock},
\end{array}
\right. \quad (e,e')\in\E\times\E.
\eqs
Finally, we set $\eta_v=\eta$ for all $v\in\V$. Regarding the initial condition, we work with
\bqs
(S_1^0,I_1^0,S_2^0,I_2^0,S_3^0,I_3^0,S_4^0,I_4^0)=\left(S_0-I_0,I_0,S_0,0,S_0-\epsilon,0,S_0+\epsilon,0\right),
\eqs
for given $(S_0,I_0,\epsilon)$, while for each $e\in\E$ we have set $u_e^0(x)=0$ on $\Omega_e$.

\begin{figure}[!t]
\centering
\includegraphics[width=0.32\textwidth]{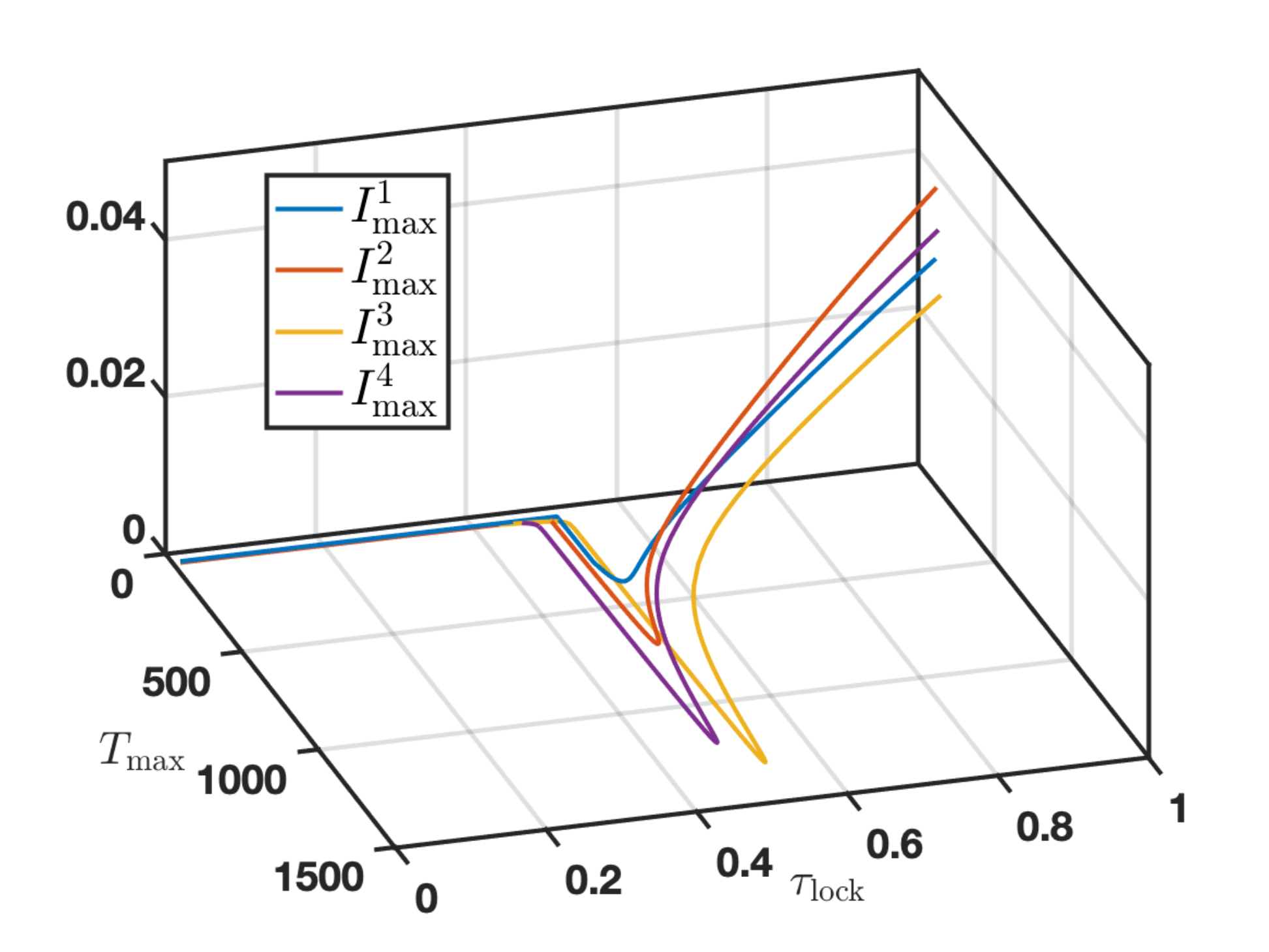}
\includegraphics[width=0.32\textwidth]{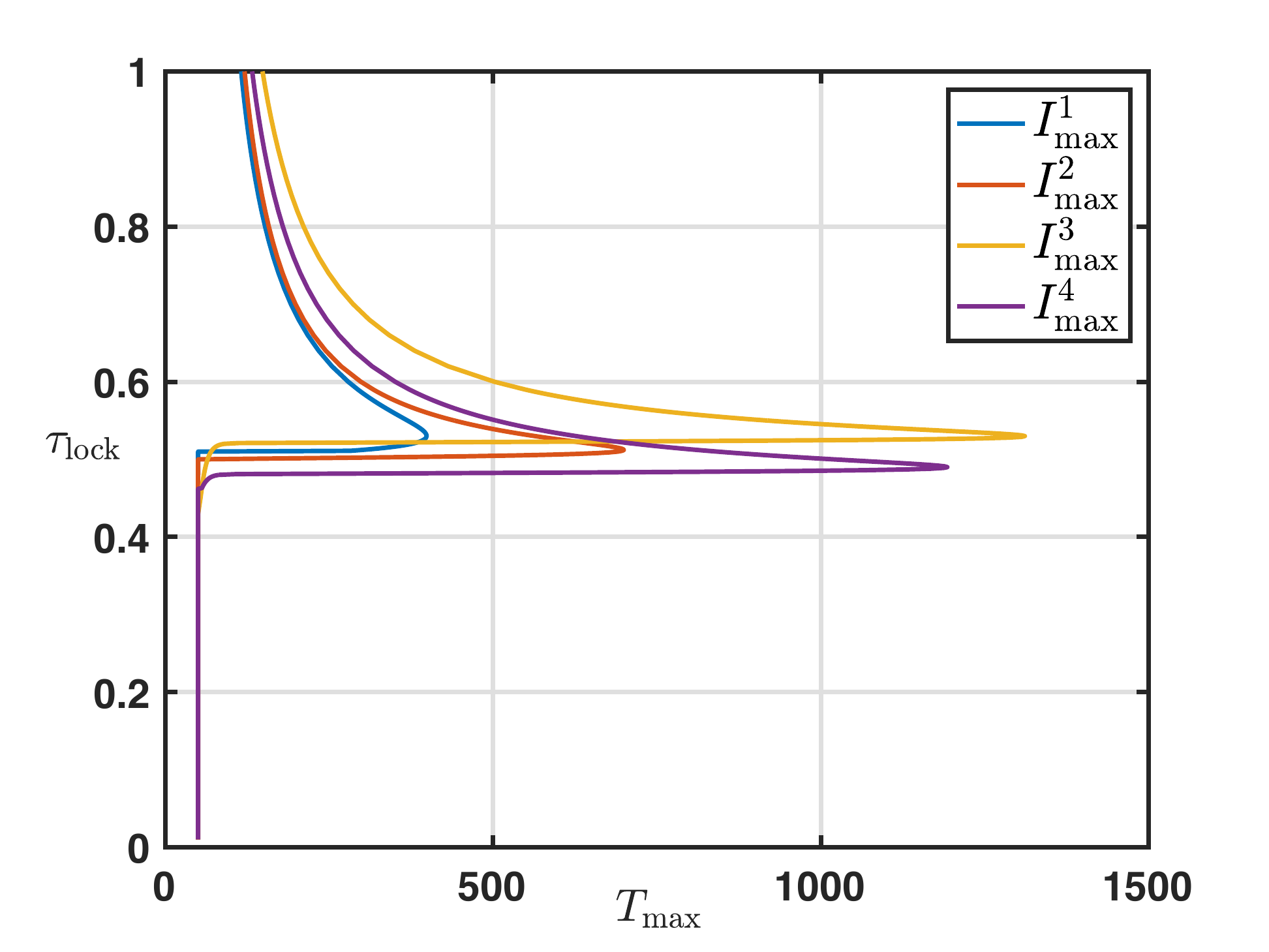}
\includegraphics[width=0.32\textwidth]{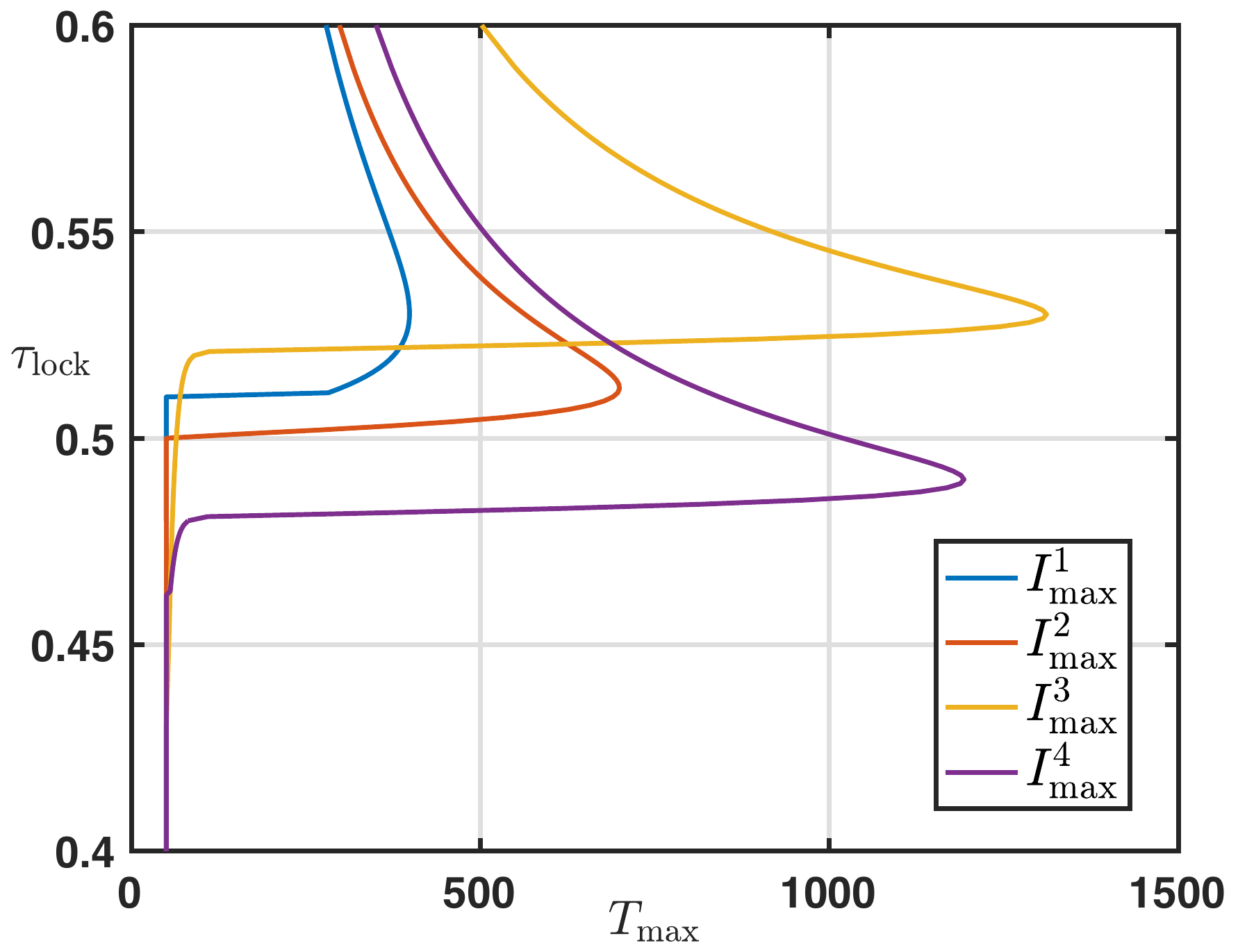}
\caption{Location of the time of maximal infection $T_{max}$ for each vertex together with the corresponding amplitude $I_{max}^j$  as a function of $\tau_{lock}$ (left) with its projection in the $(T_{max},\tau_{lock})$-plane (middle) and a zoom near the turning points (right). Other parameters are set to $\ell=1$, $d=0.1$, $\eta=1/8$, $(\alpha,\lambda,\nu)=(1/8,1/20,1/20)$, $T_{lock}=50$ and $\mu_{lock}=100$ with $(S_0,I_0,\epsilon)=(1/4,10^{-6},10^{-2})$.}
\label{fig:4c3rTau}
\end{figure}

\begin{figure}[!t]
\centering
\includegraphics[width=0.45\textwidth]{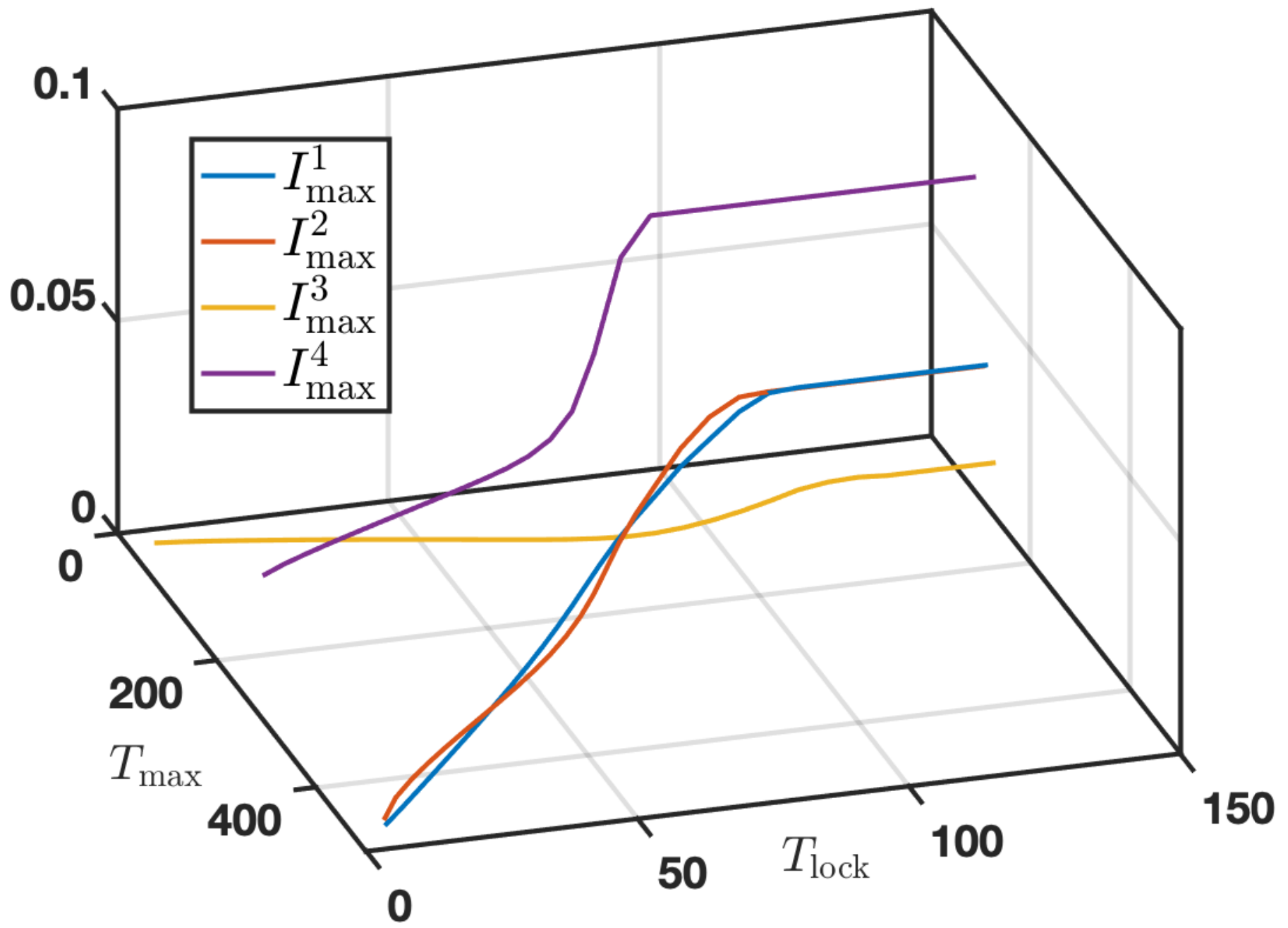}
\includegraphics[width=0.45\textwidth]{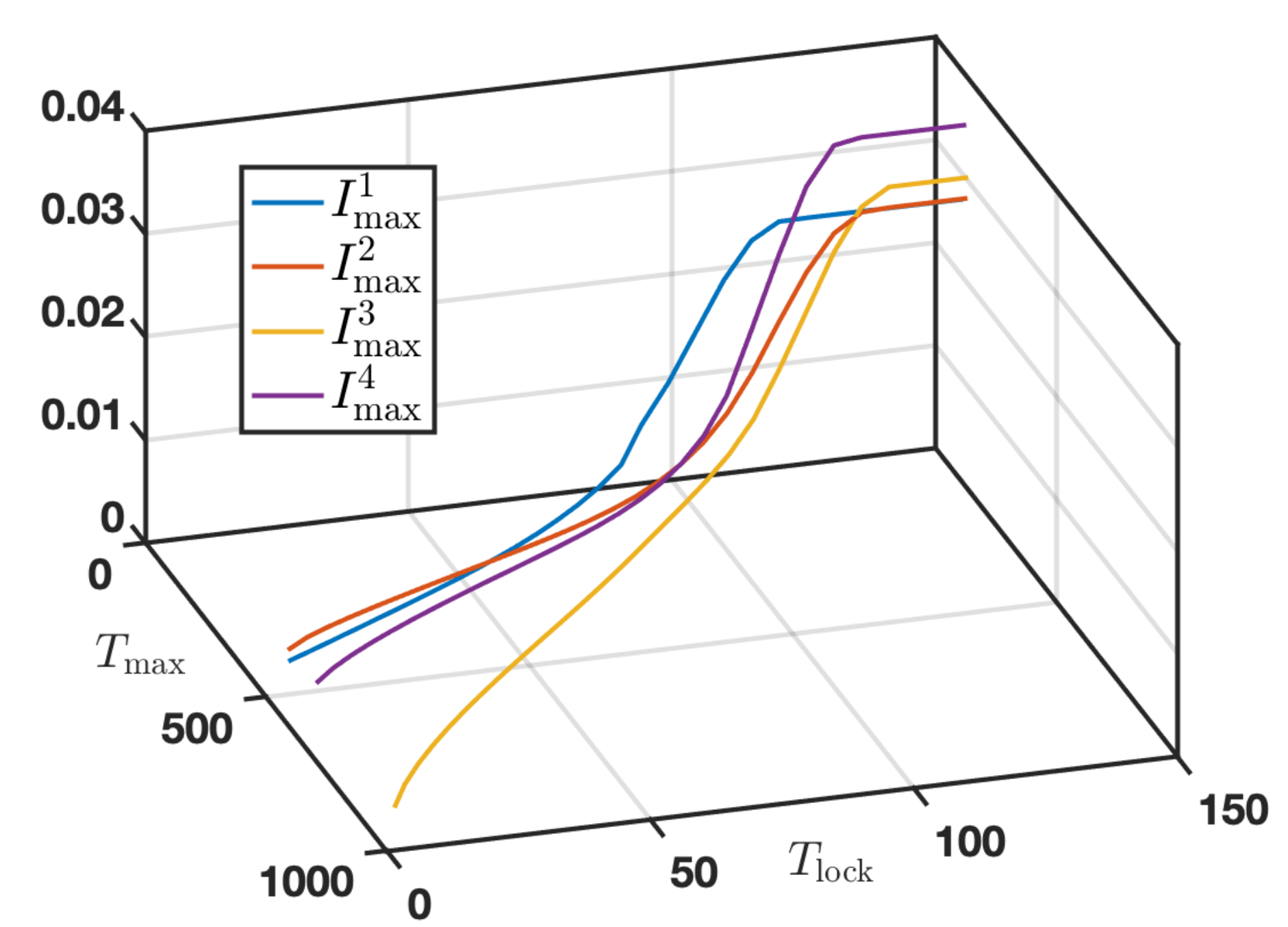}
\caption{Location of the time of maximal infection $T_{max}$ for each vertex together with the corresponding amplitude $I_{max}^j$  as a function of $T_{lock}$ for two configurations of initial susceptible populations at vertices $v_3$ and $v_4$, with $\epsilon=10^{-1}$ (left) and $\epsilon=10^{-2}$ (right). Other parameters are set to $\ell=1$, $d=0.1$, $\eta=1/8$, $(\alpha,\lambda,\nu)=(1/8,1/20,1/20)$, $T_{lock}=50$ and $\eta_{lock}=0.6$ with $(S_0,I_0)=(1/4,10^{-6})$.}
\label{fig:4c3rTauT}
\end{figure}

In Figure~\ref{fig:4c3rTau}, we report the location of the time of maximal infection $T_{max}$ for each vertex together with the corresponding amplitude $I_{max}^j=\max_{t\geq0}I_j(t)$  as a function of $\tau_{lock}$. We observe that below a critical value of $\tau_{lock}$, the time of maximal infection always occurs at $t=T_{max}=T_{lock}$ traducing the fact that the locked down strategy has no effect on the dynamics of the epidemic. At each vertex, we observe the same pattern: as $\tau_{lock}$ is decreased the corresponding $I_{max}^j$ is decreasing while $T_{max}$ is increasing up to some value of $\tau_{lock}$ where we observe a sudden turning point (see the right panel of Figure~\ref{fig:4c3rTau}). We observe that $\tau_{lock,v_k}^{tp}$, the value of the turning point, is well approximated (actually always bounded by below) by the value at which the effective reproduction number of each vertex is equal to $1$. Indeed we have $\mathscr{R}_{e,v_k}=1$ if and only if $\tau_{v_k}^{c}=\frac{\eta_{v_k}}{S^0_{v_k}}$, and we find
\bqs
\tau_{v_1}^{c}\simeq 0.5, \quad \tau_{v_2}^{c}= 0.5,\quad \tau_{v_3}^{c}\simeq 0.52, \quad \text{ and } \quad \tau_{v_4}^{c}\simeq 0.48,
\eqs
with our specific values of the initial condition, while we have computed
\bqs
\tau_{lock,v_1}^{tp}\simeq 0.53, \quad \tau_{lock,v_2}^{tp}\simeq 0.51,\quad \tau_{lock,v_3}^{tp}\simeq 0.53, \quad \text{ and } \quad \tau_{lock,v_4}^{tp}\simeq 0.49.
\eqs
We also point out that when $\tau_{lock}$ is below the turning point $\tau_{lock,v_k}^{tp}$, the corresponding value of $I_{max}^j$ is below $10^{-3}$. On the other hand, in Figure~\ref{fig:4c3rTauT}, we present similar results but this time $\tau_{lock}$ is fixed and $T_{lock}$ varies. Above some critical value of $T_{lock}$, $I_{max}^j$ saturates to a fixed value independent of $T_{lock}$ traducing the fact that the locked down strategy has no effect on the dynamics of the epidemic if it occurs to late in time. Depending on the initial configuration of susceptible populations at each vertex, we observe intricate nonlinear relationships on the location of the time of maximal infection $T_{max}$.

\subsection{Case of $N+1$ vertices and $N$ edges}

\begin{center}
\begin{tikzpicture}[scale=1.5]
    \coordinate (A1) at (0,0);
    \coordinate (A2) at (1,0);
    \coordinate (A3) at (2,0);
    \coordinate (A4) at (3,0);
    \coordinate (A5) at (5,0);
    \coordinate (A6) at (6,0);

    \node[below] at (A1) {$v_1$};
    \node[below] at (A2) {$v_2$};
    \node[below] at (A3) {$v_3$};
    \node[below] at (A4) {$v_4$};
    \node[below] at (A5) {$v_N$};
    \node[below] at (A6) {$v_{N+1}$};

    \node at (A1) {$\bullet$};
    \node at (A2) {$\bullet$};
    \node at (A3) {$\bullet$};
    \node at (A4) {$\bullet$};
    \node at (A5) {$\bullet$};
    \node at (A6) {$\bullet$};

    \draw (A1) -- (A2) -- (A3) -- (A4);
    \draw[dashed] (A4) -- (A5);
    \draw (A5) -- (A6);
\end{tikzpicture}
  \end{center}

\begin{figure}[!t]
\centering
\includegraphics[width=0.45\textwidth]{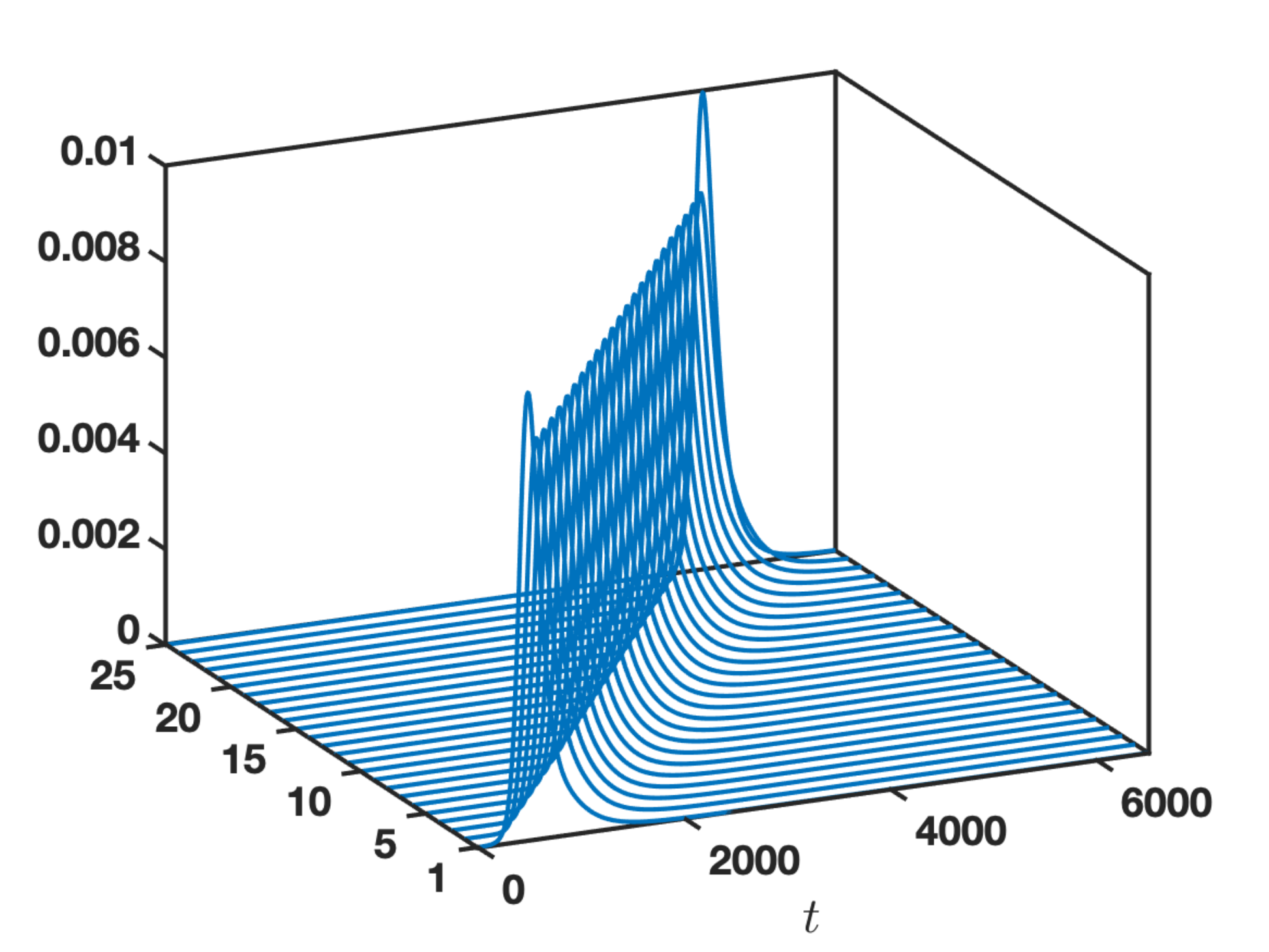}
\includegraphics[width=0.45\textwidth]{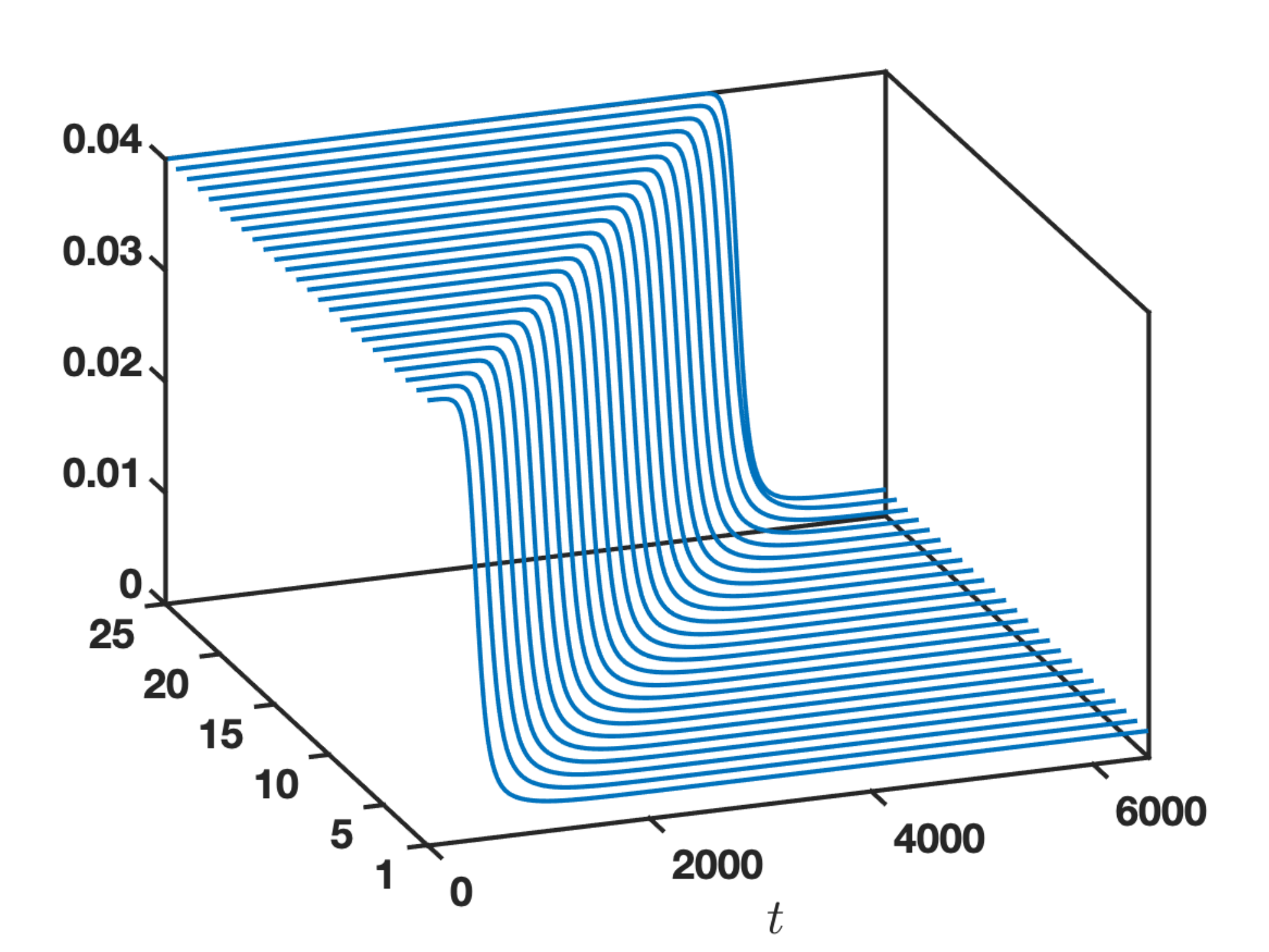}
\includegraphics[width=0.45\textwidth]{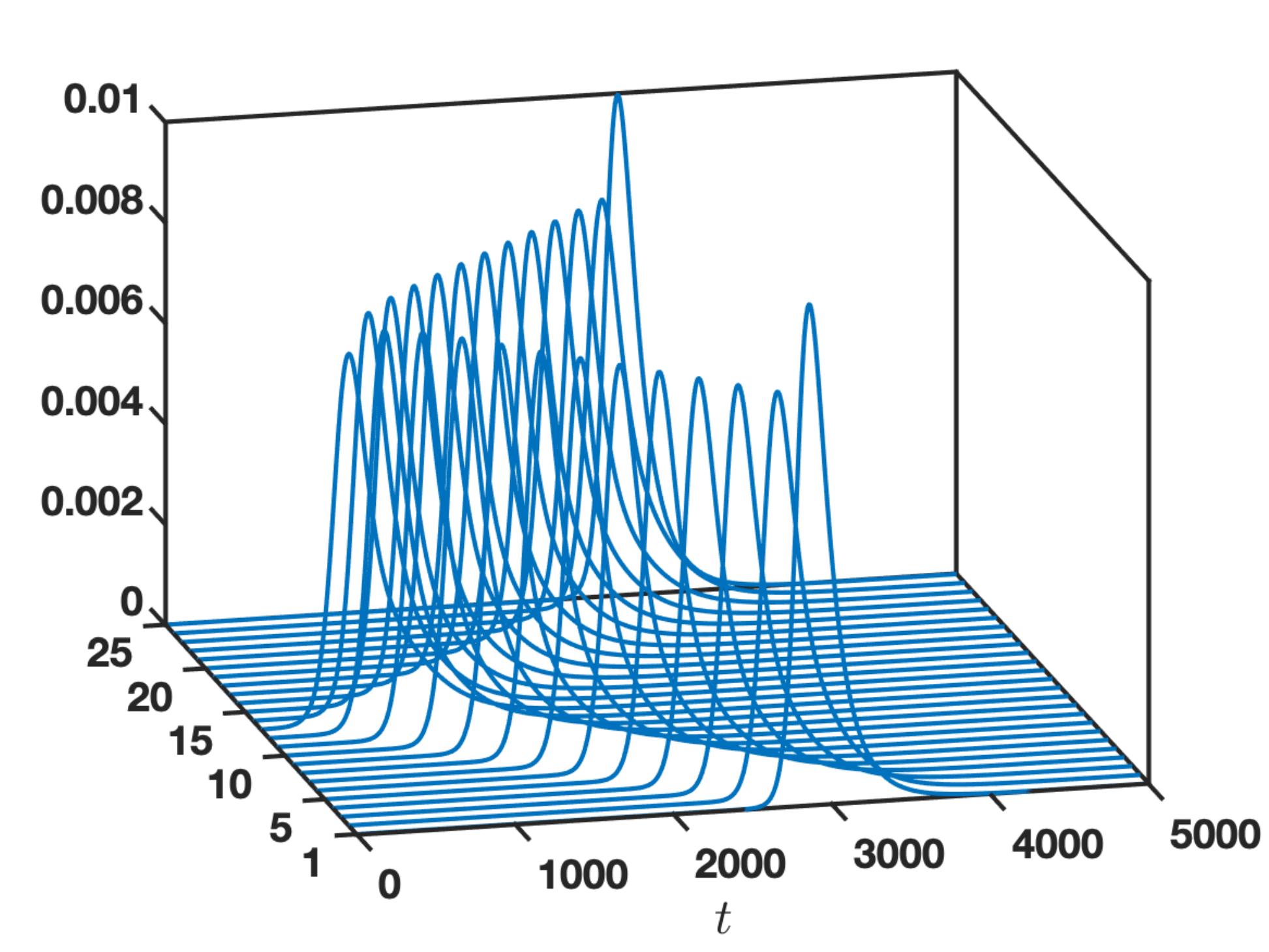}
\includegraphics[width=0.45\textwidth]{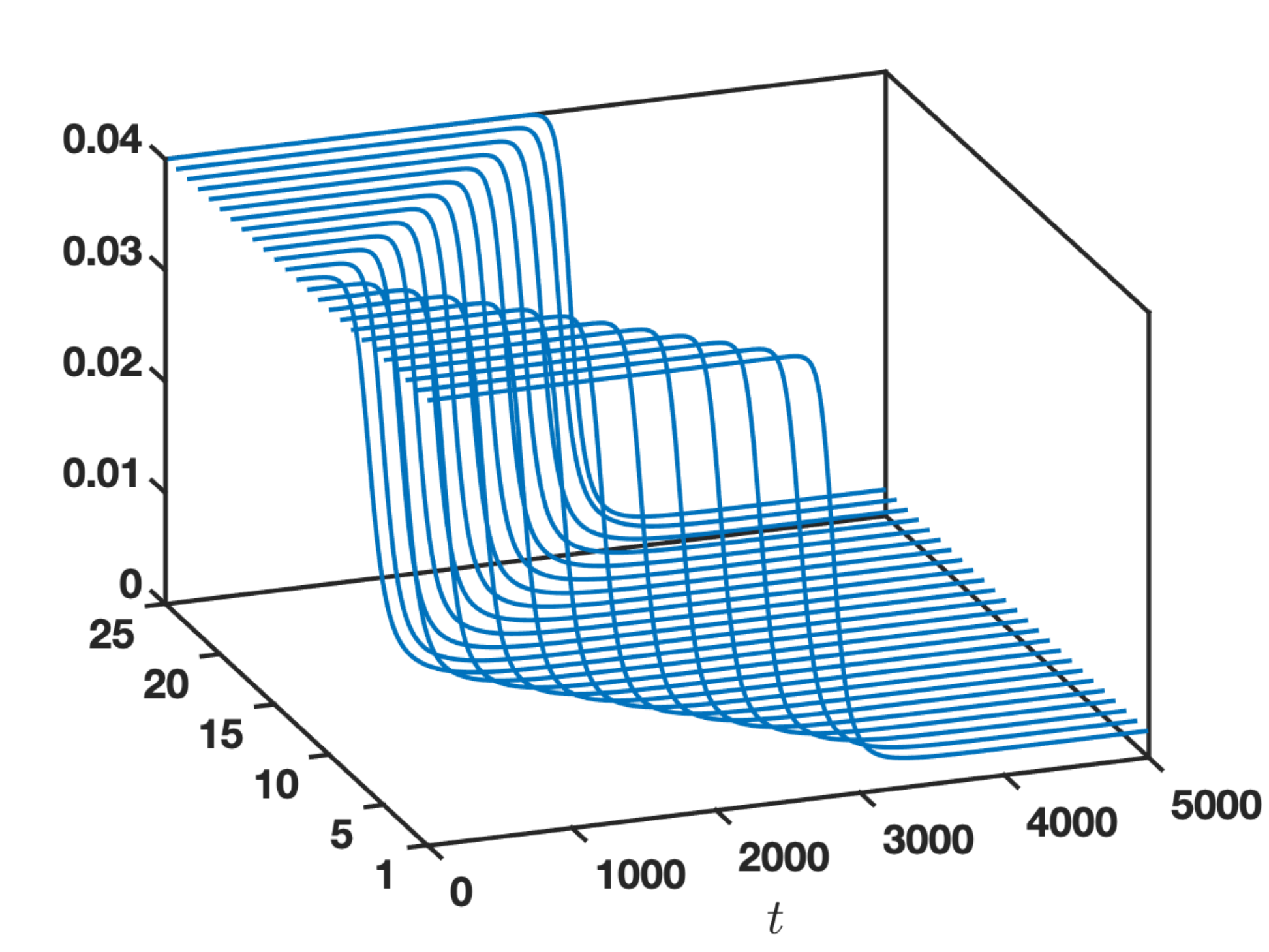}
\includegraphics[width=0.45\textwidth]{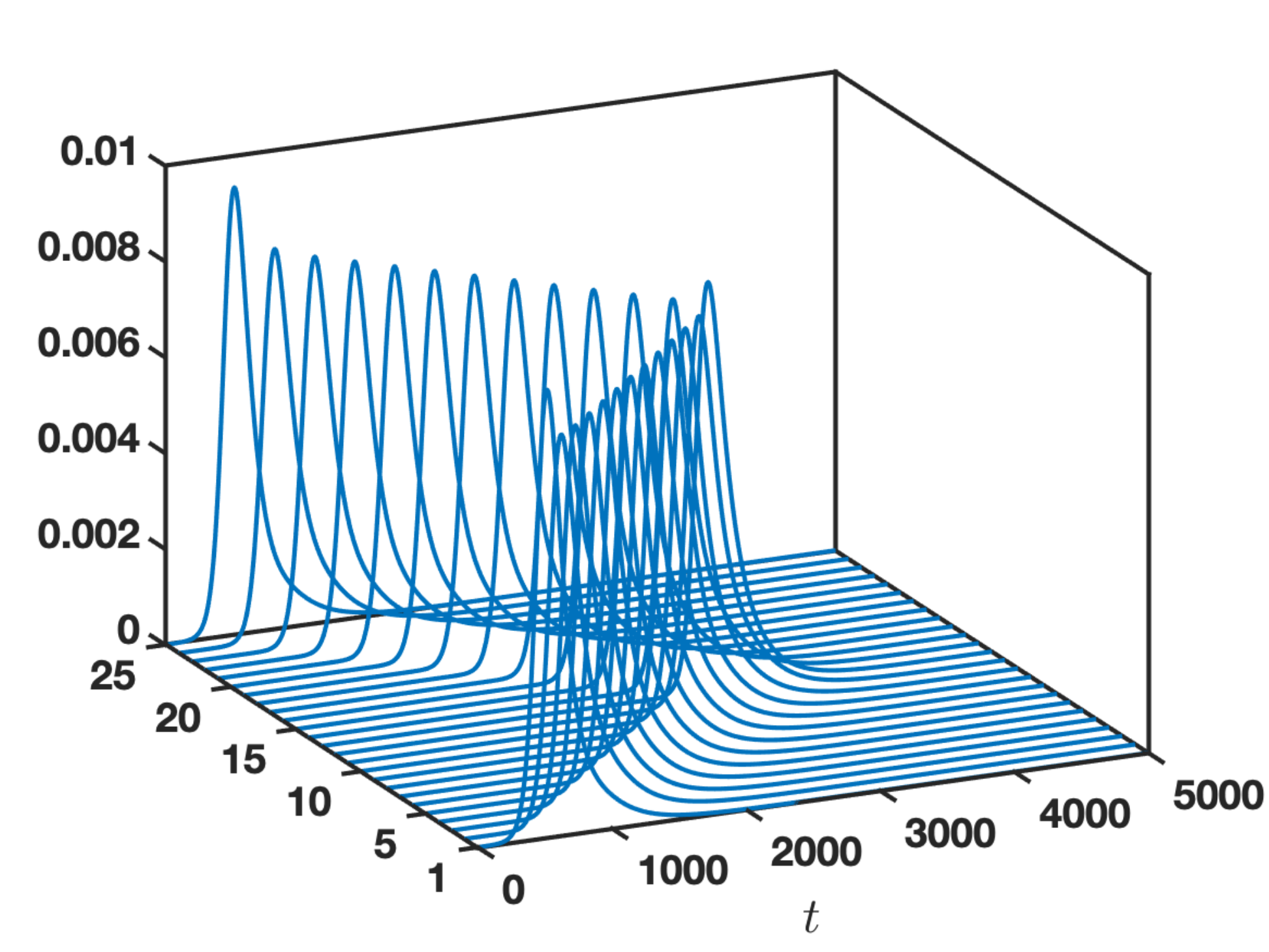}
\includegraphics[width=0.45\textwidth]{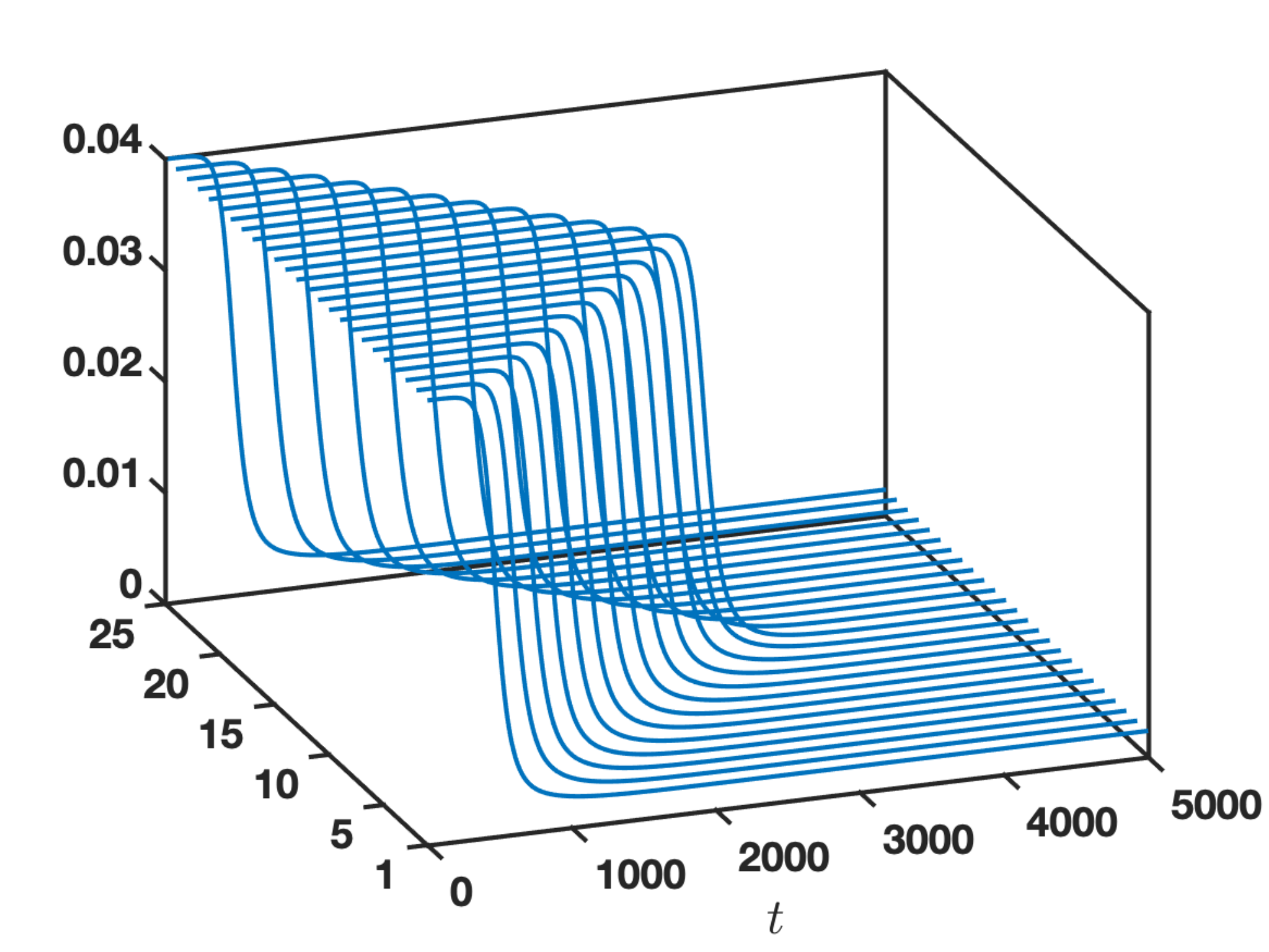}
\caption{Time evolution of the infected (left) population $I_j(t)$  and susceptible (right) populations $S_j(t)$ at each vertex for several different initial conditions and $N=24$. Top: infected individuals are initially present only at vertex. Middle:  infected individuals are initially present only at the middle vertex. Bottom: infected individuals are initially present only at the first and last vertices. We observe a traveling wave of infectious activity propagating though the vertices. Parameters were set to $\ell=1$, $d=10^{-3}$, $(\tau,\eta)=(1,1/75)$, and $(\alpha,\lambda,\nu)=(1/8,1/10,1/20)$, while the initial condition is $(S_0,I_0)=(1/25,10^{-6})$.}
\label{fig:TPvshape}
\end{figure}

\begin{figure}[!t]
\centering
\includegraphics[width=0.65\textwidth]{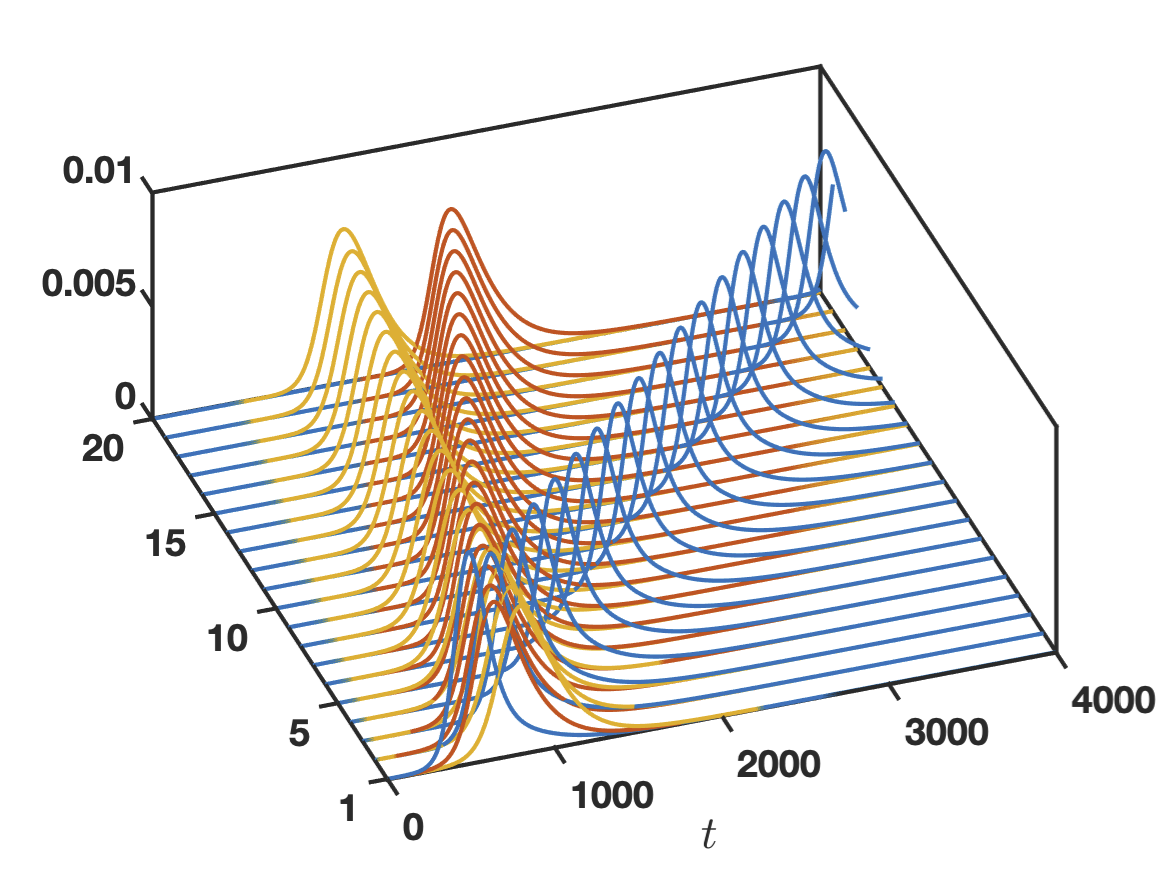}
\includegraphics[width=0.45\textwidth]{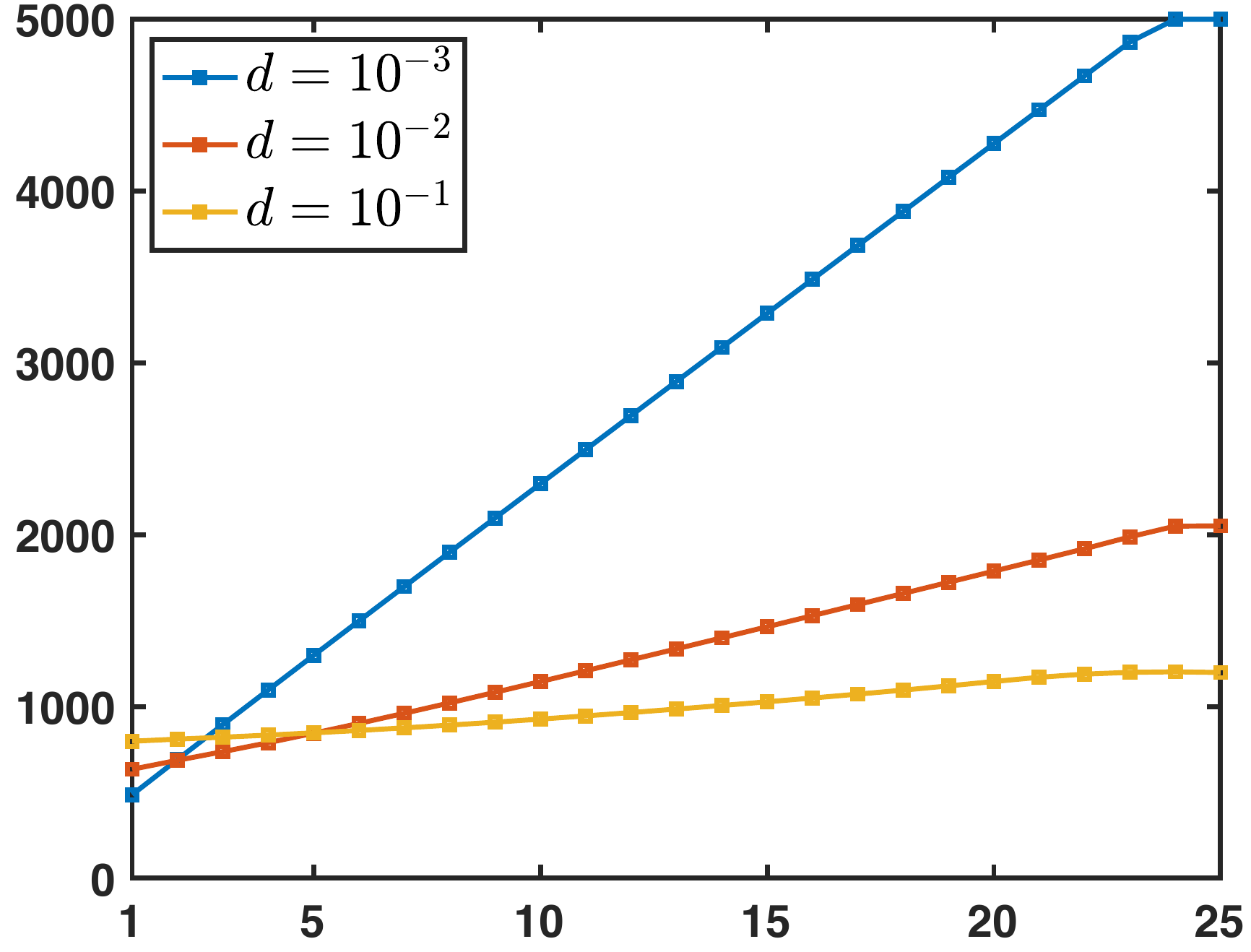}\hspace{0.2cm}
\includegraphics[width=0.45\textwidth]{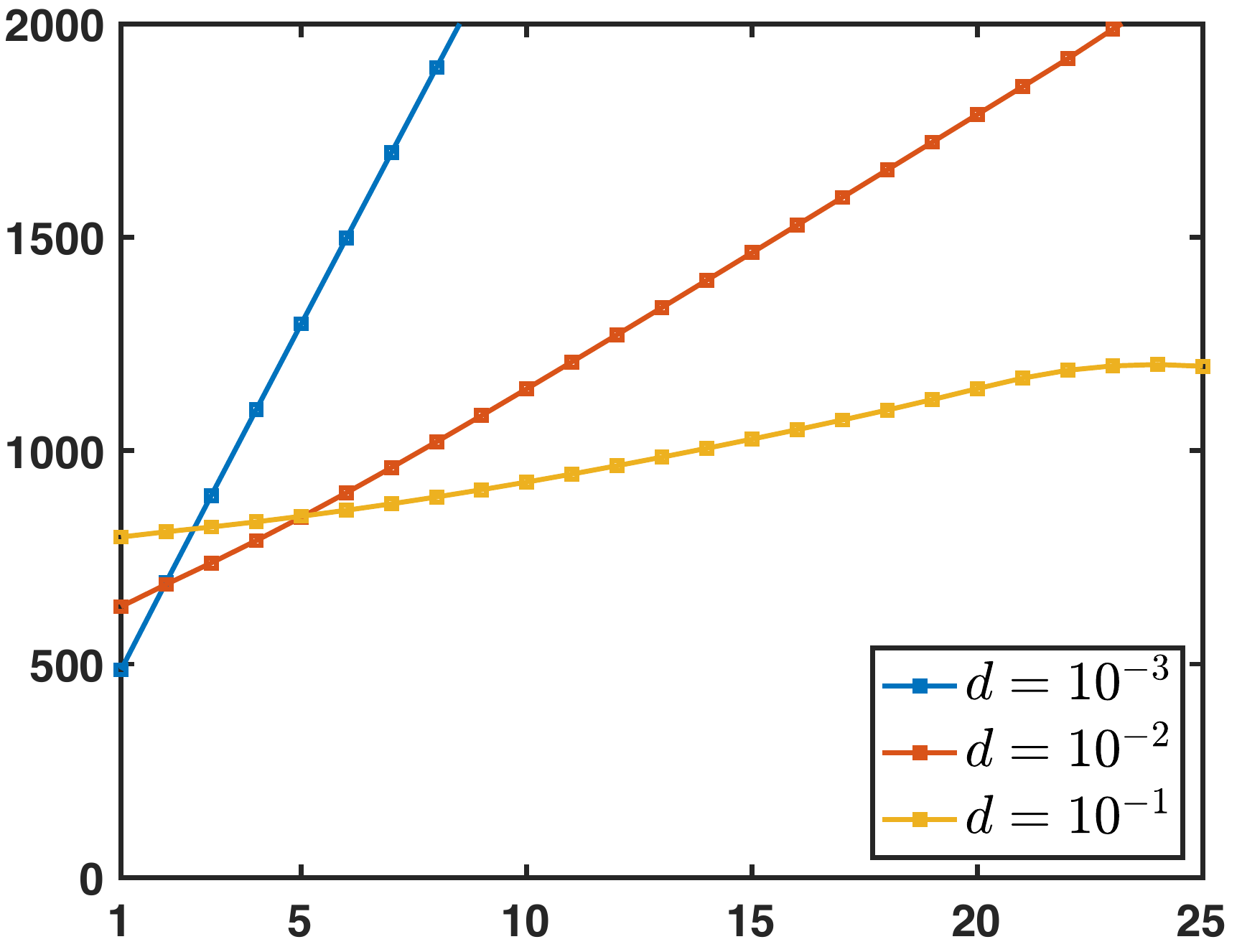}
\caption{Left: Location of the maxima of the infected population $I_j(t)$ at each vertex $j\in\left\{1,\cdots,25\right\}$ for several different diffusion coefficients $d\in\left\{ 10^{-3},10^{-2},10^{-1}\right\}$. Right: Zoom of the left figure for small times. We observe that for very small $d$ the location of the maxima is a long a line while it is curved for larger values of $d$. We also remark that if $I_{\mathrm{max},1}^d$ denotes the maximum as a function of $d$ at the first vertex, we have $I_{\mathrm{max},1}^{d_1}\leq I_{\mathrm{max},1}^{d_2}$ for $d_1\leq d_2$ while for larger vertices $j\geq 6$ we have the reverse ordering $I_{\mathrm{max},j}^{d_1}\geq I_{\mathrm{max},j}^{d_2}$ for $d_1\leq d_2$. Other parameters were set to $\ell=1$, $(\tau,\eta)=(1,1/75)$, and $(\alpha,\lambda,\nu)=(1/8,1/10,1/20)$, while the initial condition is $(S_0,I_0)=(1/25,10^{-6})$.}
\label{fig:TPmaxInfect}
\end{figure}

In our final example, we have considered a network of $N+1$ vertices and $N$ edges arranged in a lattice, in the sense that vertex $v_j$ is only connected to vertices $v_{j-1}$ and $v_{j+1}$ via two different edges. Figure~\ref{fig:TPvshape} shows the time evolution of the infected population $I_{v_j}(t)$  and susceptible populations $S_{v_j}(t)$ at each vertex for several different initial conditions when the length and diffusion coefficient of each edge are equal. In the first case (top panel), we assume that $I_{v_1}^0>0$ while $I_{v_j}^0=0$ for all other vertices, and observe a propagation of burst of activity among infected and susceptible populations. In the second case (middle panel), we assume that $I_{v_{\lfloor N/2\rfloor}}^0>0$ while $I_{v_j}^0=0$ for all other vertices, and we see the propagation of two bursts of activity among infected and susceptible populations going leftwards and rightwards. In the last case (bottom panel), we assume that $I_{v_1}^0=I_{v_{N+1}}^0>0$ while $I_{v_j}^0=0$ for all other vertices, and we note the propagation of two waves activity which collide at the middle vertex $v_{\lfloor N/2\rfloor}$. For very small values of the diffusion coefficient $d$, this burst of epidemic activity seems to travel coherently and forms a coherent traveling wave, as can be seen in Figure~\ref{fig:TPmaxInfect} where we represent the location of $\max_{t>0}I_{v_j}(t)$ at each vertex. Such a traveling wave of epidemic activity share similarities with traveling waves in excitable media such as the propagation of electrical activity along a nerve cell \cite{HH52,HS10} or calcium waves \cite{sneyd}. When $d=10^{-3}$, they are all aligned on the same line, where for smaller values $d\in\left\{10^{-1},10^{-2}\right\}$ the location is a nonlinear curve. We also demonstrate that larger diffusion coefficient leads to a faster propagation of epidemic burst across vertices. Finally, we also remark that if $I_{\mathrm{max},1}^d$ denotes the maximum as a function of $d$ at the first vertex, we have $I_{\mathrm{max},1}^{d_1}\leq I_{\mathrm{max},1}^{d_2}$ for $d_1\leq d_2$ while for larger vertices $j\geq 6$ we have the reverse ordering $I_{\mathrm{max},j}^{d_1}\geq I_{\mathrm{max},j}^{d_2}$ for $d_1\leq d_2$.

For the numerical simulations presented in Figures~\ref{fig:TPvshape}-\ref{fig:TPmaxInfect}, we have assumed full symmetry in the parameters that is
\begin{align*}
& (\ell_e,d_e)=(\ell,d), \quad e\in\E,\quad (\tau_v,\eta_v)=(\tau,\eta), \quad v\in\V,\\
& (\alpha_e^v,\lambda_e^v)=(\alpha,\lambda), \quad (e,v)\in\E\times\V, \quad \nu_{e,e'}^v=\nu, \quad (e,e',v)\in\E\times\E\times\V.
\end{align*}
 Regarding the initial condition on the edge, we have set $u_e^0(x)=0$ on $\Omega_e$ for each $e\in\E$. 
  
\section*{Acknowledgment}

This works was partially supported by Labex CIMI under grant agreement ANR-11-LABX-0040.

%\newpage


\begin{thebibliography}{10}

\bibitem{A77}
D. G. Aronson.
\newblock The asymptotic speed of propagation of a simple epidemic.
\newblock{\em Res. Notes Math.}, 14, pp. 1-23, 1977.

\bibitem{BB20}
F. Ball and T. Britton.
\newblock{ Epidemics on networks with preventive rewiring.}
\newblock{ ArXiv }, arXiv:2008.06375, 2020.

\bibitem{BRR20}
H. Berestycki, J.-M. Roquejoffre and L. Rossi. 
\newblock Propagation of epidemics along lines with fast diffusion. 
\newblock{\em Bull. Math. Biol.} to appear, 2020.

\bibitem{BG18}
L. Bonnasse-Gahot, H. Berestycki, M.-A. Depuiset, M. B. Gordon, S. Roch\'e, N. Rodriguez, and J.-P. Nadal.
\newblock Epidemiological modelling of the 2005 French riots: a spreading wave and the role of contagion.
\newblock{\em Scientific Reports}, 8, 2018.

\bibitem{BDL08}
T. Britton, M. Deijfen, M. Lindholm, and A. Nordvall Lageras.
\newblock Epidemics on random graphs with tunable clustering.  
\newblock{\em J. Appl. Prob.}, 45, 743-756, 2008.

\bibitem{CDC03}
Centers for Disease Control and Prevention . 
\newblock Severe acute respiratory syndrome -- Singapore, 2003. 
\newblock{\em Morbidity and mortality weekly report} 52.18, 405, 2003.

\bibitem{lambert}
R.M. Corless, G.H. Gonnet, D.E. Hare, D.J. Jeffrey and D.E. Knuth.
\newblock On the LambertW function.
\newblock{\em Advances in Computational mathematics,} 5(1), 329-359, 1996.

\bibitem{DIWB20}
JF David, SA Iyaniwura, MJ Ward and F Brauer.
\newblock A novel approach to modelling the spatial spread of airborne diseases: an epidemic model with indirect transmission. 
\newblock{\em Mathematical Biosciences and Engineering,} 17(4):3294, 2020.

\bibitem{D78}
O. Diekmann.
\newblock Thresholds and travelling waves for the geographical spread of infection.
\newblock{\em J. Math. Biol.}, 6, pp. 109-130, 1978.

\bibitem{DHM90}
O. Diekmann, J.A.P. Heesterbeek, J.A.J. Metz.
\newblock On the definition and the computation of the basic reproduction ratio $R0$ in models for infectious diseases in heterogeneous populations.
\newblock{\em J. Math. Biol.}, 28, p. 365, 1990.

\bibitem{GMS20}
Q. Griette, P. Magal and O. Seydi.
\newblock Unreported cases for Age Dependent COVID-19 Outbreak in Japan.
\newblock{\em Biology} 9, 132, 2020.

\bibitem{Heth00}
H.W. Hethcote.
\newblock The mathematics of infectious diseases.
\newblock{\em SIAM Rev. } 42 (4) 599--653, 2000.

\bibitem{HH52}
A.L. Hodgkin and A.F. Huxley. 
\newblock A quantitative description of membrane current and its application to conduction and excitation in nerve. 
\newblock{\em Journal of Physiology,}117, pages 500--544, 1952.

\bibitem{HS10}
H.J. Hupkes and B. Sandstede. 
\newblock Traveling pulse solutions for the discrete FitzHugh-Nagumo system. \newblock{\em SIAM J. Applied Dynamical Systems,} vol 9, no 3, pages 827--882, 2010.

\bibitem{KMK27}
W. O. Kermack, A. G. McKendrick.
\newblock  A contribution to the mathematical theory of epidemics.
\newblock{\em Proc. Roy. Sot. Ser. A}, 115, pp. 700-721, 1927.

\bibitem{LMSW20}
Z. Liu, P. Magal, O. Seydi, and G. Webb.
\newblock Predicting the cumulative number of cases for the COVID-19 epidemic in China from early data.
\newblock{\em Mathematical Biosciences and Engineering}, 17(4), 3040-3051, 2020.

\bibitem{MSB18}
P. Magal, O. Seydi and G. Webb.
\newblock Final size of a multi-group SIR epidemic model: Irreducible and non-irreducible modes of transmission.
\newblock{\em Mathematical Biosciences} 301, 59-67, 2018.

\bibitem{MSB16}
P. Magal, O. Seydi and G. Webb.
\newblock Final size of an epidemic for a two group SIR model.
\newblock{\em SIAM Journal on Applied Mathematics}, 76, 2042-2059, 2016.

\bibitem{MSS11}
S. Mandal, R.R. Sarkar and S. Sinha.
\newblock Mathematical models of malaria - a review.
\newblock{\em Malaria Journal,} 10:202, 1-19, 2011.

\bibitem{MW18}
P. Magal and G. Webb.
\newblock The parameter identification problem for SIR epidemic models: Identifying Unreported Cases.
\newblock{\em Journal of Mathematical Biology,} 77(6-7), 1629--1648, 2018.

\bibitem{NJE20}
New England Journal of Medicine.
\newblock{\em Letter to the Editor,} DOI: 10.1056/NEJMc2001468, January 30, 2020.

\bibitem{Reluga}
T. Reluga.
 \newblock A two-phase epidemic driven by diffusion.
 \newblock{\em Journal of theoretical biology,} 229.2: 249-261, 2004.

\bibitem{SE11}
M. Sekiguchi and I. Emiko. 
\newblock{ Dynamics of a discretized SIR epidemic model with pulse vaccination and time delay.}
\newblock{\em Journal of Computational and Applied Mathematics}, 236.6: 997-1008, 2011.

\bibitem{sneyd}
J. Sneyd. 
\newblock Tutorials in Mathematical Biosciences II. 
\newblock{\em Lecture Notes in Mathematics, chapter Mathematical Modeling of Calcium Dynamics and Signal Transduction,} Volume 187, Berlin Heidelberg, New York: Springer, 2005.

\bibitem{SB19}
K. Spricer and T. Britton. 
\newblock An epidemic model on a weighted network. 
\newblock{\em Network Science,} 7:556-580, 2019.

\bibitem{VanW02}
P. Van den Driessche and J. Watmough. 
\newblock Reproduction numbers and sub-threshold endemic equilibria for compartmental models of disease transmission.
\newblock{\em Mathematical biosciences,} 180.1-2: 29-48, 2002.

\end{thebibliography}
\end{document}